\documentclass[10pt]{article}
\usepackage[dvipdfmx]{graphicx}
\usepackage[dvipdfmx]{xcolor}
\usepackage{amsmath, amssymb}
\usepackage{amsthm} 
\usepackage[all]{xy}

\theoremstyle{plain} 
\newtheorem{theorem}{\indent\sc Theorem}[section] %

\newtheorem{corollary}[theorem]{\indent\sc Corollary}
\newtheorem{proposition}[theorem]{\indent\sc Proposition}

\theoremstyle{definition} 
\newtheorem{definition}[theorem]{\indent\sc Definition}
\newtheorem{remark}[theorem]{\indent\sc Remark}
\newtheorem{example}[theorem]{\indent\sc Example}

\def\even{\text{\rm even}}
\def\odd{\text{\rm odd}}

\def\L{\text{{\mathcal L}}}
\def\TT{T_M\oplus T^*_M}

\def\M{{\mathcal M}}
\def\Diff{\text{\rm Diff}}
\def\ol{\overline}
\def\part{\partial}

\def\tt{\scriptstyle  T\oplus T*}

\def\Sch{\text{\rm Sch}}

\def\L{{\mathcal L}}

\def\bgn{\begin}

\def\CL{\text{\rm CL}}

\def\J{{\mathcal J}}
\def\L{{\mathcal L}}

\def\1{{[1]}}
\def\2{{[2]}}
\def\3{{[3]}}
\def\({\left(}
\def\){\right)}
\def\s-circ{\,{\scriptstyle{\circ}}\,}
\def\<<{<\negthinspace \negthinspace<}
\def\Ad{\text{\rm Ad}}
\def\ad{\text{\rm ad}}

\def\even{\text{\rm even}}
\def\bgn{\begin}

\def\endaln{\end{align}}
\def\cou{\ss\text{\rm co}}

\def\<{<\negthinspace \negthinspace <}

\def\t{\theta}

\def\({\left(}
\def\){\right)}
\def\Im{\text{\rm Im}}
\def\Re{\text{\rm Re}}
\def\[{\big[\neg\big[}
\def\]{\big]\neg\big]}
\def\al{\al}

\def\M{{\mathcal M}}

\def\tr{\text{\rm tr}}
\def\a{\alpha}
\def\b{\beta}

\def\e{\varepsilon}
\def\gam{\gamma}
\def\Gam{\Gamma}

\def\del{\delta}
\def\lam{\lambda}
\def\ome{\omega}
\def\Ome{\Omega}
\def\sig{\sigma}

\def\A{{\mathcal A} }

\def\Diff{\text{\rm Diff}_0}

\def\R{\Bbb R}
\def\C{\Bbb C}

\def\M{\frak M}

\def\w{\wedge}
\def\({\left(}
\def\){\right)}

\def\neg{\negthinspace}
\def\h{\hat}

\def\wideh{\widehat}
\def\til{\tilde}
\def\wtil{\widetilde}

\def\ol{\overline}
\def\pa{\partial}

\def\ran{\rangle} 
\def\lan{\langle}

\def\ss{\scriptscriptstyle}
\def\trian{\triangle}

\def\arrow{\longrightarrow}

\def\bsh{\backslash}

\def\:{\, :\,}
\def\CL{\text{\rm CL}}
\def\TT{T\oplus T^*}
\def\complex{generalized complex }
\def\K\"ahler{generalized K\"ahler}
\def\vol{\text{\rm vol}}

\def\10{\displaystyle L^{10}}
\def\2{\displaystyle L^2}
\def\c0{\displaystyle C^0}

\def\10{\displaystyle L^{10}}
\def\2{\displaystyle L^2}
\def\del{\delta}
\def\del2{\displaystyle L^2_{0,\delta}}
\def\c0{\displaystyle C^0}

\def\del{\delta}
\def\cl{\text{\rm cl}}

\def\K{{\mathcal K}}
\def\M-A{\text{\rm Monge-Amp\`ere}}

\def\Ric{\text{\rm Ric}}
\def\M-A{\text{\rm Monge-Amp\`ere}}
\def\[{\big[\,}
\def\]{\,\big]}

\def\id{\text{\rm id}}

\def\ss{\scriptscriptstyle}
\def\even{\text{\rm even}}
\def\odd{\text{\rm odd}}

\def\L{\text{{\mathcal L}}}

\def\M{{\mathcal M}}
\def\Diff{\text{\rm Diff}}
\def\ol{\overline}
\def\part{\partial}

\def\ham{{\mathcal ham}}

\def\L{{\cal L}}
\def\R{\mathbb R}
\def\B{{\cal B}}
\def\frak g{\mathfrak{g}}
\def\frak h{\mathfrak{h}}
\def\C{\mathbb C}

\def\frak{\mathfrak}
\def\ham{\mathfrak ham}
\def\ol{\overline}
\def\T{\mathbb T}
\def\TT{\T M}

\makeatother
\title{Generalized scalar curvature and modified moment map in generalized K\"ahler geometry  \\
} 

\author{Ryushi Goto}

\date{} 
\begin{document}
\maketitle

\begin{abstract}
We introduce a notion of generalized scalar curvature for generalized K\"ahler manifolds using the pure spinor formalism. For an arbitrary compact generalized K\"ahler manifold, we develop an extended action of generalized Hamiltonians involving an abelian Lie algebra. This leads naturally to the notion of a modified moment map, which satisfies the moment-map condition only modulo the infinitesimal action of this abelian Lie algebra. We prove that this modified moment map is given by the generalized scalar curvature. Thus our result extends the theorem of Fujiki and Donaldson, which realizes the scalar curvature as a moment map in ordinary K\"ahler geometry.
We also study explicit examples. Although a compact connected Lie group admits a K\"ahler structure only in the torus case, every connected compact even-dimensional Lie group admits generalized K\"ahler structures with constant generalized scalar curvature. In particular, we explicitly construct such structures on the standard Hopf surface.
\end{abstract}

\numberwithin{equation}{section}
\tableofcontents

\noindent
\section{Introduction}
Let $M$ be a compact manifold of real dimension $2n$ equipped with a symplectic form $\ome$. 
We denote by $C_\ome$ the space of $\ome$-compatible almost complex structures on $(M, \ome)$. 
Then $C_\ome$ is an infinite dimensional homogeneous K\"ahler manifold, which can be identified with the space of global sections of a fibre bundle whose fibre is
the Siegel space Sp$(2n , \R)/$U$(n)$.  The group of Hamiltonian diffeomorphisms acts on $C_\ome$ preserving its K\"ahler structure. 
Fujiki and Donaldson showed that the scalar curvature arises as a moment map for this action of Hamiltonian diffeomorphisms on $C_\ome$,  
which is now known as the moment map framework in K\"ahler geometry \cite{Fu_1990, Do_1997}. \\
\indent
The main purpose of this paper is to extend the moment map framework to arbitrary generalized K\"ahler manifolds, more precisely, to
$H$-twisted generalized K\"ahler manifolds,  where $H$ denotes a real $d$-closed $3$-form on $M$.
A generalized K\"ahler structure on a manifold can be described as a triple $(g, I,J)$ consisting of a Riemannian metric $g$ compatible 
with two complex structures $I$ and $J$ satisfying  certain integrability conditions \cite{Gua_2003}. This notion has its origins in  nonlinear sigma models in mathematical physics.
However, a generalized K\"ahler structure also has a natural description in the language of nondegenerate, pure spinors \cite{Hi_2003}. Namely,  a generalized K\"ahler structure is defined by a pair $(\phi, \psi)$ of locally defined nondegenerate, pure spinors which induces 
a generalized K\"ahler structure $(\J_\phi, \J_\psi)$ \cite{Gua_2003}.
Using the pure spinor formalism, we introduce an invariant function 
$$S(\J_\phi, \J_\psi,\vol_M)$$ associated with a generalized K\"ahler structure $(\J_\phi, \J_\psi)$ together with a volume form $\vol_M$;  see Definition \ref{definition of scalar curvature} in Section 4. This function is called the generalized scalar curvature in this paper.
The generalized scalar curvature has the following properties:
\bgn{enumerate}
\item If $(\J_\phi, \J_\psi,\vol_M)$ comes from an ordinary K\"ahler structure, then the generalized scalar curvature $S(\J_\phi, \J_\psi,\vol_M)$ coincides with the ordinary scalar curvature.
If $(\J_\phi, \J_\psi)$ is a generalized K\"ahler structure of symplectic type, then $S(\J_\phi, \J_\psi,\vol_M)$ coincides with the generalized scalar curvature introduced in the previous paper \cite{Goto_2016}, \cite{Goto_2020}.
\item The generalized scalar curvature $S(\J_\phi, \J_\psi,\vol_M)$ is equivariant under the action of an extension 
$\wtil\Diff_0(M,H)$ of the diffeomorphism group of $M$ by $b$-field transformations.
\item The generalized scalar curvature is symmetric under the interchange of $\J_\phi$ and $\J_\psi$, i.e., 
$$
S(\J_\phi, \J_\psi,\vol_M)=S(\J_\psi, \J_\phi,\vol_M).
$$
\end{enumerate}
As in the moment map framework in K\"ahler geometry,  we fix an integrable generalized complex structure $\J_\psi$ and a volume form $\vol_M$ satisfying
$$\vol_M =i^{-n}\lan \psi_\a, \,\,\ol\psi_\a\ran_s.$$
We denote by $\B_{\J_\psi}(M)$ the space of $\J_\psi$-compatible almost generalized complex structures.
This is also an infinite-dimensional homogeneous K\"ahler manifold, and it can be identified with the space of global sections of a fibre bundle whose fibre is $\mathrm {U}(n,n)/(\mathrm{U}(n)\times \mathrm{U}(n))$.
The generalized complex structure $\J_\psi$ determines the Lie algebra of  generalized Hamiltonians $\ham(M, \J_\psi)$, which is given by 
$$\ham(M,\J_\psi):=C^\infty_\psi(M,\R)=\{\, f\in C^\infty(M,\R)\, |\, \pi_T(\zeta_\a)f=0, \, \}.$$
An abelian extension $\wtil\ham(M, \J_\psi)$ of $\ham(M, \J_\psi)$
acts on $\B_{\J_\psi}(M)$ by the extended action, preserving its K\"ahler structure (see Section \ref{GSCasmoment}  for more details).
In contrast to the ordinary K\"ahler case, the natural action in generalized K\"ahler geometry involves not only diffeomorphisms but also $b$-field transformations. 
Moreover, 
the infinitesimal action contains an additional abelian component, which may obstruct the existence of a moment map in the usual sense.
This leads naturally to the notion of a modified moment map; see Section \ref{Def and results of modified moment map} and Section \ref{GSCasmoment} for more details. 
Then we have
 \medskip \\
\noindent{\sc Theorem A }(see Theorem \ref{existence of moment map}).
Let $\J_\psi$ be a $H$-twisted generalized complex structure on a compact manifold $M$ equipped with the fixed volume form $\vol_M$.
Then there exists a modified moment map 
$$\mu: {\mathcal B}_{\J_\psi}(M)\to C^\infty_\psi(M,\R)^*$$ 
for the action of the Lie algebra 
$\wtil\ham(M, \J_\psi)$ such that
 $\mu(\J)$ is given by the generalized scalar curvature $S(\J, \J_\psi,\vol_M)$. More precisely,
with respect to the natural pairing, one has
$$
\langle \mu(\J),f\rangle
=
\int_M f\,S(\J,\J_\psi,\vol_M)\,\vol_M
$$
for all $f\in C^\infty_\psi(M,\mathbb R)$.
 
\smallskip
Since an earlier version of this paper [arXiv:2105.13654] was posted on arXiv, Apostolov, Streets, and Ustinovskiy have given an explicit formula for the generalized scalar curvature \cite{ASU_2026}. They also pointed out that this quantity had already appeared in several different contexts, including the string effective action \cite{CSW_2011}, Bismut's work on the index theorem in complex geometry \cite{Bi_1989},
 and generalized Ricci flow \cite{MGFJS_2020}. 
From this perspective, the derivation of the generalized scalar curvature from the spin geometry of generalized K\"ahler structures reveals a rather unexpected link between the moment map framework and these other appearances of scalar curvature in generalized geometry.

In the first version of this paper, we established a variational formula for the generalized scalar curvature under generalized Hamiltonian deformations. Apostolov, Streets, and Ustinovskiy subsequently introduced the notion of an adapted volume and showed, using this variational formula, that the generalized scalar curvature defines an ordinary moment map on the space of generalized K\"ahler structures admitting a fixed adapted volume.
The present framework provides a geometric interpretation of their adapted-volume condition. The infinitesimal generalized Hamiltonian action naturally extends by an additional abelian part. We show that the locus on which the prescribed volume form is adapted is precisely the fixed-point locus of this abelian action. Thus, on the full space of compatible almost generalized complex structures, the generalized scalar curvature is naturally viewed as a modified moment map, whereas its restriction to the adapted-volume locus becomes the ordinary moment map considered by Apostolov, Streets, and Ustinovskiy.


Even in the case of ordinary K\"ahler manifolds, replacing the standard volume form determined by the K\"ahler form with a general volume form leads to a useful variant of the generalized scalar curvature. Inoue's $\mu$-scalar curvature K\"ahler metrics, which extend the notion of K\"ahler--Ricci solitons on Fano manifolds \cite{Ino_2022} to arbitrary K\"ahler classes, as well as Lahdili's weighted scalar curvature \cite{La_2019}, may be naturally understood from the viewpoint of the generalized scalar curvature. In Subsection 2.2 of Section \ref{Def and results of modified moment map}, we discuss this important point in a general framework.
\\
\indent
A non-abelian connected compact Lie group admits no K\"ahler structure; however, every compact even-dimensional Lie group admits a remarkable bi-Hermitian structure twisted by the Cartan $3$-form, which gives rise to a twisted generalized K\"ahler structure \cite{Gua_2003}. 
Alekseev, Bursztyn, and Meinrenken provided an interesting description of Dirac structures on a Lie group $G$ in terms of 
the double of its Lie algebra $\mathfrak{g}\oplus\ol{\mathfrak g}$ equipped with an Ad-invariant metric $B$.
We apply their description to construct a family of twisted generalized K\"ahler structures on  a compact Lie group $G$, by using 
the action of the real Pin group of the double of the real Cartan subalgebra $\mathfrak h\oplus \ol{\mathfrak h}$. 
Using the description of pure spinors on $G$ in terms of pure spinors of $\mathfrak g\oplus\ol{\mathfrak g}$, we compute the generalized scalar curvature of the twisted generalized K\"ahler structures 
$(\J_\phi, \J_\psi)$ on $G$. 
We denote by $\Xi$ the Cartan-Killing $3$-form on $\mathfrak g$, and by $\hat B$ the Hermitian form induced by $B$. 
The contraction $\sqrt{-1}\Lambda_{\hat B}\Xi\in \mathfrak h$ is denoted by $P$.
We show that the generalized scalar curvature of $(\J_\phi, \J_\psi)$ on $G$ is the constant $2\|P\|^2$, 
where $\|P\|$ is the norm of $P$ with respect to $\hat B$.

This paper is organized as follows.
In Section \ref{Def and results of modified moment map}, we introduce the notion of modified moment map.
In Section \ref{GCS GK}, we fix our notation and recall preliminary results on generalized complex structures and generalized K\"ahler structures used throughout the paper.
In Section \ref{Generalized Scalar curvature}, we define the generalized scalar curvature of  generalized K\"ahler structures. 
In Section \ref{Generalized Hamiltonian diffeomorphisms}, the notion of generalized Hamiltonians is introduced. 
In Section \ref{GSCasmoment}, 
we develop the moment map framework and prove Theorem \ref{existence of moment map}. 
In Section \ref{Generalized scalar curvature of generalized Kahler structures on the standard Hopf surfaces}, we explicitly construct generalized K\"ahler structures on the standard Hopf surfaces with constant generalized scalar curvature. 
In Section \ref{Review}, we give a brief review of the results of
Alekseev, Bursztyn, and Meinrenken. 
Finally in Section \ref{GK on compact Lie group}, we construct a family of twisted generalized K\"ahler structures with constant generalized scalar curvature on 
a compact Lie group $G.$ 

\section{Modified moment maps}\label{Def and results of modified moment map}
\subsection{Definition of modified moment maps}
Let $X$ be a symplectic manifold equipped with a symplectic form $\ome$, and suppose that a Lie algebra $\mathfrak g$ acts on $X$ by symplectic vector fields.
We assume that $\mathfrak g$ is an extension of a Lie algebra $\mathfrak g_0$ by an abelian Lie algebra $\mathfrak a$, that is, there is a short exact sequence of Lie algebras
$$
0\longrightarrow \mathfrak a\longrightarrow  \mathfrak g\overset\pi{\longrightarrow}  \mathfrak g_0{\to} 0.
$$
We assume that there is a Lie algebra homomorphism $s: \mathfrak g_0\longrightarrow \mathfrak g$ such that 
$\pi\circ s=\id_{\mathfrak g_0}$. 
Then $\mathfrak g_0$ is regarded as a Lie subalgebra of $\mathfrak g$ by using $s$.

The modified moment map introduced in this section satisfies the usual moment map condition in a weak sense.  The term ``modified'' indicates that the usual moment map condition is
changed by a correction term coming from the abelian part 
$\mathfrak a$.

\begin{definition}[modified moment map]\label{modified moment map}
Let  $(X,\omega)$ be a symplectic manifold equipped with a symplectic action of a Lie algebra $\mathfrak g$ as before.
Then both $\mathfrak g_0$ and $\mathfrak a$ act on $(X,\ome)$ by symplectic vector fields.
We denote by $\h \xi$ the vector field on $X$ given by $\xi\in \mathfrak g$.
Let $\rho:\mathfrak g_0\to\mathfrak g$ be a linear map which gives the splitting of the 
sequence $0\to \mathfrak a\to \mathfrak g\overset{\pi}{\to} \mathfrak g_0\to 0,$ that is,
$$ \pi\circ\rho=\id_{\mathfrak g_0} .$$
A  map $\mu: X \to \mathfrak g_0^*$ 
is {\it a modified moment map }associated with the splitting $\rho$ if $\mu$ satisfies the following conditions: 
\begin{itemize}
\item[(1)]
$\mu$ is  infinitesimally $\mathfrak g_0$-equivariant, that is, 
$$\lan d\mu, \xi_1\ran(\h\xi_2) =\lan \mu, [\xi_1,\xi_2]\ran, \quad  
\qquad \forall \xi_1, \xi_2\in \mathfrak g_0,$$
where we write $\h \xi:=\widehat{s(\xi)}$ for $\mathfrak g_0$.
\item[(2)] $\mu$ satisfies 
$$d\lan \mu, \xi \ran=  i_{\widehat{\rho(\xi)}}\ome \quad (\xi \in \mathfrak g_0),$$
where $\widehat{\rho(\xi)}$ denotes the vector field given by $\rho(\xi)\in \mathfrak g$ for $\xi\in \mathfrak g_0$.
\end{itemize}
\end{definition}
If the splitting $\rho$ is a homomorphism of Lie algebras, then $\rho(\mathfrak g_0)$ is a Lie subalgebra of $\mathfrak g$ by using $\rho$, and $\mu$ is a moment map for the action of $\mathfrak g_0$ on $(X,\omega)$.
In the case of a modified moment map, however, $\rho$ is not a homomorphism in general, and hence $\mu$ is not a moment map in the usual sense. 
The difference from the usual moment map is entirely accounted for by the abelian action of $\mathfrak a$, that is, 
$$
[\rho(\xi_1), \rho(\xi_2)]-\rho([\xi_1, \xi_2])\in \mathfrak a, 
$$
where $\xi_1, \xi_2\in \frak g_0$.
In fact, the splitting map $\rho$ is written as 
$$
\rho(\xi) =s(\xi)+R(\xi), \qquad R(\xi)\in \mathfrak a.
$$
Then the condition (2) of Definition \ref{modified moment map} is 
$$
d\lan \mu, \xi\ran =i_{\widehat{s(\xi)}}\ome+i_{\widehat{R(\xi)}}\ome.
$$
The first term is the usual Hamiltonian term associated with the
$\mathfrak g_0$-action induced by $s$, whereas the second term is a
correction arising from the action of the abelian ideal $\mathfrak a$.
In this sense, a modified moment map can be regarded as a moment map modified by the abelian action of $\mathfrak a$.
This formulation is motivated by the moment map picture in generalized K\"ahler geometry, developed in the works \cite{Goto_2016}, \cite{Goto_2020}. In that setting, the Fujiki--Donaldson moment map formalism naturally leads to modified Hamiltonian actions. 
We conclude this subsection with a basic example of an abelian
extension. This type of extension appears in the Lie algebra of
generalized Hamiltonians considered in Section \ref{Generalized Hamiltonian diffeomorphisms} and 
\ref{GSCasmoment}.
\begin{example}\label{the abelian Lie algebra}
Let $\mathfrak g_0$ be a Lie algebra. Then the adjoint representation gives 
an abelian extension $\mathfrak g =\mathfrak g_0\rtimes_{\ad}\mathfrak a$, where the abelian Lie algebra $\mathfrak a$ is $\mathfrak g_0$ as a module and the Lie bracket of $\mathfrak g$ is given by 
$$
[(f_1, g_1), (f_2, g_2)] =([f_1, f_2],  [f_1, g_2]-[f_2, g_1]), \quad 
f_1, f_2\in \mathfrak g_0, \,\, g_1, g_2\in \mathfrak a.
$$
The associated short exact sequence
$$
0\longrightarrow\mathfrak a
\longrightarrow\mathfrak g
\longrightarrow\mathfrak g_0
\longrightarrow 0
$$
is split by the Lie algebra homomorphism
$
s(f)=(f,0).
$
\end{example}
\subsection{Moment map on the zero set $X^{\mathfrak a}$} 
Let $(X, \omega)$ and $\mathfrak g$ be as in the previous section, and assume that there exists a modified moment map $\mu: X\to \mathfrak g_0^*$ associated with a linear splitting $\rho : \mathfrak g_0\to\mathfrak g$.
We denote by $Y:=X^{\mathfrak a}$ the zero set of the infinitesimal action of the abelian Lie algebra $\mathfrak a$ on $X$, namely, 
$$
Y:=\{ x\in X \, |\,   \h a(x)=0, \,\,\, \text{\rm for every }a \in \frak a\}.
$$
Since $\mathfrak a$ is an ideal of $\mathfrak g$, $Y$ is invariant under the action of $\mathfrak g$.
Through the Lie algebra splitting $s: \mathfrak g_0\to \mathfrak g$, $Y$ is also invariant under the induced action of $\mathfrak g_0$.
We assume that $Y=X^{\mathfrak a}$ is a symplectic submanifold of $X$, with symplectic form given by the restriction $\omega|_Y=\omega_Y$.

\begin{proposition}\label{X^A}
Let  $Y=X^{\mathfrak a}$ be a symplectic submanifold as above. Then the restriction $\mu_Y:=\mu|_Y: Y\to\mathfrak g_0^* $ is an ordinary infinitesimally $\mathfrak g_0$-equivariant moment map.
\end{proposition}

\begin{proof}
Since $\rho(\xi)$ is written as $\rho(\xi)=s(\xi)+R(\xi)$, 
the vector field $\widehat{s(\xi)}|_Y$ is precisely the fundamental vector field on $Y$ generated by $\xi\in\mathfrak g_0$ under the descended $\mathfrak g_0$-action. 
Since $\rho(\xi)$ is written as $\rho(\xi)=s(\xi)+R(\xi)$ for $R(\xi)\in \mathfrak a$, we have 
that $\wideh{\rho(\xi)}|_Y =\wideh{s(\xi)}|_Y$. 
It follows from the modified moment map equation that
$$
\begin{aligned}
d\langle\mu_Y,\xi\rangle
&=d\langle\mu,\xi\rangle|_Y =\iota_{\widehat{\rho(\xi)}|_Y}\omega_Y =\iota_{\widehat{s(\xi)}|_Y}\omega_Y.
\end{aligned}
$$
Therefore $\mu_Y$ satisfies the usual moment map condition on $(Y,\omega_Y)$.
Moreover, since $\mu$ is infinitesimally $\mathfrak g_0$-equivariant, its restriction $\mu_Y$ is also infinitesimally $\mathfrak g_0$-equivariant. Hence $\mu_Y$ is an ordinary infinitesimally $\mathfrak g_0$-equivariant moment map.
\end{proof}
\begin{remark}
The notion of adapted volume in \cite{ASU_2026} can be explicitly interpreted in terms of the fixed-point set $X^{\mathfrak a}$.  In fact, $X$ corresponds to the space of $\J_\psi$-compatible almost generalized complex structures $\B_{\J_\psi}(M)$ and $\B_{\J_\psi}(M)$ admits a K\"ahler structure. 
The Lie algebra $\mathfrak g$ corresponds to the extension $\wtil\ham(M, \J_\psi)$ of $\mathfrak g_0=\ham(M, \J_\psi)$ by an abelian Lie algebra $\mathfrak a=C^\infty_\psi(M,\R)$. 
The action of the abelian Lie algebra $\mathfrak a$ is given by the action of $\sqrt{-1} gL^H_{\zeta_\a}$, where 
$g\in C^\infty_\psi(M,\R)$.
Then the subspace of $\B_{\J_\psi}(M)$ which admits a fixed volume form $\vol_M$ as the adapted volume form is given by the fixed-point set $\B^{\mathfrak a}_{\J_\psi}(M)$.
Under the identification of the canonical infinitesimal symmetry associated with the volume form with the action generated by 
$\zeta_\a$ , the adapted-volume condition is equivalent to
$$
L^H_{\zeta_\a}\J=0.
$$
 It follows from Proposition \ref{X^A} that a modified moment map is a moment map on 
$\B^{\mathfrak a}_{\J_\psi}(M)$
 (see Section \ref {GSCasmoment} for more details).
\end{remark}

\subsection{An obstruction to the existence of a moment map}
 Let a Lie algebra $\mathfrak g$ act on a symplectic manifold $X$ by symplectic vector fields as before.
 We assume that $X$ is connected and satisfies 
 $$H^1(X)=0.$$
 Then for every $\xi\in \mathfrak g$, the $1$-form
 $i_{\h\xi}\ome$ is $d$-exact. Hence there exists a Hamiltonian function $H_\xi$ such that 
 $$
 dH_\xi = i_{\h\xi}\ome \qquad (\xi\in \mathfrak g)
 $$
 Choosing the constants of integration so that the assignment $\xi\mapsto H_\xi$ is linear, we obtain a map 
 $$H: X\to\mathfrak g^*,   \quad \lan H, \xi\ran =H_\xi.
 $$
 Although $H$ satisfies the moment map equation, it need not be
infinitesimally $\mathfrak g$-equivariant.
Then the obstruction to the existence of a moment map is the well-known  
 $2$-cocycle of $\mathfrak g$
 $$
 c(\xi_1, \xi_2) =\{H_{\xi_1}, H_{\xi_2}\}- H_{[\xi_1, \xi_2]}.
 $$
Since $X$ is connected, $c(\xi_1, \xi_2)$ is a constant function on $X$. Then $c$ is regarded as an element $c$ of $Z^2(\mathfrak g, \R)$. 
This $2$-cocycle $c$ gives rise to a cohomology class $[c]\in H^2(\mathfrak g, \R)$, where $H^2(\mathfrak g, \R)$ denotes the second Lie algebra cohomology group of 
$\mathfrak g$. 
Changing the additive constants in the functions $H_\xi$ changes
$c$ by a Lie algebra coboundary. Thus, the class $[c]$ is independent
of these choices.
Then
there exists an infinitesimally $\frak g$-equivariant moment map for the action of $\mathfrak g$ on $X$ if and only if the class $[c]\in H^2(\mathfrak g, \R)$ vanishes.
\begin{proposition} \label{vanishing theorem} Let $(X,\ome)$ be a connected, symplectic manifold satisfying $H^1(X, \R)=0$, and let a Lie algebra $\mathfrak g$ act on $X$ by symplectic vector fields, where $\mathfrak g$ is an extension of a Lie algebra $\mathfrak g_0$ by an abelian Lie algebra $\mathfrak a$, that is, there is a short exact sequence of Lie algebras 
$$
0\longrightarrow \mathfrak a\longrightarrow  \mathfrak g\overset\pi{\longrightarrow}  \mathfrak g_0{\to} 0,
$$
equipped with a Lie algebra homomorphism $s: \mathfrak g_0\to \mathfrak g$ satisfying
$\pi\circ s=\id_{\mathfrak g_0}$ as before.
We assume that there exists a modified moment map $$\mu: X \to \mathfrak g_0^*$$ associated with a splitting $\rho: \mathfrak g_0 \to \mathfrak g$. If the zero set $Y=X^{\mathfrak a}$ is nonempty, then the class $[c]\in H^2(\mathfrak g, \R)$ vanishes and there exists an ordinary infinitesimally $\mathfrak g$-equivariant moment map for the action of $\mathfrak g$ on $X$.
\end{proposition}
\begin{proof}
Fix a point $x_0\in Y=X^{\mathfrak a}$. Then we choose the unique function $h_a$ which satisfies 
$$
dh_a=i_{\h a}\ome, \quad h_a(x_0)=0.
$$
Every $u\in\mathfrak g$ can be written uniquely in the form
 $$u =\rho(\xi)+a\in \mathfrak g, \quad\xi\in\mathfrak g_0,\,\, a\in\mathfrak a.$$  We define $H_u$ by 
$$
H_u=\lan\mu, \xi\ran +h_a.
$$
This assignment is linear in $u$, and the modified moment map equation
gives
$$
\begin{aligned}
dH_u
&=
d\langle\mu,\xi\rangle+dh_a
=\iota_{\widehat{\rho(\xi)}}\omega
+
\iota_{\widehat a}\omega
=\iota_{\widehat u}\omega.
\end{aligned}
$$
Let $c$ be the cocycle which is defined by this function $H_u$, 
$$
c(u_1, u_2)=\{ H_{u_1}, H_{u_2}\}-H_{[u_1, u_2]} \qquad (u_1, u_2\in \mathfrak g)
$$
Since $c(u_1, u_2)$ is a constant function on $X$, it suffices to evaluate it at $x_0\in Y$.
First, 
if $a_1, a_2\in \mathfrak a$, then it follows $[a_1, a_2]=0$ since $\mathfrak a$ is abelian.
Since $\{H_{a_1}, H_{a_2}\}(x_0)=\ome(\h a_1, \h a_2)(x_0)$ and $\h a_1(x_0)=\h a_2(x_0)=0$, we have $c(a_1, a_2)=0$.

Next, if $u_1=\rho(\xi)$, $ u_2=a\in \mathfrak a$, then we have $[\rho(\xi), a]\in \mathfrak a$ and $h_{[\rho(\xi), a]}(x_0)=0$.
Then we have 
$\{ H_{\rho(\xi)}, H_a\}(x_0)=\ome(\widehat{\rho(\xi)}, \h a)(x_0) =0$ since $\h a(x_0)=0$.

Finally, if $u_1=\rho(\xi_1)$, $u_2=\rho(\xi_2)$, then we have 
$$
[\rho(\xi_1), \rho(\xi_2)]=\rho([\xi_1, \xi_2])+F_\rho(\xi_1, \xi_2),    \qquad F_\rho(\xi_1, \xi_2)\in \mathfrak a.
$$
Since $\widehat{\rho(\xi)}(x_0)=\widehat{s(\xi)}(x_0)$ and $\mu$ is $\mathfrak g_0$-equivariant, we have 
\begin{align*}
\{H_{\rho(\xi_1)}, H_{\rho(\xi_2)}\}(x_0)=&dH_{\rho(\xi_1)}(\widehat{\rho(\xi_2)}(x_0))\\=&
d\lan\mu, \xi_1\ran(\widehat{s(\xi_2)}(x_0))\\
=&\lan \mu,  [\xi_1, \xi_2]\ran
\end{align*}
The normalization of the Hamiltonians associated
with elements of $\mathfrak a$ gives
$h_{F_\rho(\xi_1,\xi_2)}(x_0)=0.$
Therefore,
\begin{align}
H_{[\rho(\xi_1),\rho(\xi_2)]}(x_0)
&=
H_{\rho([\xi_1,\xi_2])}(x_0)+h_{F_\rho(\xi_1,\xi_2)}(x_0)\\
=&H_{\rho([\xi_1,\xi_2])}(x_0)
\end{align}
Hence we have 
\begin{align*}
c(\rho(\xi_1), \rho(\xi_2))=&
\{H_{\rho(\xi_1)}, H_{\rho(\xi_2)}\}(x_0)-H_{[\rho(\xi_1), \rho(\xi_2)]}(x_0)\\
=&\lan \mu,   [\xi_1, \xi_2]\ran (x_0)-H_{\rho ([\xi_1, \xi_2])}(x_0)\\
=&\lan \mu,   [\xi_1, \xi_2]\ran(x_0) -\lan \mu,   [\xi_1, \xi_2]\ran(x_0)=0
\end{align*}
By bilinearity, we obtain $c(u_1, u_2)=0$ for all $u_1, u_2\in \mathfrak g$.
Hence the class $[c]\in H^2(\mathfrak g, \R)$ vanishes. 
Therefore there exists an ordinary moment map for the action of $\mathfrak g$ on $X$.
\end{proof}
\section{Generalized complex structures and generalized K\"ahler structures}\label{Generalized complex structures and generalized Kahler structures}\label{GCS GK}
\subsection{Generalized complex structures and nondegenerate pure spinors}
Let $M$ be a manifold of real dimension $2n$, and denote by $TM$ and $T^*M$ its tangent and cotangent bundles, respectively. 
We set
$$
 \T M := TM \oplus T^*M .
$$
The generalized tangent bundle $\mathbb{T}M$ is equipped with the natural symmetric bilinear form
$$\lan v+\xi, u+\eta \ran_{\tt}=\frac12\(\xi(u)+\eta(v)\),\quad  u, v\in TM, \quad\xi, \eta\in T^*M .$$
Let $\mathrm{SO}(\mathbb{T}M)$ denote the subbundle of $\mathrm{End}(\mathbb{T}M)$ whose fibre is isomorphic to $\mathrm{SO}(2n,2n)$.
{\it An almost \complex structure} $\J$ is a section of $\mathrm{SO}(\mathbb{T}M)$ satisfying $\J^2=-\id.$ 
As in the case of ordinary almost complex structures, an almost generalized complex structure $\mathcal J$ yields the eigenspace decomposition
\bgn{equation}\label{eigenspace decomposition}
(\mathbb{T}M)^{\mathbb C}=\L_\J \oplus \ol \L_\J,
\end{equation} 
where $\L_{\mathcal J}$ is the $\sqrt{-1}$-eigenbundle of $\mathcal J$, and $\overline{\L}_{\mathcal J}$ is its complex conjugate.
Let $H$ be a real $d$-closed $3$-form.
The Courant bracket on sections of $\mathbb{T}M$ is defined by
$$
 [u+\xi, v+\eta]_{\cou}=[u,v]+{\mathcal L}_u\eta-{\mathcal L}_v\xi+\frac12 d\bigl(\iota_u\eta-\iota_v\xi\bigr),
 $$
 where $u, v\in \Gam(TM)$ and $\xi, \eta\in \Gam(T^*M)$.
 The $H$-twisted Courant bracket is given by
 $$
 [u+\xi, v+\eta]_{H}= [u+\xi, v+\eta]_{\cou}-\iota_v \iota_uH.
 $$
If $\L_\J$ is involutive with respect to the Courant bracket, then the almost generalized complex structure $\J$ is called {\it a generalized complex structure}, that is, $[e_1, e_2]_{\cou}\in \Gam(\L_\J)$  for any two sections 
 $e_1,e_2\in \Gam(\L_\J)$.
 If $\L_\J$ is involutive with respect to the $H$-twisted Courant bracket, the $\J$ is called {\it a $H$-twisted generalized complex structure}.
Let $\CL(\T M)$ be the Clifford algebra bundle which is 
a fibre bundle with fibre the Clifford algebra $\CL(2n, 2n)$ with respect to $\lan\,,\,\ran_{\tt}$ on $M$.
Then a vector $v$ acts on the space of differential forms $\oplus_{p=0}^{2n}\w^pT^*M$ by 
the interior product $\iota_v$ and a $1$-form $\t$ acts on $\oplus_{p=0}^{2n}\w^pT^*M$ by the exterior product $\t\w$, respectively.
Then the space of differential forms gives a representation of the Clifford algebra $\CL(\T M)$ which is 
the spin representation of $\CL(\T M)$. 
Thus
the spin representation of the Clifford algebra arises as the space of differential forms $$\w^\bullet T^*M=\oplus_p\w^pT^*M=\w^{\even}T^*M\oplus\w^{\odd}T^*M.$$ 
The spinor inner product $\lan\,,\,\ran_s$ of the spin representation is given by 
$$
\lan \a, \,\,\,\b\ran_s:=(\a\w\sig\b)_{[2n]},
$$
where $(\a\w\sig\b)_{[2n]}$ is the component of degree $2n$ of $\a\w\sig\b\in\oplus_p \w^pT^*M$ and 
$\sig$ denotes the Clifford involution which is given by 
$$
\sig\b =\bgn{cases}&+\b\qquad \deg\b \equiv 0, 1\,\,\mod 4 \\ 
&-\b\qquad \deg\b\equiv 2,3\,\,\mod 4\end{cases}
$$
We define $\ker\Phi:=\{ e\in (\T M)^\C\, |\, e\cdot\Phi=0\, \}$ for a differential form $\Phi
\in \w^{\even/\odd}T^*M.$
If $\ker\Phi$ is maximal isotropic, i.e., $\dim_\C\ker\Phi=2n$, then $\Phi$ is called {\it a pure spinor} of even/odd type.
A pure spinor $\Phi$ is {\it nondegenerate} if $\ker\Phi\cap\ol{\ker\Phi}=\{0\}$, i.e., 
$(\T M)^\C=\ker\Phi\oplus\ol{\ker\Phi}$.
Then a nondegenerate pure spinor $\Phi\in \w^\bullet T^*M$ gives an almost generalized complex structure $\J_{\Phi}$ which satisfies 
$$
\J_\Phi e =
\bgn{cases}
&+\sqrt{-1}e, \quad e\in \ker\Phi\\
&-\sqrt{-1}e, \quad e\in \ol{\ker\Phi}
\end{cases}
$$
Conversely, an almost \complex structure $\J$ locally arises as $\J_\Phi$ for a nondegenerate pure spinor $\Phi$ which is unique up to multiplication by
non-zero functions.  Thus an almost \complex structure yields {\it the canonical line bundle} $K_{\J}:=\C\lan \Phi\ran$ which is a complex line bundle locally generated by a nondegenerate, pure spinor $\Phi$ satisfying 
$\J=\J_\Phi$.
Let $d_H$ be a differential operator $d+H\w$ which acts on differential forms.
An almost \complex structure 
$\J_\Phi$ is a $H$-twisted generalized complex structure if and only if $d_H\Phi=\eta\cdot\Phi$ for a section $\eta\in \T M$. 
The {\it type number} of $\J=\J_\Phi$ is defined as the minimal degree of the differential form $\Phi$. Note that type number, denoted by Type$\J$  is a function on a manifold which is not a constant in general.
\bgn{example}
Let $J$ be a complex structure on a manifold $M$ and $J^*$ the complex structure on the dual  bundle $T^*M$ which is given by $J^*\xi(v)=\xi (Jv)$ for $v\in TM$ and $\xi\in T^*M$.
Then a \complex structure $\J_J$ is given by the following matrix
$$\J_J=\bgn{pmatrix}J&0\\0&-J^*
\end{pmatrix}$$
Then the canonical line bundle is generated by complex forms of type $(0,n)$.
Thus we have  Type $\J_J =n.$
\end{example}
\bgn{example}
Let $\ome$ be a symplectic structure on $M$ and $\h\ome$ the isomorphism from $TM$ to $T^*M$ given by $\h\ome(v):=\iota_v\ome$. We denote by $\h\ome^{-1}$ the inverse map from $T^*M$ to $TM$.
Then a \complex structure $\J_\psi$ is given by the following
$$\J_\psi=\bgn{pmatrix}0&-\h\ome^{-1}\\
\h\ome&0
\end{pmatrix}$$
Then the canonical line bundle is given by the differential form $\psi=e^{{\sqrt{-1}\ome}}$. 
Thus Type $\J_\psi=0.$
\end{example}
\bgn{example}[$b$-field action]
A real $d$-closed $2$-form $b$ acts on a \complex structure by the adjoint action of Spin group $e^b$ which provides
a \complex structure $\Ad_{e^b}\J=e^b\circ \J\circ e^{-b}$. 
\end{example}
\bgn{example}[Poisson deformations]\label{Poisson deformations}
Let $\b$ be a holomorphic Poisson structure on a complex manifold. Then the adjoint action of Spin group $e^\b$ gives deformations of new \complex structures by 
$\J_{\b t}:=\Ad_{\b^{Re} t}\J_J$.  Then Type ${\J_{\b t}}_x=n-2$ (rank of $\b_x$) at $x\in M$,
which is called the Jumping phenomena of type number.
\end{example}
Let $(M, \J)$ be a generalized complex manifold and $\ol \L_\J$ the eigenspace of eigenvalue $-\sqrt{-1}$.
Then we have the Lie algebroid complex $\w^\bullet\ol{\L}_\J$:
$$
0\arrow\w^0\ol \L_\J\overset{\ol\pa_\J}\arrow\w^1\ol \L_\J\overset{\ol\pa_\J}\arrow\w^2\ol \L_\J\overset{\ol\pa_\J}\arrow\w^3\ol \L_\J\arrow\cdots
$$
The Lie algebroid complex is the deformation complex of \complex structures. 
In fact, $\e\in \w^2\ol \L_\J$ gives deformed isotropic subbundle 
$E_\e:=\{ e+[\e, e]\, |\, e\in \L_\J\}$. 
Then $E_\e$ yields deformations of \complex structures if and only if $\e$ satisfies Generalized Maurer-Cartan equation
$$
\ol{\pa}_\J\e+\frac12[\e, \e]_{\Sch}=0,
$$
where $[\e, \e]_{\Sch}$ denotes the Schouten bracket. 
The Kuranishi space of generalized complex structures is constructed.
Then the second cohomology group $H^2(\w^\bullet\ol \L_\J)$ of the Lie algebraic complex gives the infinitesimal deformations of \complex structures and the third one 
$H^3(\w^\bullet\ol \L_\J)$ is the obstruction space to deformations of \complex structures.
\subsection{Generalized K\"ahler structures}
\bgn{definition}
{\it A generalized K\"ahler structure} is a pair $(\J_1, \J_2)$ consisting of two commuting \complex structures 
$\J_1$ and $\J_2$ such that $\h G:=-\J_1\circ\J_2=-\J_2\circ \J_1$ gives a positive definite symmetric form 
$G:=\lan \h G\,\,,  \,\,\ran$ on $\T M$. 
We call $G$ {\it a generalized metric}.
\end{definition}
Each $\J_i$ gives the decomposition $(\T M)^\C=\L_{\J_i}\oplus\ol \L_{\J_i}$ for $i=1,2$.
Since $\J_1$ and $\J_2$ are commutative, we have the simultaneous eigenspace decomposition 
$$
(\T M)^\C=(\L_{\J_1}\cap \L_{\J_2})\oplus (\ol \L_{\J_1}\cap \ol \L_{\J_2})\oplus (\L_{\J_1}\cap \ol \L_{\J_2})\oplus
(\ol \L_{\J_1}\cap \L_{\J_2}).
$$
Since $\h G^2=+\id$,
The generalized metric $\h G$ also gives the eigenspace decomposition: $\T M=C_+\oplus C_-$, 
where $C_\pm$ denote the eigenspaces of $\h G$ of eigenvalues $\pm1$. 
We denote by $\L_{\J_1}^\pm$ the intersection $\L_{\J_1}\cap C^\C_\pm$. 
Then it follows 
\bgn{align*}
&\L_{\J_1}\cap \L_{\J_2}=\L_{\J_1}^+,  \quad \ol \L_{\J_1}\cap \ol \L_{\J_2}=\ol \L_{\J_1}^+\\
&\L_{\J_1}\cap \ol \L_{\J_2}=\L_{\J_1}^-,\quad \ol \L_{\J_1}\cap \L_{\J_2}=\ol \L_{\J_1}^-
\end{align*}
\bgn{example}\label{ordinary Kahler}
Let $X=(M, J, \ome)$ be a K\"ahler manifold. Then the pair $(\J_J, \J_\psi)$ is a generalized K\"ahler where 
$\psi=\exp({\sqrt{-1}\ome})$. 
\end{example}
\bgn{example}
Let $(\J_1, \J_2)$ be a generalized K\"ahler structure. 
Then the action of $b$-fields gives a generalized K\"ahler structure 
$(\Ad_{e^b}\J_1, \Ad_{e^b}\J_2)$ for a real $d$-closed $2$-form $b.$
\end{example}
\bgn{definition}\label{Generalized Kahler structure of symplectic type}
 {\it A generalized K\"ahler structure of symplectic type} is a generalized K\"ahler structure $(\J, \J_\psi),$
where $\J_\psi$ is a generalized complex structure induced from a $d$-closed, nondegenerate, pure spinor $\psi =e^{b+\sqrt{-1}\ome}$ for a $d$-closed $2$-form $b$ and a symplectic structure $\ome.$
\end{definition}
\section{Generalized scalar curvature of generalized K\"ahler manifolds}
\label{Generalized Scalar curvature}
Let $M$ be a compact, orientable manifold of real dimension $2n$,  and $\vol_M$ a volume form of $M$. 
We denote by $H$ a real $d$-closed $3$-form. 
Let $(\J_\phi,\J_\psi)$ be an almost $H$-twisted generalized K\"ahler structure on $M,$
where $\phi=\{\phi_\a\}$ is a choice of local trivializations of the canonical line bundle $K_{\J_\phi}$ of 
generalized complex structure $\J_\phi.$ Namely,
each $\phi_\a$ is a nondegenerate pure spinor which induces $\J_\phi$ on an open set $U_\a, $ and $\{ U_\a\}$ is an open covering of $M.$
We also denote by $\psi=\{\psi_\a\}$  a choice of local trivializations of the canonical line bundle $K_{\J_\psi}$ relative to the cover $\{U_\a\}$, that is, each $\psi_\a$ is a
nondegenerate pure spinor on $U_\a$ which induces $\J_{\psi}.$ 

From now on,
we assume the following normalization for all $\a$,
\bgn{equation}\label{normalization}
\vol_M:=i^{-n}\lan\phi_\a, \ol\phi_\a\ran_s=i^{-n}\lan\psi_\a, \ol\psi_\a\ran_s,
\end{equation}
where $\vol_M$ is the fixed volume form on $M$. 
Note that a generalized K\"ahler structure gives an orientation of a manifold $M$.
Then we can achieve the normalization as in (\ref{normalization}),
if necessary, by multiplying nonzero real functions on each $\phi_\a$ and $\psi_\a$.
We denote by $d_H$ the differential operator $d+H$ which acts on differential forms on $M.$
Then the action of $d_H$ on $\phi_\a$ and $\psi_\a$ are respectively written as
$$
d_H\phi_\a=(\eta_\a+ N_\phi)\cdot\phi_\a, \quad d_H\psi_\a=(\zeta_\a+ N_\psi)\cdot\psi_\a,
$$
where $\eta_\a, \zeta_\a$ are chosen to be pure imaginary sections of $(\TT)^\C$, and $N_\phi, N_\psi\in \w^3(\TT)^\R$  are real.
With this convention, $\eta_\a, \zeta_\a$ are uniquely determined. 
If $N_\phi$ (resp. $N_{\psi})$ vanishes, then $\J_\phi$ (resp. $\J_\psi)$ is integrable as a generalized complex structure.
Thus $N_\phi, N_\psi$ are referred to as {\it the Nijenhuis tensors}. 
Note that $N_\phi, N_\psi$ are globally defined real skew-symmetric $3$-tensors of $\w^3(\TT)^\R.$
For each $e\in (\TT)^\C,$
we denote by $L^H_e$ the $H$-twisted Lie derivative  
 $L^H_e:=d_H e+ ed_H$ which acts on both differential forms and on sections of $\TT$.
We now define a real function $S(\J_\phi, \J_\psi,\vol_M)$ by taking the real part of the following:
\bgn{definition}\label{definition of scalar curvature}
\bgn{align*}
S(\J_\phi, \J_\psi,\vol_M)\vol_M:=&\Re\,\( i^{-n}\lan \psi_\a, \,\,\, L^H_{\eta_\a}\ol\psi_\a\ran_s +
     i^{-n}\lan \phi_\a, \,\,\, L^H_{\zeta_\a}\ol\phi_\a\ran_s\)\\
 &+ 2\lan \zeta_\a, \eta_\a\ran_{\scriptscriptstyle T\oplus T^*}\vol_M,
    \end{align*}
where $\ol\phi_\a$ and $\ol\psi_\a$ denote the complex conjugates of $\phi_\a$ and $\psi_\a,$
respectively.
We refer to $S(\J_\phi, \J_\psi,\vol_M)$ as {\it the generalized scalar curvature} of $(\J_\phi, \J_\psi)$.
For simplicity, we denote by $S(\J_\phi, \J_\psi)$ the generalized scalar curvature 
$S(\J_\phi, \J_\psi,\vol_M)$ for a fixed volume form $\vol_M$.
\end{definition}
Then we have the following well-definedness result:
\bgn{proposition}
The generalized scalar curvature
$S(\J_\phi, \J_\psi)$ is independent of the choices of local trivializations $\{\phi_\a\}, \{\psi_\a\}$ of  the canonical line bundles $K_{\mathcal J_\phi}$ and $K_{\mathcal J_\psi}$.
\end{proposition}
\bgn{proof}
We denote by $S(\phi_\a, \psi_\a)$ the generalized scalar curvature as in Definition \ref{definition of scalar curvature}
\bgn{align*}
S(\phi_\a, \psi_\a)\vol_M:=&\Re\,\( i^{-n}\lan \psi_\a, \,\,\, L^H_{\eta_\a}\ol\psi_\a\ran_s +
     i^{-n}\lan \phi_\a, \,\,\, L^H_{\zeta_\a}\ol\phi_\a\ran_s\)\\
 &+ 2\lan \zeta_\a, \eta_\a\ran_{\scriptscriptstyle T\oplus T^*}\vol_M,
   \end{align*}
Let $\phi_\a'$ and $\psi_\a'$ be another choice of local trivializations. 
Then it follows from the normalization condition (\ref{normalization}) that one has
 $\phi'_\a=e^{ip_\a}\phi_\a$ and $\psi'_\a=e^{iq_\a}\psi_\a$, where 
$p_\a$ and $q_\a$ are real smooth functions on $U_\a.$
Then it suffices to show $S(\phi_\a, \psi_\a)=S(\phi_\a', \psi_\a').$
In fact, we define $\eta_\a'$ and $\zeta_\a'$ by
$$
d_H\phi'_\a=\eta'_\a\cdot\phi'_\a+N'_\phi\cdot\phi', \quad d_H\psi'_\a=\zeta_\a'\cdot\psi_\a'+N'_\psi\cdot\psi'
$$
Then $\eta_\a'=\eta_\a+idp_\a$ and $\zeta_\a'=\zeta_\a+idq_\a.$
Then the Lie derivatives are given by 
$$
L^H_{\eta_\a'}\psi_\a'=L^H_{\eta_\a+ idp_\a}\psi_\a'=L^H_{\eta_\a}(e^{iq_\a}\psi_\a)
=e^{iq_\a}L^H_{\eta_\a}\psi_\a+(\eta_\a (idq_\a))e^{iq_\a}\psi_\a,
$$
where $(\eta_\a (idq_\a))$ denotes the coupling $\eta_\a$ and $i dq_\a.$
Taking the complex conjugate, we have 
$$
L^H_{\eta_\a'}\ol\psi_\a=e^{-iq_\a}L^H_{\eta_\a}\ol\psi_\a-(\eta_\a(idq_\a))e^{-iq_\a}\ol\psi_\a
$$
Thus 
$$
\lan \psi_\a', \,\,L^H_{\eta_\a'}\ol\psi_\a'\ran_s=\lan 
e^{iq_\a}\psi_\a, \,\,e^{-iq_\a}L^H_{\eta_\a}\ol\psi_\a\ran_s
-(\eta_\a (idq_\a))\lan \psi_\a, \,\,\ol\psi_\a\ran_s
$$
One also has 
$$
\lan \phi_\a', \,\,L^H_{\zeta_\a'}\ol\phi_\a'\ran_s=\lan 
e^{ip_\a}\phi_\a, \,\,e^{-ip_\a}L^H_{\zeta_\a}\ol\phi_\a\ran_s
-(\zeta_\a (idp_\a))\lan \phi_\a, \,\,\ol\phi_\a\ran_s
$$
Then we have
\bgn{align*}
S(\phi', \psi')\vol_M-2\lan\zeta'_\a, \,\,\eta'_\a\ran_{\tt}\vol_M =&
S(\phi, \psi)\vol_M-2\lan\zeta_\a, \,\,\eta_\a\ran_{\tt}\vol_M\\
-(\eta_\a (idq_\a))&i^{-n}\lan \psi_\a, \,\,\ol\psi_\a\ran_s
-(\zeta_\a(idp_\a))i^{-n}\lan \phi_\a, \,\,\ol\phi_\a\ran_s
\end{align*}
We also have
\bgn{align*}
2\lan\zeta'_\a, \,\,\eta'_\a\ran_{\tt}-2\lan\zeta_\a, \,\,\eta_\a\ran_{\tt}
=&2\lan idq_\a, \,\,\eta_\a\ran_{\tt}+2\lan \zeta_\a, \,\,idp_\a\ran_{\tt}\\
=&(\eta_\a (idq_\a))
+(\zeta_\a(idp_\a))\end{align*}
From the normalization condition (\ref{normalization}), we obtain 
$$
S(\phi', \psi')=S(\phi, \psi).
$$
\end{proof}
\bgn{lemma}\label{another description of GS}
Let $S(\J_\phi, \J_\psi)$ be the generalized scalar curvature of $(\J_\phi, \J_\psi)$ as in Definition \ref{definition of scalar curvature}.
Then $S(\J_\phi, \J_\psi)$ is also given by 
$$S(\J_\phi, \J_\psi)\vol_M:=\Re\( i^{-n}\lan \psi_\a, \,\,\, d_H({\eta_\a}\cdot\ol\psi_\a)\ran_s 
+i^{-n}\lan \phi_\a, \,\,d_H(\zeta_\a\cdot\ol\phi_\a)\ran_s \)
$$
\end{lemma}
\bgn{proof} 
Since the Lie derivative $L^H_{\eta_\a}$ is given by $d_H\circ \eta_\a+\eta_\a\circ d_H$, 
it follows 
\bgn{align*}
\Re\( i^{-n}\lan \psi_\a, \,\,\, L^H_{\eta_\a}\ol\psi_\a\,\ran_s\)=&
\Re\( i^{-n}\lan \psi_\a, \,\,\, d_H({\eta_\a}\cdot\ol\psi_\a)\ran_s\)\\
+&\Re\( i^{-n}\lan \psi_\a, \,\,\, {\eta_\a}\cdot d_H\ol\psi_\a\ran_s\)
\end{align*}
Since $d_H\ol\psi_\a=-\zeta_\a\cdot\ol\psi_\a+N_\psi\cdot\ol\psi_\a$ and 
$\lan \psi_\a, \,\,\eta_\a\cdot N_\psi\cdot\ol\psi_\a\ran_s=0$,  one has 
\bgn{align*}
\Re\( i^{-n}\lan \psi_\a, \,\,\, L^H_{\eta_\a}\ol\psi_\a\,\ran_s\)=&
\Re\( i^{-n}\lan \psi_\a, \,\,\, d_H({\eta_\a}\cdot\ol\psi_\a\)\ran_s\\
-&\Re\( i^{-n}\lan \psi_\a, \,\,\, {\eta_\a}\cdot\zeta_\a\cdot\ol\psi_\a\ran_s\)
\end{align*}
Since $d_H \ol\phi_\a=-\eta_\a\cdot\ol\phi_\a+N_\phi\cdot\ol\phi_\a$, 
one also has 
\bgn{align*}
\Re\( i^{-n}\lan \phi_\a, \,\,\, L^H_{\zeta_\a}\ol\phi_\a\,\)\ran_s=&
\Re\( i^{-n}\lan \phi_\a, \,\,\, d_H({\zeta_\a}\cdot\ol\phi_\a\)\ran_s\\
-&\Re\( i^{-n}\lan \phi_\a, \,\,\, {\zeta_\a}\cdot\eta_\a\cdot\ol\phi_\a\ran_s\)
\end{align*}
Since $\eta_\a$ is pure imaginary, it follows $\eta_\a=\eta_\a^{1,0}+\eta_\a^{0,1},$ where 
$\eta_\a^{1,0}\in \L_{\J_\psi}$,  $\eta_\a^{0,1}\in\ol L_{\J_\psi}$ and $\ol{\eta_\a^{0,1}}=-\eta_\a^{1,0}.$
Thus $2\,\Re\lan\eta_\a^{0,1}, \,\,\zeta_\a^{1,0}\ran_{\tt}=\Re\lan\eta_\a, \,\,\zeta_\a\ran_{\tt}$ 
Then it follows
\bgn{align*}
\Re\( i^{-n}\lan \psi_\a, \,\,\, {\eta_\a}\cdot\zeta_\a\cdot\ol\psi_\a\ran_s\)=&
\Re\( i^{-n}\lan \psi_\a, \,\,\, {\eta_\a^{0,1}}\cdot\zeta_\a^{1,0}\cdot\ol\psi_\a\ran_s\)\\
=&2\Re\lan\eta_\a^{0,1}, \,\,\zeta_\a^{1,0}\ran_{\tt} i^{-n}\lan \psi, \,\,\ol\psi\ran_s\\
=&\lan\eta_\a, \,\,\zeta_\a\ran_{\tt} i^{-n}\lan \psi, \,\,\ol\psi\ran_s\\
\end{align*}
We also have 
$$
\Re\( i^{-n}\lan \phi_\a, \,\,\, {\zeta_\a}\cdot\eta_\a\cdot\ol\phi_\a\ran_s\)=
\lan\zeta_\a, \,\,\eta_\a\ran_{\tt} i^{-n}\lan \phi, \,\,\ol\phi\ran_s
$$
Thus it follows 
\bgn{align*}
S(\J_\phi, \J_\psi)\vol_M=
&\Re\( i^{-n}\lan \psi_\a, \,\,\, L_H({\eta_\a}\cdot\ol\psi_\a)\ran_s+ i^{-n}\lan \phi_\a, \,\,\, L_H({\zeta_\a}\cdot\ol\phi_\a)\ran_s\)\\
&+2\lan \zeta_\a, \,\,\eta_\a\ran_{\tt}\\
=&\Re\( i^{-n}\lan \psi_\a, \,\,\, d_H({\eta_\a}\cdot\ol\psi_\a)\ran_s
+ i^{-n}\lan \phi_\a, \,\,\, d_H({\zeta_\a}\cdot\ol\phi_\a)\ran_s\)
\end{align*}
Hence the result follows.
\end{proof}
Let $\Diff_0(M)$ be the identity component of $\Diff(M).$
Define $\wtil\Diff_0(M)$ to be the extension of $\Diff_0(M)$ by $2$-forms
$$
0\to \Ome^2(M)\to \wtil\Diff_0(M)\to \Diff_0(M)\to 0
$$
The group $\wtil\Diff_0(M)$ acts on differential forms by pullback of diffeomorphisms together with the exterior product of $e^b$ of
the $b$-field.
We define a subgroup of $\wtil\Diff_0(M)$ by
\bgn{definition}\label{invarinat}
$$\wtil\Diff_0(M, H)=\{ g=e^bF\in \wtil\Diff_0(M)\, |\, d_H\circ g=g\circ d_H\, \}$$
\end{definition}
Note that $e^bF\in\wtil\Diff_0(M, H)$ if and only if 
$db=F^*H-H$.
Each $g=e^bF\in \wtil\Diff_0(M, H)$ acts on $\phi=\{\phi_\a\}$ by 
$g\cdot\phi=\{e^b\w F^*\phi_\a\}$ and also acts on $\psi=\{\psi_\a\}$ by 
$g\cdot\psi=\{e^b\w F^*\psi_\a\}.$
Then $(g\cdot\phi, g\cdot\psi)$ gives a generalized K\"ahler structure 
$(\J_{g\cdot\phi}, \,\J_{g\cdot\psi}).$
\bgn{proposition}\label{coequivariant }
$S(\J_\phi, \J_\psi,\vol_M)$ is equivariant  under the action of $\wtil\Diff_0(M, H),$ in the following sense, 
$$
S(\J_{g\cdot\phi}, \J_{g\cdot\psi}, F^*\vol_M)F^*\vol_M =F^*\(S(\J_\phi, \J_\psi,\vol_M)\vol_M\)
$$
for all $g=e^bF\in \wtil\Diff_0(M,H),$
\end{proposition}
\bgn{proof}
Since the action of $e^b$ preserves the spinor inner product $\lan\,,\,\ran_s$, 
we have 
\bgn{align}
i^{-n}\lan g\cdot\phi_\a, \ol{g\cdot\phi}\ran_s=&i^{-n}\lan e^b\w F^*(\phi_\a),  e^b\w F^*(\ol\phi_\a)\ran_s\\
=&i^{-n} F^*\lan \phi, \ol\phi\ran_s=F^*\vol_M
\end{align}
We also have $i^{-n}\lan g\cdot\psi_\a, \ol{g\cdot\psi_\a}\ran_s=F^*\vol_M$.
Thus we have the following normalization condition :
$$
F^*\vol_M =i^{-n}\lan g\cdot\phi_\a, \ol{g\cdot\phi}\ran_s=i^{-n}\lan g\cdot\psi_\a, \ol{g\cdot\psi_\a}\ran_s,
$$
we have
\bgn{align*}d_H(g\cdot \phi_\a)=&g d_H\phi_\a=g\circ\cdot (\eta_\a\cdot\phi_\a+N_\phi\cdot\phi_\a)\\
=&(g\circ\eta_\a\circ g^{-1})\cdot g\cdot\phi_\a+g\circ N_\phi\circ g^{-1}\circ g\cdot\phi_\a,
\end{align*}
where $g\circ\eta_\a\circ g^{-1}\in \sqrt{-1}(\TT).$
We also have
$$
d_H(g\cdot \psi_\a)=(g\circ\zeta_\a\circ g^{-1})\cdot g\cdot\psi_\a+g\circ N_\psi\circ g^{-1}\circ g\cdot\psi_\a
$$
Thus $\eta_\a$ and $\zeta_\a$ are changed by 
$(g\circ\eta_\a\circ g^{-1})$ and $(g\circ\zeta_\a\circ g^{-1})$, respectively.
From Lemma \ref{another description of GS},  we have
\bgn{align*}
S(\J_{g\cdot\phi}, \J_{g\cdot\psi}, F^*\vol_M)F^*\vol_M=&
\Re\( i^{-n}\lan g\cdot\psi_\a, \,\,\, d_H(g\circ{\eta_\a}\circ g^{-1}\cdot g\cdot\ol{\psi_\a})\ran_s \)\\
+&\Re\(i^{-n}\lan g\cdot\phi_\a, \,\,d_H(g\circ\zeta_\a\circ g^{-1}\cdot g\cdot\ol\phi_\a)\ran_s \)
\end{align*}
Since $d_H\circ g=g\circ d_H$, we have 
\bgn{align*}
S(\J_{g\cdot\phi}, \J_{g\cdot\psi}, F^*\vol_M)F^*\vol_M=&\Re\( i^{-n}\lan g\cdot\psi_\a, \,\,\,g\circ d_H({\eta_\a}\cdot\ol{\psi_\a})\ran_s \)\\
+&\Re\(i^{-n}\lan g\cdot\phi_\a, \,\,g\circ d_H(\zeta_\a\cdot\ol\phi_\a)\ran_s \)\\
=&\Re\,\, F^*\( i^{-n}\lan \psi_\a, \,\,\, d_H({\eta_\a}\cdot\ol{\psi_\a})\ran_s \)\\
+&\Re\,\, F^*\(i^{-n}\lan \phi_\a, \,\,d_H(\zeta_\a\cdot\ol\phi_\a)\ran_s \)\\
=&F^*S(\J_\phi, \J_\psi,\vol_M)F^*\vol_M
\end{align*}
\end{proof}
\bgn{proposition}\label{the ordinaryKahler}
(1) Let $(\J_J, \J_\ome)$ be a generalized K\"ahler structure which comes from the ordinary K\"ahler structure $(J, \ome)$ as in Example \ref{ordinary Kahler}. Then 
$S(\J_\phi, \J_\psi)$ coincides with the ordinary scalar curvature, where $vol_M=\frac {\ome^n}{n!}$\\
(2) If $(J_\phi, \J_\psi)$ is a GK of symplectic type, then 
$S(\J_\phi, \J_\psi)$ is a moment map as in \text{\rm\cite{Goto_2020}}, where $\J_\psi$ is given by 
a pure spinor $\psi=e^{b+\sqrt{-1}\ome}$ and
$\vol_M=\frac{\ome^n}{n!}$\\
(3) $S(\J_\phi, J_\psi)=S(\J_\psi, \J_\phi)$ is symmetric under the interchange of $\J_1$ and $\J_2$.
\end{proposition}
\bgn{proof}
(1) Let $(M, J, \ome)$ be a K\"ahler manifold. Then $\ome$ is locally written as 
$$
\ome=\sqrt{-1}\sum_{i,j} g_{i,\ol j} dz_i\w d\ol z_{ j}.
$$
Set $\rho=\log\det g_{i,\ol j}$.  Then 
$\phi=e^{\frac{\rho}{ 2}}d\ol {z_{1}}\w\cdots\w d\ol {z_{ n}}$ and 
$\psi=2^{\frac {-n}2}e^{i\ome}$ are non-degenerate, pure spinors which induce the generalized complex structure $\J_J$ and $\J_\psi$, respectively. Moreover they satisfy the normalization condition (\ref{normalization}): 
$$
i^{-n}\lan \phi, \ol\phi\ran =i^{-n}\lan \psi, \ol\psi\ran =\vol_M =\frac{\ome^n}{n!}
$$
Then we have 
$$d\phi=\frac 12(d\rho)\w\phi =\frac{\sqrt{-1}}2 \J_\phi d\rho \w\phi, \quad d\psi =0.
$$
Thus $\eta=\frac{\sqrt{-1}}2 \J_\phi d\rho$ and $\zeta=0$.
By the definition of the generalized scalar curvature, we obtain
$$
S(\J_\phi, \J_\psi )\vol_M =\Re\, i^{-n}\lan \psi, L_{\eta}\ol\psi\ran
$$
Then we have 
\bgn{align*}
\Re\, i^{-n}\lan \psi, L_{\eta}\ol\psi\ran=&\Re\,\( i^{-n}\psi\w\frac{\sqrt{-1}}2 \sig(dJ_\phi d\rho \w\ol\psi)\)_{[2n]}\\
=& -\frac12 d\J_\phi d\rho \w\frac{\ome^{n-1}}{(n-1)!}
\end{align*}
Since $J_\phi d\rho=\J_\phi(\pa\rho +\ol\pa\rho)=-\sqrt{-1}\pa\rho+\sqrt{-1}\,\ol\pa\rho$,  we have 
$$
-\frac12 d\J_\phi d\rho=\sqrt{-1}\,\ol\pa\pa\rho =\Ric_\ome.
$$
Hence
$$
S(\J_\phi, \J_\psi)\vol_M =\Ric_\ome \w\frac{\ome^{n-1}}{(n-1)!}
=\text{\rm Scal}_\ome \,\,\vol_M
$$
Statements (2) and (3) follow directly from our definition.
\end{proof}
Even in the case of an ordinary K\"ahler manifold, if the volume form is replaced by 
$e^f \frac{\omega^n}{n!}$ for a function $f$,
then the generalized scalar curvature is given as follows. 
\bgn{proposition}
Let $(\J_\phi, \J_\psi)$ be the generalized K\"ahler structure induced by the ordinary K\"ahler structure $(J, \ome, g)$ as before.
The generalized scalar curvature $S(\J_\phi, \J_\psi,\vol_M )$ with respect to $\vol_M =e^f\frac{\ome^n}{n!}$ is given by 
$$
S(\J_\phi, \J_\psi, \vol_M )\vol_M =e^f(\Ric_\ome+\sqrt{-1}\ol\pa\pa f)\w \frac{\ome^{n-1}}{(n-1)!}= \(\text{\rm Scal}_\ome+\frac14 \trian_g f \) \vol_M
$$
\end{proposition}

\bgn{remark}
Since a generalized K\"ahler structure gives the Born-Infeld volume form $\vol_{BI}:=\sqrt{\det(g+B)}$, one can define 
$S(\J_\phi, \J_\psi, \vol_{BI})$. 
However, in order to obtain the moment map framework, we fix $\J_\psi$ and choose
$\vol_M$  to be $i^{-n}\lan \psi_\a, \ol\psi_\a\ran_s$. (See Section \ref{GSCasmoment}).
\end{remark}
\bgn{remark}
Boulanger also obtained the moment map in the cases of toric generalized K\"ahler manifolds of symplectic type by using 
a description of toric geometry \text{\rm\cite{Bou_2019}}. 
Although Boulanger's description of the moment map seems to be different from the one in  \text{\rm\cite{Goto_2020},} 
these should match each other since the moment map is unique modulo constant.
In fact, 
Wang, Yicao actually shows that these are the same by using explicit calculations \text{\rm\cite{WY_2020}}.
\end{remark}
\bgn{remark}
There is a one to one correspondence between bihermitian structures $(I_+, I_-, g, b)$ and generalized K\"ahler structures
$(\J_1,\J_2).$
In Supergravity \text{\rm\cite{CSW_2011}}
, the notion of generalized scalar curvature is introduced, which is given by 
depending an arbitrary function 
$f$, \text{\rm\cite{Fer_2014}}
$$
GS^f(J)=s_g+4\trian_gf-4|df|^2-\frac12|db|^2
$$
In dimension $4$, Boulanger shows that 
our scalar curvature $S(\J, \J_\ome)=GS^f(J)$ if $f=-\frac12\log(1-p),$ where $p=-\frac14 tr (I_+I_-)$ is the angle function. 
\end{remark}
\bgn{remark}
J.~Streets studies problems of generalized K\"ahler structures by using  pluriclosed flows \cite{St_2016}. 
In the cases of generalized K\"ahler structures of type $(0,0)$, (which is also called a degenerate generalized K\"ahler structure),
his definition of 
generalized K\"ahler structure with constant scalar curvature is the same as the one in our paper 
(see also \cite{Goto_2016}, for generalized K\"ahler structures of type $(0,0)$).
The Calabi-Yau type problem of generalized K\"ahler manifolds of type $(0,0)$ was discussed by Apostolov and Streets in
\cite{AS_2017}.
A generalized K\"ahler Ricci flow has been explored
\cite{MGFJS_2020}, \cite{ASU_2021}. It is a remarkable problem to investigate a relation between the scalar curvature as moment map in this paper and the generalized K\"ahler Ricci flow.
\end{remark}
\begin{remark}
Apostolov, Streets, and Ustinovskiy have given an explicit formula for the generalized scalar curvature in terms of a bi-Hermitian structure $(g,I, J, b)$ \cite {ASU_2026}.
$$
S(\J_\phi, \J_\psi, e^{-f}\vol_g)=\frac14\(R-\frac1{12}|H|^2+2\trian f -|df|^2\),
$$
where $H=-d_I^C\ome_I=d_J^c\ome_J=H_0+db$ and $R$ denotes the scalar curvature.
\end{remark}

\section{Generalized Hamiltonians and an extended action}
\label{Generalized Hamiltonian diffeomorphisms}


We shall introduce generalized Hamiltonian diffeomorphisms on a $H$-twisted generalized complex manifold $M$. 
All manifolds are assumed to be connected unless otherwise stated.
Let $\J_\psi$ be a $H$-twisted generalized complex structure on $M$ which satisfies 
$$
d_H\psi_\a =\zeta_\a\cdot\psi_\a.
$$
Note that in this section we assume that $\J_\psi$ is integrable. 
\bgn{definition}
Let $\pi_T :\TT\to TM$  be the projection to the tangent $TM$.
We define a function $\{ f_1, f_2\}_{\psi}$  by 
$$\{ f_1, f_2\}_{\psi}:= \frac12 \(L_{\pi_T(\J_\psi df_1)}f_2-L_{\pi_T(\J_\psi df_2)}f_1\), $$
where $ f_1, f_2\in C^\infty(M)$
 \end{definition}
This bracket is a generalization of the Poisson bracket of symplectic geometry.
In fact, if $\J_\psi$ is induced from a symplectic structure, then $\{\,,\,\}_\psi$ coincides with the Poisson bracket.
We also have 
\bgn{lemma}
$$\{f_1, f_2\}_\psi =L_{\pi_T(\J_\psi df_1)} f_2 =-L_{\pi_T(\J_\psi df_2)} f_1.$$
\end{lemma}
\bgn{proof}
Since $\lan\, ,\,\ran_{\tt}$ is $\J_\psi$-invariant, we have
$$
L_{\pi_T(\J_\psi df_1)} f_2=\lan \J_\psi df_1, df_2\ran_{\tt}=\lan df_1,  -J_\psi df_2\ran_{\tt}=-L_{\pi_T(\J_\psi df_2)}f_1
$$
Then the result follows.
\end{proof}
\bgn{remark}
In general,  note that $\{f_1, f_2\}_\psi \neq L_{\J_\psi df_1} f_2 \neq -L_{\J_\psi df_2} f_1$
\end{remark}
 
We denote by $\pi_T: \TT\to T_M$ the projection to the tangent bundle. 
Then it  follows that $\pi_T(\zeta_\a)$ does not depend on the choice of trivializations of $K_{\J_\psi}$ and thus 
one has a vector field $\pi_T(\zeta_\a)$ on $M,$ which is globally defined.
We consider a set of real smooth functions $ C^\infty_\psi(M, \R)$ which is given by 
 $$
 C^\infty_\psi(M, \R):=\{\, f\in C^\infty(M,\R)\, |\, \pi_T(\zeta_\a)f=0,  \}
 $$
 \bgn{remark}
We do not impose the condition $\int_M f\vol_M=0$ in the definition of $C^\infty_\psi(M,\R)$.
Thus a nonzero constant function is included in $C^\infty_\psi(M,\R)$.
\end{remark}
 \bgn{lemma}
 We define $\ham_0(M,\J_\psi)$ by 
 $$\ham_0(M, \J_\psi):=\{\, \J_\psi df\in \TT\, |\, f\in C^\infty_\psi(M, \R)\, \}$$
Then $\ham_0(M, \J_\psi)$ is closed under the Courant bracket.
\end{lemma}
\bgn{proof}
We denote by $[\,,\,]_H$ the $H$-twisted Courant bracket. Since $\J_\psi$ is integrable, we have
$$
[\J_\psi df_1, \J_\psi df_2]_{H}=[df_1, \,\, df_2]_H+\J_\psi[\J_\psi df_1, \,\ df_2]_H
+\J_\psi[df_1, \,\, \J_\psi df_2]_H
$$
Since $[df_1, \,\, df_2]_H=[df_1, \,\, df_2]_{\cou}=0,$
we have 
\bgn{align}\label{hamiltonian of Jpsi}
[\J_\psi df_1, \J_\psi df_2]_{H}=&\J_\psi[\J_\psi df_1, \,\ df_2]_{\cou}
+\J_\psi[df_1, \,\, \J_\psi df_2]_{\cou}\\
=&\frac12\J_\psi[\J_\psi df_1, \,\ df_2]_{Dor}
-\frac12\J_\psi[ \J_\psi df_2,\,\, df_1]_{Dor}\notag\\
=&\frac12\J_\psi d (L_{\pi_T(\J_\psi df_1)}f_2-L_{\pi_T(\J_\psi df_2)}f_1)\notag\\
=&\J_\psi d\{ f_1, f_2\}_{\psi}\notag,
\end{align}
where $[\,,\,]_{Dor}$ denotes the Dorfman bracket.
Since $[\zeta_\a, \J_\psi df]_{Dor,H}=[L^H_{\zeta_\a}, \J_\psi df]=0,$
one also has 
\bgn{align*}
2\pi_T\zeta_\a\{f_1, f_2\}_\psi=&2L_{\zeta_\a}^H\{f_1, f_2\}_\psi =L^H_{\zeta_\a}
(L_{\J_\psi df_1}^Hdf_2-L^H_{\J_\psi df_2}df_1)\\
=&L^H_{[\zeta_\a, \J_\psi df_1]_H}df_2-L^H_{[\zeta_\a, \J_\psi df_2]_H}df_1=0
\end{align*}
Thus we obtain the result.
\end{proof}
Since $\ham_0(M,\J_\psi)$ is isotropic, it follows that $\ham_0(M,\J_\psi)$ is a Lie algebra.
Then $C^\infty_\psi(M,\R)$ is also a Lie algebra with the bracket $\{\,,\,\}_\psi$ and the map $C^\infty_\psi(M,\R)$ to $\ham_0(M,\J_\psi)$ is given by $f\mapsto \J_\psi df$ is a Lie algebra homomorphism.
The Lie algebra $\ham_0(M,\J_\psi)$ is identified with $C^\infty_\psi(M)/\R$. 
\bgn{definition}
We define $\ham(M,\J_\psi)$ to be a Lie algebra $C^\infty_\psi(M,\R)$ with the bracket $\{\,,\}_\psi$.
\end{definition}

As in Example \ref{the abelian Lie algebra}, 
we obtain the abelian extension $\wtil{\ham}(M,\J_\psi)$ of 
the Lie algebra $\ham(M,\J_\psi)$,
$$
0\to \mathfrak a \to \wtil{\ham}(M,\J_\psi)\to \ham(M,\J_\psi)\to 0,
$$
where the abelian Lie algebra $\mathfrak a$ is given by 
$C^\infty_\psi(M)\cong \ham(M,\J_\psi)$ as a module.
We denote by $(f, g)$ an element of $\wtil{\ham}(M,\J_\psi)$, where $f\in {\ham}(M,\J_\psi)$ and $g\in \mathfrak a$. 
The bracket of $\wtil{\ham}(M,\J_\psi)$ is given by 
$$
[(f_1, g_1), (f_2, g_2)] =(\{f_1, f_2\},  \{f_1, g_2\}-\{f_2, g_1\}), \quad 
f_1, f_2\in \ham(M,\J_\psi),  g_1, g_2\in \mathfrak a.
$$
In order to obtain a moment map framework, we need an extension of  the action of generalized Hamiltonians
in terms of $\zeta_\a,$ 
\bgn{definition}\label{def:lifted action}
Each $(f, g)\in \wtil{\ham}(M,\J_\psi)$ gives rise to the following 
action:
\bgn{equation}\label{lifted action}
\wtil L^H_{f,g}:= L^H_{\J_\psi df}+\sqrt{-1}gL^H_{\zeta_\a}
\end{equation}
Then $(f, g)$ acts on both sections of $\TT$ and differential forms on $M$ by $H$-twisted Lie derivative, together with multiplication by a function $\sqrt{-1}g$. \end{definition}
The action  $\wtil L^H_{f,g}$ is called the extended action.
If $\zeta_\a=0$, it follows $\wtil L^H_{f,g}:= L^H_{\J_\psi df}$ which is the ordinary $H$-twisted Lie derivative.
However, if $\zeta_\a\neq 0,$ we need to modify the ordinary definition of the action of Hamiltonians.
Since $\zeta_\a$ is pure imaginary, $\wtil L^H_{f,g}$ is defined as a real differential operator.
Since 
$\zeta_\a =\zeta_\b+\sqrt{-1}dq_{\a, \b}$, it follows 
$L_{\zeta_\a}^H=L_{\zeta_\b}^H$. Thus $\wtil L^H_{f,g}$ is well-defined on $M$. 
We need the following lemmas in order to show that the extended action gives a representation of the Lie algebra $\wtil{\ham}(M,\J_\psi)$.
\bgn{lemma}
The operator bracket of  $L_{e_1}$ and $L_{e_2}$ is denoted by
$[{L_{e_1}}, \,\, L_{e_2}]:={L_{e_1}}{L_{e_2}}-{L_{e_2}}{L_{e_1}}$, for $e_1, e_2\in \TT$.
Then one has
$$
[L_{e_1}, \,\, L_{e_2}]=L_{[e_1, \,\, e_2]_{\cou}}
$$
\end{lemma}
\bgn{proof}
The operator bracket is given in terms of  the Dorfman bracket of $e_1, e_2,$
$$
[L_{e_1}, \,\, L_{e_2}]=L_{[e_1, \,\, e_2]_{Dor}}
$$
Since the difference between the Dorfman bracket and the Courant bracket is given by 
${[e_1, \,\, e_2]_{Dor}}-{[e_1, \,\, e_2]_{\cou}}=d\lan e_1, e_2\ran_{\tt},$
we have 
$L_{[e_1, \,\, e_2]_{Dor}}=L_{[e_1, \,\, e_2]_{\cou}}.$
Thus we have $[L_{e_1}, \,\, L_{e_2}]=L_{[e_1, \,\, e_2]_{\cou}}$.
\end{proof}
\bgn{lemma}\label{[LHe1, LHe_2]}
$[L^H_{e_1},\,\, L^H_{e_2}]=L_{[e_1, e_2]_H}^H$
\end{lemma}  
\bgn{proof}
For $e_1=v_1+\t_1, e_2=v_2+\t_2\in \TT$, one has 

\bgn{align*}
[L^H_{e_1}, \,\, L^H_{e_2}]=&L_{[e_1, e_2]_{\cou}}+di_{v_1}i_{v_2}H+ i_{[v_1, v_2]}H\\
=&L_{[e_1, e_2]_H}+i_{[v_1, v_2]}H\\
=&L^H_{[e_1, e_2]_H}
\end{align*}
\end{proof}
\bgn{lemma}\label{LHzetaapsiJpsi=0}
$L^H_{\zeta_\a}\J_\psi=0$
\end{lemma}
\bgn{proof}
It suffices to show that 
$L^H_{\zeta_\a}\psi_\a\in K_{\J_\psi}$.
Since $d_H\psi_\a=\zeta_\a\cdot\psi_\a$ and $d_H\circ d_H=0,$ we have 
\bgn{align*}
L^H_{\zeta_\a}\psi_\a=&d_H (\zeta_\a\cdot\psi_\a)+\zeta_\a\cdot d_H\psi_\a\\
=&d_H d_H \psi +\zeta_\a \cdot\zeta_\a \cdot\psi_\a\\
=&\lan \zeta_\a, \,\, \zeta_\a\ran_{\tt} \psi_\a
\end{align*}
Since $\lan \zeta_\a, \,\, \zeta_\a\ran_{\tt}$ is a function, 
it implies that $L_{\zeta_\a}^H\psi_\a\in K_{\J_\psi}.$
Thus $L_{\zeta_\a}^H$ preserves $\J_\psi.$
Then we have $L^H_{\zeta_\a}\J_\psi=0.$
\end{proof}
\bgn{proposition}\label{the commutator of two differential operators}
For $(f_1,g_1)$, $(f_2, g_2)\in \wtil{\ham}(M,\J_\psi)$, 
the commutator of two differential operators $\wtil L^H_{f_1, g_1}$
and $\wtil L^H_{f_2, g_2}$ is given by the following
$$
[\wtil L^H_{f_1, g_1}, \,\, \wtil L^H_{f_2,g_2}]
=L^H_{\J_\psi d\{ f_1, \, f_2\}_\psi}+ \sqrt{-1}\( \{f_1, g_2\}_\psi-\{f_2, g_1\}_\psi \)L^H_{\zeta_\a},
$$
Namely we have 
$$
[\wtil L^H_{f_1,g_1},\wtil L^H_{f_2,g_2}]
=
\wtil L^H_{f_3, g_3}, 
$$
where $f_3=  \{f_1,f_2\}_\psi,\, g_3=
 \{f_1,g_2\}_\psi-\{f_2,g_1\}_\psi.$
\end{proposition} 
\bgn{proof}
From Lemma \ref{[LHe1, LHe_2]} and (\ref{hamiltonian of Jpsi}), we have 
$$
[L^H_{\J_\psi df_1}, \,\, L^H_{\J_\psi df_2}]=L^H_{[\J_\psi df_1, \,\, \J_\psi df_2]_H}
=L^H_{\J_\psi d\{f_1, \,\,f_2\}_\psi}
$$
Since $f_1, g_2\in C_\psi^\infty(M,\R)$, we have 
\bgn{align*}
[L^H_{\J_\psi df_1}, \,\, \sqrt{-1}g_2L^H_{\zeta_\a}]=&\sqrt{-1}(L_{\pi_T(\J_\psi df_1)}g_2)L_{\zeta_\a}^H+\sqrt{-1}g_2[L^H_{\J_\psi df_1}, \,\, L^H_{\zeta_\a}]\\
=&\sqrt{-1}\pi_T(\J_\psi df_1)g_2 L^H_{\zeta_\a}+\sqrt{-1}g_2L^H_{[\J_\psi df_1, \,\, \zeta_\a]_H}
\end{align*}
Since $L^H_{[\J_\psi df_1, \,\, \zeta_\a]_{co}}=
-L^H_{[ \zeta_\a,\,\,\J_\psi df_1 \,]_{Dor}}$, 
we have 
\bgn{align*}
L^H_{[\J_\psi df_1, \,\, \zeta_\a]_H}=
&-L^H_{[\zeta_\a,\,\,\J_\psi df_1 \, ]_H}\\
=&-L^H_{[ \zeta_\a,\,\,\J_\psi df_1 \,]_{Dor}}
+L^H_{(\pi_T (\J_\psi d f_1)\cdot\pi_T(\zeta_\a)\cdot H
)}
\end{align*}
Then from Lemma \ref{LHzetaapsiJpsi=0} and  $L_{\zeta_\a}f_1=0$, we have 
$[L^H_{\zeta_\a},\J_\psi d f_1]=0$. 
Since
$L_{\zeta_\a}^H=L_{\zeta_\a}+\pi_T(\zeta_\a)\cdot H$, we have 
\bgn{align*}
[\zeta_\a, \,\,\J_\psi df_1]_{Dor}=&[L_{\zeta_\a}, \J_\psi d f_1]=[L^H_{\zeta_\a},\J_\psi d f_1]
+\pi_T (\J_\psi d f_1)\cdot\pi_T(\zeta_\a)\cdot H\\
=&\pi_T (\J_\psi d f_1)\cdot\pi_T(\zeta_\a)\cdot H
\end{align*}
Thus we have
\bgn{align}
L^H_{[\J_\psi df_1, \,\, \zeta_\a]_H}=0
\end{align}
Hence we obtain 
\bgn{align}
[L^H_{\J_\psi df_1}, \,\, \sqrt{-1}g_2L^H_{\zeta_\a}]=&\sqrt{-1}(L_{\pi_T(\J_\psi df_1)}g_2)L^H_{\zeta_\a}
\end{align}
We also have 
\bgn{align}
 \sqrt{-1}[g_1L^H_{\zeta_\a}, \,\, L^H_{\J_\psi df_2}]=-\sqrt{-1}&(L_{\pi_T(\J_\psi df_2)}g_1)L^H_{\zeta_\a}
\end{align}
Then it follows 
$$
[L^H_{\J_\psi df_1}, \,\, \sqrt{-1}g_2L^H_{\zeta_\a}]+ [\sqrt{-1}g_1L^H_{\zeta_\a}, \,\, L^H_{\J_\psi df_2}]
=\sqrt{-1}\(\{f_1, g_2\}_\psi-\{f_2,g_1\}_\psi \)L^H_{\zeta_\a}
$$
Since $L_{\zeta_\a}^H f_1=L_{\zeta_\a}^H f_2=0,$ we have 
\bgn{align*}
[\wtil L^H_{f_1,g_1}, \,\, \wtil L^H_{f_2,g_2}]
=&L^H_{\J_\psi d\{ f_1, f_2\}_\psi}+\sqrt{-1}\(\{f_1, g_2\}_\psi-\{f_2,g_1\}_\psi \) L^H_{\zeta_\a}\\
=&\wtil L^H_{[(f_1, g_1), (f_2, g_2)]}.
\end{align*}
\end{proof} 
\begin{corollary}
The extended action $\wtil L^H_{f,g}$ gives a representation of the Lie algebra $\wtil{\ham}(M,\J_\psi)$.
\end{corollary}
\begin{proof}
The result follows from Proposition \ref{the commutator of two differential operators} and Example \ref{the abelian Lie algebra}.
\end{proof}
\bgn{proposition}\label{LHJPSIdfJpsi}
We also have
$$
L^H_{\J_\psi df}\J_\psi=0
$$
\end{proposition} 
 \bgn{proof}
 Since $d_H \psi_\a =\zeta_\a\cdot\psi_\a$ and $(\J_\psi df)\cdot\psi=-\sqrt{-1}df\cdot\psi$, 
 we have 
 \bgn{align*}
 L^H_{\J_\psi df}\psi_\a=&d_H (\J_\psi df)\cdot \psi_\a+ (\J_\psi df) d_H \psi_\a\\
 =&-d_H (\sqrt{-1} df)\cdot \psi_\a+ (\J_\psi df) d_H \psi_\a\\
 =&(\sqrt{-1} df)\cdot d_H \psi_\a+ (\J_\psi df) d_H \psi_\a\\
 =&(\sqrt{-1}df +\J_\psi df)\cdot \zeta_\a\cdot\psi_\a
 \end{align*}
 Since $(\sqrt{-1}df +\J_\psi df)\in \L_\psi,$ we have 
 $(\sqrt{-1}df +\J_\psi df)\cdot \zeta_\a\cdot\psi_\a
\in K_{\J_\psi}.$
Thus $L_{\J_\psi df}^H \psi_\a\in K_{\J_\psi}.$
Hence we have $L_{\J_\psi df}^H \J_\psi=0.$
 \end{proof} 
 \bgn{lemma}\label{wtil LHJpsi df=0}
 $$
\wtil L^H_{f,g}\J_\psi=0
$$
 \end{lemma}
 \bgn{proof}
 The result follows from Proposition \ref{LHJPSIdfJpsi} and
 Lemma \ref{LHzetaapsiJpsi=0}.
 \end{proof}
From Proposition \ref{the commutator of two differential operators}, the modified action $(f,g)\mapsto \wtil L^H_{f,g}$ is a Lie algebra homomorphism. 

\section{Generalized scalar curvature as modified moment map}\label{GSCasmoment}
\subsection{Moment map framework}
Let $M$ be a $2n$-dimensional manifold with a volume form $\vol_M$.
We denote by ${\mathcal B}(M)$ the set of almost \complex structures on $M$,
that is, 
$${{\mathcal B}}(M):=\{ \J\,:\,\text{\rm almost \complex structure on }M\,\}.$$
We also define ${\mathcal B}^{\rm int}(M)$ as the set of $H$-twisted \complex structures on $M$, i.e., integrable ones
$${\mathcal B}^{\text{\rm int}}(M):=\{ \J\, :\text{\rm  $H$-twisted \complex structure on }M\,\}.$$
We fix a $H$-twisted \complex structure $\J_\psi$ which is defined by 
a set of nondegenerate, pure spinors $\psi:=\{\psi_\a\}$ relative to a cover $\{U_\a\}$ of $M.$
Then we have 
\bgn{equation}\label{dpsia=zetaacdotpsi}
d_H\psi_\a =\zeta_\a\cdot\psi_\a,
\end{equation}
where $\zeta_\a\in \sqrt{-1}(\TT).$
We can take $\{\psi_\a\}$ which satisfies 
$i^{-n}\lan \psi_\a, \,\,\ol\psi_\a\ran_s=
i^{-n}\lan \psi_\b, \,\,\ol\psi_\b\ran_s=\vol_M$ if $U_\a\cap U_\b\neq\emptyset$. 

An almost \complex structure $\J$ is {\it $\J_\psi$-compatible} if and only if the pair $(\J, \J_\psi)$
is an almost generalized K\"ahler structure.
Let ${\mathcal B}_{\J_\psi}(M)$ be the set of $\J_\psi$-compatible almost \complex structure,   
$${\mathcal B}_{\J_\psi}(M):=\{\, \J\in {{\mathcal B}(M)}\, :\, (\J,\J_\psi)\,\,\text{\rm is an almost generalized K\"ahler structure}\, \}.$$
We also define ${\mathcal B}^{\rm int}_{\J_\psi}(M)$ to be the set of $\psi$-compatible \complex structures.
For each point $x\in M$, we define ${\mathcal B}_{\J_\psi}(M)_x$ to be the set of $\psi_x$-compatible almost \complex structures on $T_xM\oplus T^*_xM$ ,  
 $${{\mathcal B}}_{\J_\psi}(M)_x:=\{\, \J_x\, |(\J_x, \J_{\psi, x}): \text{\rm almost generalized K\"ahler structure at } x \, \}.$$
 Then it follows that ${\mathcal B}_{\J_\psi}(M)_x$ is given by the Riemannian Symmetric space of type  AIII
$$U(n,n)/(U(n)\times U(n))$$ which is biholomorphic to the complex bounded domain 
 $\{\, h\in M_n(\C)\, |\, 1_n-h^*h>0\, \},$ where $M_n(\C)$ denotes the set of complex matrices of $n\times n.$ 
 \bgn{remark} In K\"ahler geometry, the set of almost complex structures compatible with a symplectic structure $\ome$ is given by 
 the Riemannian symmetric space Sp$(2n, \R)/U(n)$ which is biholomorphic to the Siegel upper half plane 
 $$\{\, h\in M_n(\C)\, |\, 1_n-h^*h>0,\, h^t=h \,\}$$
 \end{remark}
Let $P_{\J_\psi}$ be the fibre bundle over $M$ with fibre ${{\mathcal B}}_{\J_\psi}(M)_x$, that is, 
$$P_\psi:=\bigcup_{x\in M}{{\mathcal B}_{\J_\psi}(M)_x}\to M,$$
Then ${\mathcal B}_{\J_\psi}(M)$ is given by smooth sections $\Gam (M, P_{\J_\psi})$ which contains the integrable generalized complex structures ${\mathcal B}^{\text{\rm int}}_{\J_\psi}(M)$. 
We can introduce a Sobolev norm on ${\mathcal B}_{\J_\psi}(M)$ such that ${\mathcal B}_{\J_\psi}(M)$ becomes a Banach manifold by the standard method.
The tangent bundle of ${\mathcal B}_{\J_\psi}(M)$ at $\J$ is given by 
$$T_{\J}{\mathcal B}_{\J_\psi}(M)=\{\, \dot{\J}\in\text{\rm so}(T_M\oplus T^*_M)\,:\, \dot{\J}\J+\J\dot{\J}=0,\, \dot{\J}\J_\psi=\J_\psi\dot{\J}\, \},$$
where so$(\TT)$ denotes the set of sections of Lie algebra bundle of SO$(\TT)$.
Then it follows that there exists an almost complex structure $J_{\mathcal B}$ on 
${\mathcal B}_{\J_\psi}(M)$which is given by 
$$
J_{\mathcal B}(\dot{\J}):=\J\dot{\J}, \qquad (\,\,\dot{\J}\in T_{\J}{\mathcal B}_{\J_\psi}(M) \,\,)
$$
We also have a Riemannian metric $g_{\mathcal B}$ and a $2$-form $\Ome_{\mathcal B}$ on 
${\mathcal B}_{\J_\psi}(M)$ by 
\bgn{align}\label{Apsi}
&g_{\mathcal B}(\dot{\J_1},\dot{\J_2}):=\int_M \tr(\dot{\J_1}\dot{\J_2})\,\vol_M
\notag\\
&\Ome_{\mathcal B}(\dot{\J_1},\dot{\J_2}):=\int_M \tr(\J\dot{\J_1}\dot{\J_2})
\,\vol_M
\end{align}
for $\dot{\J_1}, \dot{\J_2}\in T_{\J}{\mathcal B}_{\J_\psi}(M)$. 
Then we have
\bgn{proposition}\text{\rm (\cite{Goto_2020}  Proposition 4.2)}
$J_{\mathcal B}$ is integrable almost complex structure on ${\mathcal B}_{\J_\psi}(M)$ and
$\Ome_{{\mathcal B}}$ is a K\"ahler form on ${\mathcal B}_{\J_\psi}(M).$
\end{proposition}
Let $\wtil{\Diff}(M)$ be an extension of diffeomorphisms of $M$ by the additive group of real $2$-forms which is defined in Section \ref{Generalized Scalar curvature},
$$
\wtil{\Diff}(M):=\{\, e^b F\,:\, F\in \Diff(M),\,\, b\in \Omega^2(M, \R)\, \,\}.
$$
Note that the product in $\wtil{\Diff}(M)$ is given by 
$$
(e^{b_1}F_1)( e^{b_2}F_2) :=e^{b_1+F_1^*(b_2)}F_1\circ F_2,
$$
where $F_1, F_2\in \Diff(M)$ and $b_1, b_2\in\Omega^2(M,\R)$.
The action of $\wtil{\Diff}(M)$ on ${\mathcal B}(M)$ by 
\bgn{equation}\label{ebFcircJcircF-1}
e^{b} F_\#\circ \J\circ F_\#^{-1} e^{-b}, 
\end{equation}
 where $F\in \Diff(M)$ acts on $\J$ by $F_\#\circ \J\circ F_\#^{-1}$ 
 and $e^b$ is regarded as an element of SO$(\TT)$ and $F_\#$ denotes the bundle map of $\TT$ which is the lift of $F.$
For a generalized complex structure $\J_\psi,$
we define $\wtil{\Diff}_{\J_\psi}(M)$ to be a subgroup consisting  of elements of $\wtil{\Diff}(M)$ which preserves $\J_\psi$, 
$$
\wtil{\Diff}_{\J_\psi}(M)=\{\, e^bF\in \wtil{\Diff}(M)\, : 
e^{b} F_\#\circ \J_\psi\circ F_\#^{-1} e^{-b}=\J_\psi\, \}.
$$
We also define 
$$
\wtil{\Diff}_{\J_\psi}(M,\vol_M)= \{  e^b F\in \wtil{\Diff}_{\J_\psi}(M)\, :\,  F^*\vol_M=\vol_M\}
$$
Then from (\ref{Apsi}), we have the following,
\bgn{proposition}\label{The symplectic structure Ome B is invariant under}
The symplectic structure $\Ome_{\mathcal B}$ 
is invariant under
the action of the group $\wtil{\Diff}_{\J_\psi}(M,\vol_M)$.\end{proposition}
\bgn{proof}
The result follows from (\ref{Apsi}) and (\ref{ebFcircJcircF-1}) since $\vol_M$ is invariant under the action of $\wtil{\Diff}_{\J_\psi}(M,\vol_M)$.
\end{proof}
\bgn{proposition}\label{LHJpsidfvolMbean}
Let $L^H_{\J_\psi df}$ be an infinitesimal action of $\ham(M,\J_\psi)$. Then $L^H_{\J_\psi df}\vol_M =0$.
\end{proposition}
\bgn{proof}
Since $\vol_M =i^{-n}\lan \psi_\a, \ol\psi_\a\ran_s$, we have 
\bgn{equation}\label{L normalization condition}
L^H_{\J_\psi df}\vol_M = i^{-n}\lan L_{\J_\psi df}\psi, \ol\psi\ran_s+i^{-n}\lan \psi,  L_{\J_\psi df}\ol\psi\ran_s.
\end{equation}
From Lemma \ref{LHJPSIdfJpsi},  we have 
$$
L_{\J_\psi df}\psi =2\lan \J_\psi df , \zeta_\a\ran_{\tt}\psi_\a.
$$
Since $\lan \J_\psi df , \zeta_\a\ran_{\tt}$ is  a pure imaginary function, 
we have 
$$
L_{\J_\psi df}\ol\psi =-2\lan \J_\psi df , \zeta_\a\ran_{\tt}\ol\psi_\a.
$$
it follows from (\ref{L normalization condition}) that 
\bgn{align*}L^H_{\J_\psi df}\vol_M=&
i^{-n}\lan L_{\J_\psi df}\psi, \ol\psi\ran_s+i^{-n}\lan \psi,  L_{\J_\psi df}\ol\psi\ran_s\\
=&2\lan \J_\psi df , \zeta_\a\ran_{\tt}\vol_M -2\lan \J_\psi df , \zeta_\a\ran_{\tt}\vol_M\\=&0
\end{align*}
\end{proof}
\bgn{proposition}\label{GJpsi preserves Ome}

The action of $\ham(M, \J_\psi)$  preserves $\Ome_{\mathcal B}.$
 \end{proposition}
\bgn{proof}

The Lie algebra of $\wtil\Diff(M)$ is given by $TM\oplus \Ome^2(M)$, where 
$\Ome^2(M)$ denotes the $2$-forms on $M.$
Then the Lie bracket of $TM\oplus \Ome^2(M)$ is given by 
$$
[L_{v_1}+b_1, \,\,L_{ v_2}+b_2]=L_{[v_1, v_2]}+L_{v_1}b_2-L_{v_2}b_1
$$
Let $\iota$ be the map $ \TT
\to TM\oplus \Ome^2(M)$ which is given by 
$$
\iota(v, \t):= (v,\,\, i_vH+d\t),
$$
where $v\in TM$ and $\t\in T^*M.$
Then as shown in Lemma \ref{[LHe1, LHe_2]}, one has 
$$
L^H_{[v_1+\t_1,\,\, v_2+\t_2]_H}=[L^H_{v_1+\t_1},\,\,L^H_{v_2+\t_2}],
$$
for $v_1+\t_1, v_2+\t_2\in \TT.$
Since $\ham_0(M,\J_\psi)=\{\, \J_\psi df,|\, f\in C^\infty_\psi(M, \R)\}\subset 
\TT$, then one sees that the map $\iota$ restricted to $\ham(M,\J_\psi)$ gives 
a Lie algebra homomorphism from $\ham(M,\J_\psi)$ to the Lie algebra of $\wtil\Diff_{\J_\psi}(M)$. Its kernel consists of the locally constant functions, and its image is naturally identified with
$\ham_0(M,\J_\psi)$.
Then it follows from Proposition \ref{LHJpsidfvolMbean} that 
the image
$\ham(M, \J_\psi)$ is a Lie subalgebra of the Lie algebra of $\wtil\Diff_{\J_\psi}(M)$ preserving $\vol_M$.
Then the result follows from Proposition \ref{The symplectic structure Ome B is invariant under}.
\end{proof}
The extended action $\wtil L^H_{f,g}:= L^H_{\J_\psi df}+\sqrt{-1}gL^H_{\zeta_\a}$  defined in Definition \ref{def:lifted action} also preserves the volume form $\vol_M$.  
\bgn{proposition}
$\wtil L^H_{f,g}\vol_M=0$.
\end{proposition}
\bgn{proof}
From Proposition \ref{LHJpsidfvolMbean},
it suffices to show that $ \sqrt{-1} gL^H_{\zeta_\a}\vol_M=0$, where $g \in C^\infty_{\psi}(M,\R)$. 
From Lemma \ref{LHzetaapsiJpsi=0}, we have 
$$
\sqrt{-1} gL^H_{\zeta_\a}\psi_\a =2\sqrt{-1} g\lan \zeta_\a, \zeta_\a\ran_{\tt}\psi_\a.
$$
Since $\zeta_\a$ is pure imaginary and $g$ is real, we see that $\sqrt{-1} g\lan \zeta_\a, \zeta_\a\ran_{\tt}$ is pure imaginary.
As in proof of Proposition \ref{LHJpsidfvolMbean}, we have 
 $$\sqrt{-1} gL^H_{\zeta_\a}\vol_M = \Re\,  i^{-n}\lan \sqrt{-1} gL^H_{\zeta_\a}\psi_\a, \ol\psi_\a\ran_s =0.$$
 Then the result follows.
\end{proof}
The extension $\wtil\ham (M,\J_\psi)$ as in Section \ref{Generalized Hamiltonian diffeomorphisms} acts
on $\Ome_{\mathcal B}.$ 
Then we also have 
\bgn{proposition}
The extended action of $\wtil\ham(M, \J_\psi)$ as in Definition \ref{def:lifted action}
preserves $\Ome_{\mathcal B}.$
\end{proposition}
\bgn{proof}
The extended action $\wtil L^H_{f,g}$ is given by 
\bgn{align}
\wtil L^H_{f,g}:=&L_{\J_\psi df}^H+\sqrt{-1}gL^H_{\zeta_\a}\\
=&L^H_{\J_\psi df} +L^H_{\sqrt{-1}g\zeta_\a}-\sqrt{-1}dg\cdot\zeta_\a
\end{align}
Since $g\in C^\infty_{\psi}(M,\R)$,
the Clifford action $dg\cdot\zeta_\a=\frac12(dg\cdot\zeta_\a-\zeta_\a\cdot dg)$  preserves $\vol_M$ and $\Ome_{\mathcal B}$. 
As in the proof of  Proposition \ref{GJpsi preserves Ome}, 
$L^H_{\J_\psi df}$ and $ L^H_{\sqrt{-1}g\zeta_\a}$ preserves $\Ome_{\mathcal B}$ infinitesimally. 
Then the result follows.
\end{proof}
Recall that
the Lie algebra $\ham(M,\J_\psi)=C_\psi^\infty(M,\R),$
where $C^\infty_\psi(M,\R)=\{\, f\in C^\infty(M,\R)\, |\,\pi_T(\zeta_\a) f=0\,\}.$
Then we denote by $e$ the element $\J_\psi(df)\in \TT$, which is referred to as {\it a generalized Hamiltonian element}. 
Since $e+\sqrt{-1} df=\J_\psi df+\sqrt{-1}df\in L_{\J_\psi}$,  we have $e\cdot\psi_\a=-\sqrt{-1}df\cdot\psi_\a$. 

Let $\J_\psi$ be a $H$-twisted generalized complex structure on a compact manifold $M$ with a fixed volume form 
$\vol_M$. The generalized scalar curvature is regarded as a  function $S_{\J_\psi}$ on $\B_{\J_\psi}(M)$ by 
$$S_{\J_\psi}: \J\mapsto S(\J, \J_\psi, \vol_M).$$
For $f\in\ham(M,\J_\psi)$, we have the coupling 
$$
\lan S_{\J_\psi}(\J), f\ran =\int_M  S_{\J_\psi}(\J)f \vol_M
$$
Then $S_{\J_\psi}$ gives $\ham(M,\J_\psi)^*$-valued function on $\B_{\J_\psi}(M)$.
Then we obtain
\begin{proposition}\label{main proposition}
\begin{itemize}
\item [(1)]  $S_{\J_\psi}: \B_{\J_\psi}(M)\to\ham(M,\J_\psi)^*$ is an infinitesimally  $\ham(M,\J_\psi)$-equivariant map 
\item[(2)]   $S_{\J_\psi}$ satisfies the following 
$$
d\lan S_{\J_\psi}(\J),  f\ran (\dot \J)=\Ome_\B( L^H_{\J_\psi df} \J,  \dot{\J}) -2\sqrt{-1}
\Ome_\B ( fL^H_{\zeta_\a}J, \dot{\J})
$$
\end{itemize}
\end{proposition}
\bgn{theorem}\label{existence of moment map} 
Let $\J_\psi$ be a $H$-twisted generalized complex structure on a compact manifold $M$ equipped with the fixed volume form $\vol_M$.
We denote by $\mathfrak g$ the Lie algebra $\wtil\ham(M, \J_\psi)$ and 
by $\mathfrak g_0$ the Lie subalgebra $\ham(M,\J_\psi)$. 
Let $\rho: \mathfrak g_0 \to \mathfrak g$ be a splitting which is given by 
$\rho(f)=(f, -2f)$.
Then there exists a modified moment map associated with $\rho$
$$\mu: {\mathcal B}_{\J_\psi}(M)\to C^\infty_\psi(M,\R)^*$$ 
for the action of the Lie algebra 
$\wtil\ham(M, \J_\psi)$ such that
 $\mu(\J)$ is given by the generalized scalar curvature $S(\J, \J_\psi,\vol_M)$. More precisely,
with respect to the natural pairing, one has
$$
\langle \mu(\J),f\rangle
=
\int_M f\,S(\J,\J_\psi,\vol_M)\,\vol_M
$$
for all $f\in C^\infty_\psi(M,\mathbb R)$.
\end{theorem}
In the remainder of this section, we prove Proposition \ref{main proposition}  and Theorem \ref{existence of moment map}. 
\bgn{remark}
In the previous paper \cite{Goto_2016}, \cite{Goto_2020}, the existence of the moment map was shown in the rather restricted cases of generalized K\"ahler manifolds of symplectic type or the case of $\zeta_\a=0$. 
However, when $\pi_T(\zeta_\a)\neq 0$, the existence of an ordinary moment map may be obstructed, as discussed in Section 2. It would be interesting to investigate this obstruction further.
\end{remark}


\subsection{Preliminary results}
In order to give a proof of  Theorem \ref{existence of moment map}, we prove several Lemmas, in particular on the Nijenhuis tensor $N$. 
Let $\J\in\mathcal B_{\J_\psi}(M)$ be an almost generalized complex structure which is induced from 
a set of nondegenerate, pure spinors $\phi=\{\phi_\a\}.$
We normalize $\{\phi_\a\}$ such that $i^{-n}\lan \phi_\a,\,\,\ol\phi_\a\ran_s =\vol_M$ for each $\a.$
Then $d_H\phi_\a$ is given by  
\bgn{equation}\label{dphia=eta+N}
d_H\phi_\a=(\eta_\a+N_\a)\cdot\phi_\a,
\end{equation}
where $\eta_\a\in \sqrt{-1}(\TT)$ and $N_\a\in (\w^3\L_\J\oplus \w^3\ol\L_\J)^\R$ as before.
\bgn{remark} Note that $N_\a$ is a real element. $N_\a=N_\b$ for all $\a, \b$. 
Then $N_\a$ defines a global element $N$, which is called {\it Nijenhuis tensor.}
\end{remark}
\bgn{lemma}\label{Ncdotpsia=0}
$N\cdot\psi_\a=0$ \end{lemma}
\bgn{proof}
Since $N$ is uniquely determined by (\ref{dphia=eta+N}), 
for $e_1, e_2, e_3\in \L_\J$, we have
\bgn{align*}
N(e_1, e_2, e_3)\lan \phi_\a, \,\,\ol\phi_\a\ran_s=&
\lan d_H\phi_\a,\,\, e_1\cdot e_2\cdot e_3\cdot\ol\phi_\a\ran_s=
-\lan e_1\cdot e_2\cdot d_H\phi_\a,\,\, e_3\cdot\ol\phi_\a\ran_s\\
=&-\lan [e_1, e_2]_{H}\cdot\phi_\a,\,\,e_3\cdot\ol\phi_\a\ran_s
\end{align*}
Thus we have 
\bgn{equation}\label{N(e1,e2,e3)}
N(e_1, e_2, e_3)=2\lan [e_1, e_2]_{H}, \,\, e_3\ran_{\tt}
\end{equation}
This implies that $N=0$ if and only if $\J$ is integrable. 
By using $\J_\psi,$ we have the decomposition 
$\L_\J =\L_\J^+\oplus \L_\J^-$ and $\ol\L_\J =\ol\L_\J^+\oplus\ol\L_\J^-.$
Since $\ker \psi_\a =\L_\J^+\oplus\ol\L_\J^-$ and $N\in (\w^3\L_\J\oplus \w^3\ol\L_\J)^{\Re},$ we have 
$N\cdot\psi =(\ol N^{+}+ N^-)\cdot\psi,$
where $\ol N^+\in \w^3\ol\L_\J^+$ and $N^-\in \w^3\L_\J^-.$
From (\ref{N(e1,e2,e3)}), we see 
$$
N(e_1^-, e_2^-, e_3^-)=\lan [e_1^-, e_2^-]_{H},\,\, e_3^-\ran_{\tt}
$$
Since $\J_\psi$ is integrable, it follows that 
$[e_1^-, e_2^-]_{H}\in \ol\L_{\J_\psi}$. Since $e_3^-\in \ol\L_{\J_\psi},$ we have 
$N(e_1^-, e_2^-, e_3^-)=0.$
Thus $N^-=0$. Since $N$ is real, we have $\ol N^+=0.$
Hence $N\cdot\psi =0.$
\end{proof}
\bgn{remark}
In the cases of the ordinary K\"ahler manifolds, $\psi$ is given by the exponential 
$e^{\sqrt{-1}\ome},$ where $\ome$ denotes a symplectic form and  $\sqrt{-1}\ome=\sum_i \t^i\w\ol\t^i,$ where 
$\{\t^i\}_{i=1}^n$ denotes the local unitary frame of $\w^{1,0}$. 
We also denote by $\{v_i\}_{i=1}^n$ the dual basis of $T^{1,0}_M.$
Then the Nijenhuis tensor $N$ is given by $N=\sum N_{i,j}^k\ol\t^i\w \ol\t^j\otimes v_k
+\ol N_{i,j}^k\t^i\w \t^j\otimes \ol v_k,$ where $N_{i,j}^k\ol\t^i\w \ol\t^j$ is the $(0,2)$-component of $d\t^i$.
 The Nijenhuis tensor$N$ acts on $\ome$ by the interior and exterior product, that is, the action of Clifford algebra
Since $d\ome=0$ implies $N\cdot \ome =0.$ Thus one has $N\cdot\psi =0.$
\end{remark}
\bgn{lemma}\label{lem: N notoriatukai 3}
We have
$$
\lan e\cdot N\cdot \phi_\a,\,\,\, h\cdot\ol\phi_\a\ran_s=0
$$
\end{lemma}
\bgn{proof}
Since $h=h^{2,0}+ h^{0,2}\in \w^2 \L_\J\oplus \w^2\ol\L_\J$ and $N=N^{3,0}+N^{0,3}\in \w^3\L_\J\oplus \w^3\ol\L_\J$, we have 
\bgn{align*}\lan e\cdot N\cdot\phi_\a,\,\,\, h\cdot\ol\phi_\a\ran_s
=&-\lan  N\cdot\phi_\a,\,\,\, e\cdot h\cdot\ol\phi_\a\ran_s\\
=&-\lan N^{0,3}\cdot \phi_\a,\,\,\,e^{1,0}\cdot h^{2,0}\cdot\ol\phi_\a\ran_s 
\end{align*}
Since $e^{1,0}\cdot h^{2,0}=h^{2,0}\cdot e^{1,0}$, 
we have 
\bgn{align*}
\lan e\cdot N\cdot \phi_\a,\,\,\,h\cdot\ol\phi_\a\ran_s 
=&-\lan   N^{0,3}\cdot \phi_\a,\,\,\,h^{2,0}\cdot e^{1,0}\cdot\ol\phi_\a\ran_s \\
=&\lan h^{2,0}\cdot N^{0,3}\cdot\phi_\a,\,\,\,e^{1,0}\cdot\ol\phi_\a\ran_s
\end{align*}
We denote by $[h^{2,0}, N^{0,3}]^{0,1}\in \ol\L_\J$ the component  of $[h^{2,0}, N^{0,3}]$. 
Then we obtain
\bgn{align*}
\lan h^{2,0}\cdot N^{0,3}\cdot\phi_\a,\,\,\,e^{1,0}\cdot\ol\phi_\a\ran_s=&
\lan [h^{2,0}\cdot N^{0,3}]^{0,1}\cdot\phi_\a,\,\,\,e^{1,0}\cdot\ol\phi_\a\ran_s\\
=&2\lan  [h^{2,0},\,\, N^{0,3}]^{0,1},\,\,\, e^{1,0}\ran_{\scriptscriptstyle{T\oplus T^*}}
\lan\phi_\a,\,\,\,\ol\phi_\a\ran_s\\
=&2\lan  [h^{2,0},\,\, N^{0,3}]^{0,1},\,\,\, e^{1,0}\ran_{\scriptscriptstyle{T\oplus T^*}}
\lan\psi_\a,\,\,\,\ol\psi_\a\ran_s\\
=&\lan [h^{2,0},\,\,N^{0,3}]^{0,1}\cdot\psi_\a,\,\,\,e^{1,0}\cdot\ol\psi_\a\ran_s
\\&-\lan e^{1,0}\psi_\a,\,\,\,[h^{2,0},\,\, N^{0,3}]^{0,1}\cdot\ol\psi_\a\ran_s
\end{align*}

From Lemma \ref{Ncdotpsia=0} and $h\cdot\psi=0$, we have 
$[h, N]\cdot\psi=0$. 
Thus we have $[h^{2,0},\,\, N^{0,3}]^{0,1}\cdot\ol\psi=0$ and
$[h^{2,0},\,\, N^{0,3}]^{0,1}\cdot\psi=0$. 
Hence we obtain 
$\lan e\cdot N\cdot\phi_\a,\,\,\, h\cdot\ol\phi_\a\ran_s=0$.
\end{proof}
\bgn{lemma}\label{lem: N notoriatukai 4}

We also have
$$
\lan N\cdot e\cdot \phi_\a,\,\,\, h\cdot\ol\phi_\a\ran_s=0
$$
\end{lemma}
\bgn{proof}
The result follows as in the  proof of Lemma \ref{lem: N notoriatukai 3}

\end{proof}
Let $\J_t$ be deformations of $\J$ such that $(\J_t, \J_\psi)$ is an almost generalized K\"ahler structure. 
Then we have a family $\{\phi_\a(t)\}$ which gives rise to $\J_t$ depending smoothly on parameter $t$.
Then 
$d_H\phi_\a(t)=(\eta_\a(t)+N(t))\cdot\phi_\a(t),$
where $\eta_\a(t)\in \sqrt{-1}(\TT)$ and $N \in (\w^3\L_{\J_t}\oplus \w^3\ol\L_{\J_t})^{\Re}.$
\bgn{lemma}\label{dotN(t)cdotpsi=0}
 Let $\dot N=\frac{d}{dt}N(t)|_{t=0}$. Then we have 
$$
\dot N\cdot \psi_\a=0
$$
\end{lemma}
\bgn{proof}From Lemma \ref{Ncdotpsia=0}, we have $N(t)\cdot\psi_\a=0$ for all $t.$
Since $\psi_\a$ is fixed,
then we have the result.
\end{proof}
\bgn{lemma}\label{lanecdotphia, dot Ncdotolphiaran=0}
$\lan e\cdot\phi_\a, \,\,\dot N\cdot\ol\phi_\a\ran_s=0.$
\end{lemma}
\bgn{proof}
The space $\w^4(\TT)$ is decomposed into $\w^4 T_M \oplus (\w^3 TM \otimes TM^*)
\oplus (\w^2TM\otimes \w^2TM^*)\oplus (TM\otimes \w^3TM^*)\oplus \w^4TM^*.$
We denote by Cont$^{2,2}$ the contraction of the component $(\w^2TM\otimes \w^2TM^*)$
which yields a map from $\w^4(\TT)$ to $C^\infty(M).$
Then it follows 
\bgn{align}
\lan e\cdot\phi_\a, \,\,\dot N\cdot\ol\phi_\a\ran_s=&-\lan \phi_\a, \,\,e\cdot \dot N\cdot\ol\phi_\a\ran_s\\
=&-\text{\rm Cont}^{2,2}
(e\cdot \dot N)\lan\phi_\a, \,\,\ol\phi_\a\ran_s
\end{align}
Since $\lan \phi_\a, \,\,\ol\phi_\a\ran_s=\lan \psi_\a,\,\,\ol\psi_\a\ran_s,$
we have 
\bgn{align}
\lan e\cdot\phi_\a, \,\,\dot N\cdot\ol\phi_\a\ran_s=&-\text{\rm Cont}^{2,2}
(e\cdot \dot N)\lan\psi_\a, \,\,\ol\psi_\a\ran_s\\
=&\lan e\cdot\psi_\a, \,\, \dot N\cdot \ol\psi_\a\ran_s
\end{align}
Since $\dot N$ is real, it follows from Lemma \ref{dotN(t)cdotpsi=0} that 
$\dot N\cdot\ol\psi_\a=0.$
Hence we have
$\lan e\cdot\phi_\a, \,\,\dot N\cdot\ol\phi_\a\ran_s=0.$
\end{proof}

\subsection{Proof of the generalized scalar curvature as a modified moment map}
This subsection is devoted to proving Proposition \ref{main proposition} and  our main theorem, i.e., Theorem  \ref{existence of moment map}.
\bgn{proof}[Proof of Proposition \ref{main proposition}]
Let $\J_\phi$ be an almost generalized complex structure which is induced from $\phi=\{\phi_\a\}$ as before.
We denote by $\J_t$ small deformations of $\J_\phi.$
Any infinitesimal deformation $\dot\J:=\frac{d}{dt}\J_t|_{t=0}$ of $\J$ is given by the adjoint action of $h$ 
$$
\dot\J_h :=[h,\,\J]
$$
where $h$ denotes a real $C^\infty$ global section of $ (\w^2\L_{\J}\oplus\w^2\ol\L_{\J})^\R.$
Then the corresponding infinitesimal deformation of $\phi$ is given by the Clifford action of $h$ on each $\phi_\a,$
$$
\dot{\phi}_\a=h\cdot{\phi}_\a
$$
Note that we can assume $\dot\phi_\a$ is transversal to the canonical line bundle $K_{\J}$ in the following sense, 
$$\lan \phi_\a, \ol{\dot{\phi}_\a}\ran_s=0.$$
A Hamiltonian element $e=\J_\psi df$ gives the infinitesimal deformation $L^H_e\J$ of $\J.$
Since $d_H\phi_\a=\eta_\a\cdot\phi_\a+N\cdot\phi_\a,$
then the corresponding infinitesimal deformations of $\phi_\a$ is given by $L^H_e\phi_\a,$ 
Then it follows
$$
L^H_e\phi_\a=d_H e\cdot\phi_\a+e\cdot d_H\phi_\a=d_H e\cdot\phi_\a+e\cdot (\eta_\a+N)\cdot\phi_\a,
$$
where $d_H$ denotes $d+H.$
Recall a modified action $\wtil L^H_{f,g}$ (see  Definition \ref{def:lifted action}), 
$$
\wtil L^H_{f,g}:= L^H_{e}+\sqrt{-1}gL^H_{\zeta_\a},
$$
Then in order to show the existence of a modified moment map associated with $\rho(f)=(f, -2f)$, we shall calculate 
$\Ome_\B(\wtil L^H_{f,-2f}\J, \,\, \dot\J_h)$.
Recall the formula (see Section 7 in \cite{Goto_2016}
). 
\bgn{equation}\label{OmeB(dotJh1,dotJh2)}\Ome_\B(\dot\J_{h_1}, \,\,\dot\J_{h_2})=
 \Im\( i^{-n}\int_M \lan h_1\cdot\phi_\a, \,\,h_2\cdot\ol\phi_\a \ran_s\),
 \end{equation}
Applying (\ref{OmeB(dotJh1,dotJh2)}), we obtain 
\bgn{align}
\Ome_\B(\til L^H_{f,-2f}\J, \,\, \dot\J_h)=& \Im\(i^{-n}\int_M \lan\til L^H_{f,-2f}\phi_\a, \,\, h\cdot\ol\phi_\a\ran_s\)\\
=&\Im\(i^{-n}\int_M\lan d_H e\cdot\phi_\a+e\cdot (\eta_\a+N)\cdot\phi_\a, \,\, h\cdot\ol\phi_\a\ran_s\)
\label{First we calculate}\\
-&\Im \(i^{-n}\int_M \lan 2\sqrt{-1}fL^H_{\zeta_\a}\phi_\a, \,\,h\cdot\ol\phi_\a\ran_s\)
\end{align}
First we shall calculate the term  (\ref{First we calculate}), that is,  $\Ome_\B(L^H_e\J, \,\, \dot\J_h)$.
Since $h\in(\w^2\L_{\J}\oplus\w^2\ol\L_{\J})^\R, $
we have 
$\lan \phi_\a, \,\, h\cdot\ol\phi_\a\ran=0.$
Since $e\cdot\eta_\a+\eta_\a\cdot e =2\lan e, \,\,\eta_\a\ran_{\tt},$ we have  
$$
\lan e\cdot\eta_\a\cdot\phi_\a,\,\, h\cdot\ol\phi_\a\ran_s =-\lan \eta_\a\cdot e\cdot\phi_\a, \,\,h\cdot\ol\phi_\a\ran_s
$$
Applying Lemma \ref{lem: N notoriatukai 3}, we have 
\bgn{align}\label{OmeBLHeJdotJh=cnIM}
\Ome_\B(L^H_e\J, \,\, \dot\J_h)=&\Im\(i^{-n}\int_M\lan (d_H-\eta_\a) e\cdot\phi_\a, \,\, h\cdot\ol\phi_\a\ran_s\)
\end{align}
Since $d\lan (e\cdot\phi_\a),\,\,h\cdot\ol\phi_\a\ran_{[2n-1]}=\lan d_H(e\cdot\phi_\a),\,\,h\cdot\ol\phi_\a\ran_s
-\lan (e\cdot\phi_\a),\,\,d_H(h\cdot\ol\phi_\a)\ran_s$ and 
$\lan (e\cdot\phi_\a),\,\,h\cdot\ol\phi_\a\ran_{[2n-1]}$ gives a globally defined $(2n-1)$-form on $M,$
then applying the Stokes Theorem, we have 
$$
\int_M \lan d_H(e\cdot\phi_\a),\,\,h\cdot\ol\phi_\a\ran_s=\int_M \lan e\cdot\phi_\a, \,\,d_H(h\cdot\ol\phi_\a)\ran_s,
$$
where we are also applying $\lan H\cdot(e\cdot\phi_\a),\,\,h\cdot\ol\phi_\a\ran_s= \lan e\cdot\phi_\a, \,\,H\cdot(h\cdot\ol\phi_\a)\ran_s.
$
Since $\eta_\a$ is in $\sqrt{-1}(\TT)$, it follows $\ol{\eta}_\a=-\eta_\a$. Then we have 
$$
\lan \eta_\a\cdot e\cdot \phi_\a, \,\, h\cdot\ol\phi_\a\ran_s=-\lan e\cdot\phi_\a, \,\, \eta_\a\cdot h\cdot\ol\phi_\a\ran_s.
$$
Substituting them into (\ref{OmeBLHeJdotJh=cnIM}), we obtain 
$$
\Ome_\B(L^H_e\J, \,\, \dot\J_h)=\Im\(i^{-n}\int_M\lan e\cdot\phi_\a, \,\, (d_H+\eta_\a)\cdot (h\cdot\ol\phi_\a)\ran_s\)
$$
Let $\phi(t)=\{\phi_\a(t)\}$ be a one parameter family of nondegenerate, pure spinors which gives 
$$
\frac{d}{dt}\phi_\a(t)|_{t=0}=h\cdot\phi_\a.
$$
We assume that $\phi_\a(t)$ depends smoothly on $t$ .
Then $d_H\phi_\a(t)$ is given by
$$
d_H\phi_\a(t)=(\eta_\a(t)+N(t))\cdot\phi_\a(t).
$$
Taking the differential of both sides at $t=0$, we have 
$$
d_H(h\cdot\phi_\a) =(\dot\eta_\a+\dot N)\cdot\phi_\a+(\eta_\a+N)\cdot (h\cdot\phi_\a),
$$
where $\dot{\eta}_\a=\frac{d}{dt}\eta_\a|_{t=0}$ and $\dot N=\frac{d}{dt}N(t)|_{t=0}.$
Since $\eta_\a(t)$ is pure imaginary and $N(t)$ is real, we have 
$$
d_H(h\cdot\ol\phi_\a) =(-\dot\eta_\a+\dot N)\cdot\ol\phi_\a+(-\eta_\a+N)\cdot (h\cdot\ol\phi_\a)
$$
Thus we obtain 
$$
\lan e\cdot\phi_\a, \,\, (d_H+\eta_\a-N)\cdot (h\cdot\ol\phi_\a)\ran_s =-\lan e\cdot\phi_\a, \,\,(\dot\eta_\a-\dot N)\cdot\ol\phi_\a\ran_s
$$
From Lemma \ref{lanecdotphia, dot Ncdotolphiaran=0}, we have 
$\lan e\cdot\phi_\a, \,\,\dot N\cdot\ol\phi_\a\ran_s=0.$
Applying Lemma \ref {lem: N notoriatukai 4}, we obtain 
$$
\lan e\cdot\phi_\a, \,\, (d_H+\eta_\a)\cdot (h\cdot\ol\phi_\a)\ran_s =-\lan e\cdot\phi_\a, \,\,\dot\eta_\a\cdot\ol\phi_\a\ran_s
$$

We decompose $e$ as $e^{1,0}+e^{0,1},$ where $e^{1,0}\in \L_\J$ and $e^{0,1}\in \ol\L_\J.$
We also decompose $\dot\eta_\a=\dot\eta_\a^{1,0}+\dot\eta_\a^{0,1},$ where $\dot\eta_\a^{1,0}\in \L_\J$ and $\dot\eta_\a^{0,1}\in \ol\L_\J.$

Then we have 
\bgn{align*}
-\Im \(i^{-n}\lan e\cdot\phi_\a, \,\, \dot\eta_\a\cdot\ol\phi_\a\ran_s \)=&
-\Im\(i^{-n}\lan e^{0,1}\cdot\phi_\a, \,\, \dot\eta_\a^{1,0}\cdot\ol\phi_\a\ran_s \)\\
=&\Im\( i^{-n}\lan \dot\eta_\a^{0,1}\cdot e^{1,0}\cdot\phi_\a, \,\,\ol\phi_\a\ran_s\)\\
=&\Im\(i^{-n}2\lan e^{0,1}, \,\,\dot\eta_\a^{1,0}\ran_{\tt}\lan \phi_\a, \,\,\ol\phi_\a\ran_s\)
\end{align*}
Since $e$ is real and $\eta_\a$ is pure imaginary, we have 
\bgn{align*}
\lan\dot\eta_\a, \,\,e\ran_{\tt}=&\lan \dot\eta_\a^{1,0}, \,\, e^{1,0}\ran_{\tt}+\lan \dot\eta_\a^{0,1}, \,\,
e^{0,1}\ran_{\tt}\\
=&\lan \dot\eta_\a^{1,0}, \,\, e^{1,0}\ran_{\tt}+\lan -\ol{\dot\eta_\a^{1,0}}, \,\,\ol{e^{1,0}}\ran_{\tt}\\
=&\lan \dot\eta_\a^{1,0}, \,\, e^{1,0}\ran_{\tt}-
\ol{\lan \dot\eta_\a^{1,0}, \,\, e^{1,0}\ran}_{\tt}\\
=&2\sqrt{-1}\,\Im\,\lan \dot\eta_\a^{1,0}, \,\, e^{1,0}\ran_{\tt}
\end{align*}
Since $i^{-n}\lan\phi_\a,\,\,\ol\phi_\a\ran_s =i^{-n}\lan \psi, \ol\psi\ran_s=\vol_M,$
we have 
$$
-\Im \(i^{-n}\lan e\cdot\phi_\a, \,\, \dot\eta_\a\cdot\ol\phi_\a\ran_s \)=\Im\,\lan\dot\eta_\a, \,\,e\ran_{\tt}\vol_M
$$
Since $i^{-n}\lan\psi_\a, \,\,\ol\psi_\a\ran_s =\vol_M,$ we have
$$
\Im\, \lan\dot\eta_\a, \,\,e\ran_{\tt}\,\vol_M=\Im\,\lan\dot\eta_\a, \,\,e\ran_{\tt}i^{-n}\lan \psi_\a\, \,\,\ol\psi_\a\ran_s
$$
Then as before, we obtain 
\bgn{align*}
\Im\, \lan\dot\eta_\a, \,\,e\ran_{\tt}\,\vol_M=
&\Im\,\lan\dot\eta_\a, \,\,e\ran_{\tt}i^{-n}\lan \psi_\a\, \,\,\ol\psi_\a\ran_s\\
=&-\Im\, \(i^{-n}\lan e\cdot\psi_\a, \,\,\dot\eta_\a\cdot\ol\psi_\a\ran_s\)
\end{align*}
Since $e$ is a generalized Hamiltonian element,
we have $e\cdot\psi_\a=-\sqrt{-1}df\cdot\psi_\a$. Applying $d_H\psi_\a=\zeta_\a\cdot\psi_\a,$ we have 
\bgn{align*}
e\cdot\psi_\a=&(-\sqrt{-1}df)\cdot\psi_\a=-\sqrt{-1}\(d_H(f\psi_\a)-f\zeta_\a\cdot\psi_\a\)\\
=&-\sqrt{-1}(d_H-\zeta_\a)(f\psi_\a)
\end{align*}
Then we obtain 
\bgn{align}
\Ome_\B(L^H_e\J, \,\, \dot\J_h)=&-\Im\(i^{-n}\int_M \lan e\cdot\psi_\a, \,\,\dot\eta_\a\cdot\ol\psi_\a\ran_s\)\\
=&\Im\,\(\int_M   i^{-n+1}\lan (d_H-\zeta_\a)f\psi_\a,  \,\,\dot\eta_\a\cdot\ol\psi_\a\ran_s \)
\end{align}
Since infinitesimal deformations are given by the action of a global section $h,$
it follows that $\dot\eta_\a=\dot\eta_\b.$
Thus $\lan \psi_\a, \,\, \dot\eta_\a\cdot\ol\psi_\a\ran_{[2n-1]}$ defines a global $(2n-1)$-form on $M.$
We have 
$$d\lan f\psi_\a, \,\, \dot\eta_\a\cdot\ol\psi_\a\ran_{[2n-1]}=\lan d_H (f\psi_\a), \,\, \dot\eta_\a\cdot\ol\psi_\a\ran_s
+\lan f\psi_\a, \,\, d_H(\dot\eta_\a\cdot\ol\psi_\a)\ran_s
$$
Applying the Stokes theorem again, the first term (\ref{First we calculate}) is given by 
 \bgn{align}
 \Ome_\B(L^H_e\J, \,\, \dot\J_h)=
 &\Im\(i^{-n+1}\int_M\lan f\psi_\a, \,\, (d_H+\zeta_\a)\cdot(\dot\eta_\a\cdot\ol\psi_\a)\ran_s\)
 \end{align}
Then $\Ome_\B(\wtil L^H_{f,g}\J, \,\, \dot\J_h)$ is written as
\bgn{align*}
\Ome_\B(\wtil L^H_{f,g}\J, \,\, \dot\J_h)
=&\Im\int_M\, i^{-n}
\lan \sqrt{-1}(f\psi_\a), \,\,(d_H+\zeta_\a)(\dot{\eta}_\a\cdot\ol\psi_\a)\ran_s\\
&+\Im \int_Mi^{-n}\lan\, \sqrt{-1}gL^H_{\zeta_\a}\phi_\a,\,\,\ol{\dot\phi}_\a\ran_s 
\\
=&\Re\int_M\, i^{-n}
\lan (f\psi_\a), \,\,(d_H+\zeta_\a)(\dot{\eta}_\a\cdot\ol\psi_\a)\ran_s\\
&+\,\Re \int_Mi^{-n}\lan gL^H_{\zeta_\a}\phi_\a,\,\,\ol{\dot\phi}_\a\ran_s, 
\end{align*}
where we are changing from the imaginary part to the real part.
\\
Since one has 
\bgn{align*}
2\Re\,i^{-n}\lan L^H_{\zeta_\a}\phi_\a,\,\,\ol{\dot\phi}_\a\ran_s=&i^{-n}\lan L^H_{\zeta_\a}\phi_\a,\,\,\ol{\dot\phi}_\a\ran_s+
\ol{i^{-n} \lan L^H_{\zeta_\a}\phi_\a,\,\,\ol{\dot\phi}_\a\ran_s}\\
=&i^{-n}\lan L^H_{\zeta_\a}\phi_\a,\,\,\ol{\dot\phi}_\a\ran_s
-{(-i)^{-n} \lan L^H_{\zeta_\a}\ol\phi_\a,\,\,{\dot\phi}_\a\ran_s}\\
=&i^{-n}\lan L^H_{\zeta_\a}\phi_\a,\,\,\ol{\dot\phi}_\a\ran_s
-i^{-n} \lan {\dot\phi}_\a, \,\, L^H_{\zeta_\a}\ol\phi_\a\ran_s,
\end{align*}
then we have
\bgn{align*}
2\frac{d}{dt}\Big|_{t=0}\,\Re \,i^{-n}\lan L^H_{\zeta_\a}\phi_\a(t),\,\,\ol{\phi}_\a(t)\ran_s
=&i^{-n}\lan L^H_{\zeta_\a}\dot\phi_\a,\,\,\ol{\phi}_\a\ran_s
-i^{-n} \lan {\dot\phi}_\a, \,\, L^H_{\zeta_\a}\ol\phi_\a\ran_s\\
+&i^{-n}\lan L^H_{\zeta_\a}\phi_\a,\,\,\ol{\dot\phi}_\a\ran_s
-i^{-n} \lan {\phi}_\a, \,\, L^H_{\zeta_\a}\ol{\dot\phi}_\a\ran_s
\end{align*}

Since $\lan \phi_\a, \,\,\,\ol{\dot\phi}_\a\ran_s=0, $
one has
$$
0=L^H_{\zeta_\a}\lan \phi_\a, \,\,\ol{\dot\phi}_\a\ran_s
=\lan L^H_{\zeta_\a}\phi_\a, \,\,\ol{\dot\phi}_\a\ran_s
+\lan \phi_\a, \,\,L^H_{\zeta_\a}\ol{\dot\phi}_\a\ran_s
$$
Taking the complex conjugate, we have also
$$
\lan L^H_{\zeta_\a}\ol\phi_\a, \,\,{\dot\phi}_\a\ran_s
+\lan \ol\phi_\a, \,\,L^H_{\zeta_\a}{\dot\phi}_\a\ran_s=0
$$
Thus we have
\bgn{align*}
\frac{d}{dt}\Big |_{t=0}\Re\,i^{-n} \lan L^H_{\zeta_\a}\phi_\a(t),\,\,\ol{\phi}_\a(t)\ran_s=
&-i^{-n} \lan {\dot\phi}_\a, \,\, L^H_{\zeta_\a}\ol\phi_\a\ran_s
+i^{-n}\lan L^H_{\zeta_\a}\phi_\a,\,\,\ol{\dot\phi}_\a\ran_s\\
=&2\, \Re\,i^{-n} \lan L^H_{\zeta_\a}\phi_\a,\,\,\ol{\dot\phi}_\a\ran_s 
\end{align*}
Since our infinitesimal deformations of $\J_\phi$ are fixing $\J_\psi,$
both $\psi_\a$ and $\zeta_\a$ do not change.
Then we obtain 
\bgn{align*}
\frac{d}{dt}\Big |_{t=0}&\int_M\Re\(i^{-n}
\lan (f\psi_\a), \,\,(d_H+\zeta_\a)({\eta}_\a\cdot\ol\psi_\a)\ran_s
-\,i^{-n}\lan fL^H_{\zeta_\a}\phi_\a,\,\,\ol{\phi}_\a\ran_s \)\\
=&\int_M\Re\(i^{-n}
\lan (f\psi_\a), \,\,(d_H+\zeta_\a)(\dot{\eta}_\a\cdot\ol\psi_\a)\ran_s
-2\,i^{-n}\lan fL^H_{\zeta_\a}\phi_\a,\,\,\ol{\dot\phi}_\a\ran_s \)\\
=&\Ome_\B( L^H_{\J_\psi df}\J, \,\, \dot\J_h)-\int_M\Im \, i^{-n}\lan 2\sqrt{-1}fL^H_{\zeta_\a}\phi_\a,\,\,\ol{\dot\phi}_\a\ran_s \\
=&\Ome_\B(\wtil L^H_{f, -2f}\J, \,\, \dot\J_h)
\end{align*}
Then we have
\bgn{align*}
\Re\,i^{-n}\lan &(f\psi_\a), \,\,(d_H+\zeta_\a)({\eta}_\a\cdot\ol\psi_\a)\ran_s
-\Re\,i^{-n}\lan fL^H_{\zeta_\a}\phi_\a,\,\,\ol{\phi}_\a\ran_s \\
=&\Re\, i^{-n} f\lan \psi_\a, \,\, L^H_{\eta_\a}\ol\psi_\a\ran_s+
\Re\, i^{-n}f\lan\phi_\a, \,\,L^H_{\zeta_\a}\ol\phi_\a\ran_s
+2\lan \zeta_\a, \,\,\eta_\a\ran_{\tt}f\vol_M\\
=&fS(\J_\phi, \J_\psi)\vol_M
\end{align*}
Thus we obtain
\bgn{equation}\label{OmewtilJJpsi,dotJ}
\Ome_\B(\wtil L^H_{f,-2f}\J_\phi, \,\,\,\dot\J)
=\frac{d}{dt}|_{t=0}\int_M fS(\J_\phi(t), \J_\psi)\vol_M
\end{equation}
Let 
$\mu: \mathcal B_{\J_\psi}(M)\to C^\infty_\psi(M)^*$ be a map which is given by the scalar curvature 
$$
\mu(\J)=S(\J, \J_\psi),
$$
where $C^\infty_\psi(M)^*$ is identified with $C^\infty_\psi(M)$ by using the $L^2$-metric.
Then $\mu$ satisfies 
$$d\lan \mu, f\ran(\dot\J)=\Ome_{\B}(\wtil L^H_{f, -2f}\J_\phi, \,\,\dot\J).$$
\\
Finally we show that $\mu: \mathcal B_{\J_\psi}(M)\to C^\infty_\psi(M)^*$ satisfies the equivariant  condition of moment map under the action of $\ham(M,\J_\psi)$. 
Since the infinitesimal action of
$\ham(M,\J_\psi)$ defines a Lie algebra homomorphism into the Lie algebra of volume preserving 
generalized diffeomorphisms,
 it follows that the generalized scalar curvature satisfies the equivariant condition under the action of $\ham(M,\J_\psi)$. 
Thus the map $\mu$ satisfies  the infinitesimal equivariance of moment map under the action of $\ham(M,\J_\psi)$. 
\end{proof}
\bgn{proof}[Proof of Theorem \ref{existence of moment map}]
Let $\rho : \mathfrak g_0 \to \mathfrak g$ be a map given by 
$\rho(f)= (f, -2f)$. 
The result follows from Proposition \ref{main proposition}. 
\end{proof}
Apostolov, Streets, and Ustinovskiy introduced the notion of an adapted volume and, using the variational formula established in the first version of this paper [arXiv:2105.13654], proved that the generalized scalar curvature defines an ordinary moment map on the adapted-volume locus \cite{ASU_2026}. In the present framework, this locus coincides with the fixed-point locus of the abelian action.
\bgn{corollary}[cf. \cite{ASU_2026} Theorem 1.1]
Let $\B_{\J_\psi}^{\mathfrak a}(M)$ be the fixed-point set which is given by 
$$
\B_{\J_\psi}^{\mathfrak a}(M)=\{ \J \in \B_{\J_\psi}(M)\, |\, L^H_{\zeta_\a}\J =0\, \}
$$
We assume that $\B^{\mathfrak a}_{\J_\psi}(M)$ is a symplectic submanifold of $\B_{\J_\psi}(M)$.
Then the restriction of the modified moment map $\mu$ as in Theorem \ref{existence of moment map} to $\B_{\J_\psi}^{\mathfrak a}(M)$ gives an ordinary moment map for the action of $\ham(M,\J_\psi)$.
\end{corollary}
\begin{proof}
Since the constant function $1$ is an element of the abelian part $C^\infty_\psi(M,\R)$,  
the result follows from Proposition \ref{X^A}.
\end{proof}
\begin{proposition}
If the fixed point set $\B_{\J_\psi}^{\mathfrak a}$ is not empty, then there exists an ordinary moment map 
for the extended action of $\wtil\ham(M,\J_\psi)$ on $\B_{\J_\psi}(M)$. 
\end{proposition}
\begin{proof}
Since $\B_{\J_\psi}(M)$ is contractible, we have $H^1(\B_{\J_\psi}(M))=0$. 
Thus the result follows from Proposition \ref{vanishing theorem}.
\end{proof}
\section{Generalized scalar curvature of generalized K\"ahler structures on the standard Hopf surfaces}
\label{Generalized scalar curvature of generalized Kahler structures on the standard Hopf surfaces}

Two kinds of generalized K\"ahler structures are known on the standard Hopf surfaces \cite{Gua_2003}. 
One is odd type, which is a generalized K\"ahler structure consisting of 
a pair of generalized complex structures of odd type. 
The other is of even type. 
Our explicit descriptions of nondegenerate, pure spinors of generalized K\"ahler structures  enable us to calculate their generalized scalar curvature.
Moreover, an action of the Pin group converts the odd one to the even one. 
First we discuss the generalized K\"ahler structure of odd type.
\subsection{Generalized K\"ahler structures of odd type on the standard Hopf surfaces}\label{Hopf odd}
Let $(z_1, z_2)$ be complex coordinates of $\C^2$. 
We define two differential forms on $\C^2\bsh \{0\}$ by
$$
\phi=dz_1\w e^{-\sqrt{-1}\frac{\ome}{r^2}}, \qquad 
\psi=dz_2\w e^{-\sqrt{-1}\frac{\ome}{r^2}},
$$
where 
$$
\ome=\sqrt{-1}\(dz_1\w d\ol{z_1}+dz_2\w d\ol{z_2}\),\qquad 
r^2=|z_1|^2+|z_2|^2,  
$$
Then $\phi$ and $\psi$ are nondegenerate, pure spinors which give rise to almost generalized complex structures 
$\J_\phi, \J_\psi$, respectively. 
In fact, $\J_\phi$ comes from the ordinary complex structure along   $z_1$-direction and 
$\J_\phi$ is induced from 
the real symplectic structure $\exp(\frac{dz_2}r\w \frac{d\ol z_2}r)$ along $z_2$-direction. 
 Thus one has
\bgn{align*}
&\J_\phi(r\frac{\pa}{\pa z_2})=\sqrt{-1}\frac{d\ol z_2}r,\qquad  
\J_\phi(\frac {d\ol z_2}r)=\sqrt{-1}r\frac{\pa}{\pa z_2}
\end{align*}
On the other hand, 
$\J_\psi$ comes from the ordinary complex structure along  $z_2$-direction and 
$\J_\psi$ is induced from the real symplectic structure $\exp(\frac{dz_1}r\w \frac{d\ol z_1}r)$  along $z_1$-direction
Thus one also has
\bgn{align*}
&\J_\psi(r\frac{\pa}{\pa z_1})=\sqrt{-1}\frac{d\ol z_1}r,\qquad  
\J_\psi(\frac {d\ol z_1}r)=\sqrt{-1}r\frac{\pa}{\pa z_1}
\end{align*}
We denote by $\vol$ the volume form $-\frac 1{2! }\frac{\ome\w\ome}{r^4}$. 
The radial vector field $r\frac{\pa}{\pa r}$ and the logarithmic $1$-form $\frac{dr}r$are given by
\bgn{align*}
\frac{dr}r =&\frac{d r^2}{2r^2}=\frac1{2r^2}(z_1d\ol z_1+\ol z_1 dz_1+z_2d\ol z_2+ \ol z_2 d z_2)\\
r\frac{\pa}{\pa r}=&z_1\frac{\pa}{\pa z_1}+\ol z_1\frac{\pa}{\pa \ol z_1}+
z_2\frac{\pa}{\pa z_2}+\ol z_2\frac{\pa}{\pa \ol z_2}
\end{align*}
We define a real $3$-form $H$ by 
\bgn{align*}
 H:=-i_{r\frac{\pa}{\pa r}}(\vol)=&-\frac1{r^4}( z_1 d\ol z_1 \w dz_2\w d\ol z_2- \ol z_1 d z_1\w dz_2\w d\ol z_2\\&+z_2 dz_1\w d\ol z_1 \w d\ol z_2-\ol z_2 dz_1\w d\ol z_1\w dz_2)
\end{align*}
 Then one has 
 $$
d\phi+H\w \phi=0, \qquad  d\psi+ H\w\psi=0.
$$
 Then it follows that the pair $(\J_\phi, \J_\psi)$ gives a $H$-twisted generalized K\"ahler structure on $\C^2\bsh\{0\}$.
We see that $(\J_\phi, \J_\psi)$ gives the bihermitian structure $(g, I^+, I^-)$, 
where $g$ is $\frac1{r^2} g_{eu}$ and $I^+$ is the standard complex structure on $\C^2$ and $I^-$ has 
holomorphic coordinates $(z_1, \ol z_2)$, where $g_{eu}$ denotes the Euclidean metric on $\C^2\cong \R^4$. 
We consider an action of $\mathbb Z$ on $\C^2$ which is generated by the diagonal
 multiplication of $\a$, $(|\a|\neq 1)$, that is , 
$(z_1,z_2)\mapsto (\a z_1, \a z_2)$. 
Then the quotient $X:=\C^2\bsh\{0\}/\mathbb Z$ is called the standard Hopf surface.
Since $\phi$, $\psi$ are equivariant  under the action of $\a$, it follows that 
the generalized K\"ahler structure descends to the underlying differential manifold of the standard Hopf surface $X$.

In order to normalize $\phi$ and $\psi$, we replace $(\phi, \psi)$ with 
$$
(\til\phi, \til\psi):=
(\frac1{\sqrt 2}\frac{\phi}{ r}, \,\, \frac1{\sqrt 2}\frac{\psi}r).
$$
Then we see 
$$
 i^{-2}\lan \frac1{\sqrt 2}\frac{\phi}r, \,\, \frac1{\sqrt 2}\frac{\ol\phi}r\ran_s=i^{-2}\lan \frac1{\sqrt 2}\frac{\psi}r, \,\,\frac1{\sqrt 2}\frac{\ol\psi}r\ran_s =\frac1{r^4}(dz_1\w d \ol z_1\w d z_2\w d\ol z_2)
 =\vol
 $$
Since $ \sqrt{-1}\J_\phi\frac{dr}r+ \frac{dr}r\in L_{\J_\phi}=\ker\phi$, one has
 \bgn{align*}
 d_H(\frac{\phi}r)=&-\frac{dr}r \frac{\phi}r+(d_H\phi)\frac1r
 =(\sqrt{-1}\J_\phi \frac{dr}r) \frac{\phi}r, \qquad 
 \end{align*}
We also have
$$
d_H(\frac{\psi}r)=(\sqrt{-1}\J_\psi\frac{dr}r)\frac{\psi}r
$$
Since $\eta$ and $\zeta$ should be pure imaginary as in Definition of the scalar curvature, we obtain
\bgn{equation}\label{etasqrt-1Jphifrac}
\eta =(\sqrt{-1}\J_\phi \frac{dr}r) ,  \quad\zeta=(\sqrt{-1}\J_\psi\frac{dr}r).
\end{equation}
The following lemma is a key to calculate the scalar curvature of $(\J_\phi, \J_\psi)$, 
\bgn{lemma}\label{key lemma 1}
$$-\J_\psi\J_\phi\frac{dr}r=r\frac{\pa}{\pa r}
$$
\end{lemma}
\bgn{proof}
In fact, we have
\bgn{align*}
-\J_\psi\J_\phi\frac{dr}r=&-\J_\psi\(\frac{\sqrt{-1}}{r^2}(z_1 d\ol z_1 -\ol z_1 d z_1)+\sqrt{-1}\frac1r(z_2 r\frac{\pa}{\pa z_2}
-\ol z_2 r\frac{\pa}{\pa \ol z_2})\)\\
=&z_1\frac{\pa}{\pa z_1}+\ol z_1\frac{\pa}{\pa \ol z_1}+z_2\frac{\pa}{\pa z_2}+\ol z_2\frac{\pa}{\pa\ol z_2}\\
=&r\frac{\pa}{\pa r}
\end{align*}
Then the result follows.
\end{proof}
\bgn{lemma}\label{dHetacdotolpsi=}
One has 
 \bgn{align}
 \eta\cdot\ol\psi=&-(\J_\psi \J_\phi\frac{dr}r)\cdot\ol\psi=(r\frac{\pa}{\pa r})\cdot \ol\psi
 \end{align}
\end{lemma}
\bgn{proof}
Since  $(\eta-\sqrt{-1}\J_\psi \eta)\in \ol L_\psi$,
it follows
$(\eta-\sqrt{-1}\J_\psi \eta)\cdot\ol\psi =0$. 
From (\ref{etasqrt-1Jphifrac}) and Lemma \ref{key lemma 1}, we have
 $\eta\cdot\ol\psi=\sqrt{-1}\J_\psi\eta\cdot\ol\psi =-(\J_\psi \J_\phi\frac{dr}r)\cdot\psi.$
 Then the result follows.
\end{proof}
\bgn{proposition} 
 The generalized scalar curvature  of $(\J_\phi, \J_\psi)$ is a constant, that is, 
 $$
 S(\J_\phi, \J_\psi)=1
 $$
\end{proposition}
\bgn{proof}
The generalized scalar curvature $S(\J_\phi, \J_\psi)$ is given by the real part of the following:
\bgn{align}
 S(\J_\phi, \J_\psi)\vol_M=
 \Re\(i^{-n}\lan \til\phi, \,\,\,d_H(\zeta\cdot\ol{\til\phi})\ran_s
 +i^{-n}\lan \til\psi, \,\,\,d_H(\eta\cdot\ol{\til\psi})\ran_s\)
\end{align}
Since $r\frac{\pa}{\pa r}$ is the vector field defined from the dilation action, one has
 $$L_{r\frac{\pa}{\pa r}}\frac{\ol\psi}r=0.$$
 One also has  $i_{r\frac{\pa}{\pa r}}H =-i_{r\frac{\pa}{\pa r}}i_{r\frac{\pa}{\pa r}}\vol=0$.
 Then it follows 
 $$d_H (r\frac{\pa}{\pa r})\cdot =L^H_{r\frac{\pa}{\pa r}}- (r\frac{\pa}{\pa r}) \cdot d_H
 =L_{r\frac{\pa}{\pa r}}- (r\frac{\pa}{\pa r})\cdot d_H.$$
 From Lemma \ref{dHetacdotolpsi=}, one has 
 \bgn{align*}
 d_H(\eta\cdot\frac{\ol\psi} r)=&d_H( r\frac{\pa}{\pa r}\cdot\frac{\ol\psi} r)\\
 =-&(r\frac{\pa}{\pa r})\cdot d_H\frac{\ol\psi} r=(r\frac{\pa}{\pa r})\cdot(\frac{dr}r)\cdot\frac{\ol\psi} r
 \end{align*}
 Thus one has
 \bgn{align*}
 i^{-2}\lan \frac{\psi}r, \,\,\,  d_H(\eta\cdot\frac{\ol\psi} r)\ran_s=&i^{-2}\lan 
 \frac{\psi}r, \,\,\, (r\frac{\pa}{\pa r})\cdot(\frac{dr}r)\cdot\frac{\ol\psi} r
\ran_s
 \end{align*}
 Taking the real part,  we have the coupling $(r\frac{\pa}{\pa r})$ and $(\frac{dr}r)$, 
 \bgn{align*}
\Re\( i^{-2}\lan \frac\psi r, \,\,\,  d_H(\eta\cdot\frac{\ol\psi} r)\ran_s\)=&\frac12i^{-2}\lan \frac\psi r, \,\,\,\frac{\ol\psi}r\ran_s
 \end{align*}
 Thus we have
 $$
 \Re\( i^{-2}\lan \til\psi, \,\,\,d_H(\eta\cdot\ol{\til\psi})\ran_s\)=\frac12 \vol_M
 $$
 We also have 
 $$
  \Re\(i^{-2}\lan \til\phi, \,\,\,d_H(\zeta\cdot\ol{\til\phi})\ran_s\)=\frac12 \vol_M
  $$
 Thus we obtain 
 $S(\J_\phi, \J_\psi)=1$. 
 \end{proof}
 
 \subsection{Generalized K\"ahler structures of even type on the standard Hopf surfaces}
We use the same notation as in Subsection \ref{Hopf odd}.
We define $E_\pm$ by 
$E_\pm:=r\frac{\pa}{\pa r}\pm\frac{dr}r$. Since $E_\pm\cdot E_\pm=\pm1$,   real $E_\pm$ are elements of the real Pin group. Since the real Pin group acts on the set of almost generalized K\"ahler structures, 
the pair $(E_\pm\cdot \phi, E_\pm\cdot \psi)$ gives an almost generalized K\"ahler structure.
We shall show that $(E_\pm\cdot \phi, E_\pm\cdot \psi)$ is a generalized K\"ahler structure.
 \bgn{lemma}\label{key lemma 2}
 $E_\pm\cdot (\J_\phi\frac{dr}r)+ (\J_\phi\frac{dr}r)\cdot E_\pm=0$
 \end{lemma}
 \bgn{proof}
 It suffices to show that
 $\lan E_\pm, \,\, \J_\phi\frac{dr}r\ran_{\tt}=0$.
 Since  one has $\lan\frac{dr}r, \J_\phi \frac{dr}r\ran_{\tt}=0$, 
it follows 
 $\lan E_\pm, \,\,\J_\phi\frac{dr}r\ran_{\tt}=\lan r\frac{\pa}{\pa r},\,\,\,\J_\phi\frac{dr}r\ran_{\tt}$. From Lemma \ref{key lemma 1}, one has
$$-\J_\psi\J_\phi\frac{dr}r=r\frac{\pa}{\pa r}.
$$
Thus one has
\bgn{align*}
\lan r\frac{\pa}{\pa r},\,\,\,\J_\phi\frac{dr}r\ran_{\tt}=&
\lan -\J_\psi\J_\phi\frac{dr}r, \,\,\,\J_\phi\frac{dr}r\ran_{\tt}\\
=&\lan -\J_\psi\frac{dr}r, \,\,\,\frac{dr}r\ran_{\tt}=0
\end{align*}
Hence one has $\lan E_\pm, \,\, \J_\phi\frac{dr}r\ran_{\tt}=0$.
\end{proof}
\bgn{proposition}
$(E_\pm\cdot \phi, E_\pm\cdot \psi)$ are $H$-twisted generalized K\"ahler structures.
\end{proposition}
\bgn{proof}
We need to calculate $d_H(E_\pm\cdot\frac\phi r)$. Since $d_H\circ E_\pm =L_{E_\pm}^H- E_\pm \circ d_H$ and 
$L^H_{E_\pm}(\frac\phi r)=L_{E_\pm}(\frac\phi r)+(i(r\frac{\pa}{\pa r})H)\w (\frac\phi r)=0
$, we have
\bgn{align*}
d_H(E_\pm\cdot\frac\phi r) =&-E_\pm d_H(\frac\phi r)=E_\pm\cdot\frac{dr}r\cdot(\frac\phi r)\\
=&-\sqrt{-1}E_\pm\cdot \J_\phi\frac{dr}r\cdot(\frac\phi r) 
\end{align*}
Applying Lemma \ref{key lemma 2}, we have
\bgn{align}
d_H(E_\pm\cdot\frac\phi r)=&\sqrt{-1} (\J_\phi\frac{dr}r)\cdot E_\pm\cdot(\frac\phi r) 
\end{align} 
Hence 
both $(E_\pm\cdot\phi)$ give integrable generalized complex structures.
It also follows that  $(E_\pm\cdot\psi)$ give integrable generalized complex structures.
Hence 
$(E_\pm\cdot \phi, E_\pm\cdot\psi)$
are generalized K\"ahler structures. 
\end{proof}
Next we shall calculate the scalar curvature of  the generalized K\"ahler structures 
$(\J_{E_\pm\cdot\phi}, \,\, \J_{E_\pm\cdot\psi})$ which are given by $(E_\pm\cdot\phi, E_\pm\cdot\psi)$.
We define a volume form $\vol_{E_\pm}$ by $\vol_{E_\pm}=\mp\vol.$ Then we have the normalization condition
$$
i^{-2}\lan E_\pm \til\phi, \,\, E_\pm \ol{\til\phi}\ran_s=i^{-2}\lan E_\pm \cdot\psi, \,\,E_\pm \cdot\ol{\til\psi}\ran_s=\vol_{E_\pm}
$$
We already have
 \bgn{align*}
d_H(E_\pm\cdot\frac\phi r)=&\sqrt{-1} (\J_\phi\frac{dr}r)\cdot E_\pm\cdot(\frac\phi r) \\
d_H(E_\pm\cdot\frac\psi r)=&\sqrt{-1} (\J_\psi\frac{dr}r)\cdot E_\pm\cdot(\frac\psi r) 
\end{align*} 
Thus we have 
 $\eta=\sqrt{-1} (\J_\phi\frac{dr}r),$ $\zeta=\sqrt{-1} (\J_\psi\frac{dr}r)$.
Applying Lemma \ref{key lemma 1} and  Lemma \ref{key lemma 2} again, we have  
 \bgn{align*}
 d_H (\eta\cdot E_\pm\cdot\frac\psi r )=&\sqrt{-1}d_H(\J_\phi\frac{dr}r)\cdot E_\pm\cdot \frac\psi r
 =-\sqrt{-1}d_H E_\pm\cdot(\J_\phi\frac{dr}r)\cdot\frac\psi r\\
 =&-d_H E_\pm\cdot(\J_\psi \J_\phi\frac{dr}r)\cdot\frac\psi r
 =d_H E_\pm\cdot(r\frac{\pa}{\pa r})\cdot\frac\psi r
 \end{align*}
Since 
 $L_{E_\pm}^H((r\frac{\pa}{\pa r})\cdot\frac\psi r)=L_{r\frac{\pa}{\pa r}}((r\frac{\pa}{\pa r})\cdot\frac\psi r)=
 ((r\frac{\pa}{\pa r})\cdot L_{r\frac{\pa}{\pa r}}\frac\psi r)=0$, one has
 \bgn{align*}
 d_H (\eta\cdot E_\pm\cdot\frac\psi r )=&-E_\pm\cdot d_H (r\frac{\pa}{\pa r})\cdot\frac\psi r=
 E_\pm\cdot(r\frac{\pa}{\pa r})\cdot d_H \frac\psi r\\
 =&-E_\pm\cdot(r\frac{\pa}{\pa r})\cdot \frac{dr}r\frac{\psi}r
\end{align*}
Taking the complex conjugate $\ol\psi$, we have 
\bgn{align*}
i^{-2}\lan E_\pm\cdot\frac\psi r, \,\, d_H(\eta\cdot E_\pm\cdot\frac{\ol\psi}r)\ran_s=&
i^{-2}\lan E_\pm\cdot\frac\psi r, \,\,\,E_\pm\cdot(r\frac{\pa}{\pa r})\cdot \frac{dr}r\frac{\ol\psi} r\ran_s\\
=&\mp i^{-2}\lan \frac\psi r, \,\,\,(r\frac{\pa}{\pa r})\cdot \frac{dr}r\frac{\ol\psi}r\ran_s\\
=&\mp\frac12 i^{-2}\lan\frac\psi r, \,\,\,\frac{\ol\psi}r\ran_s
\end{align*}
Hence we obtain 
\bgn{align*}
\Re \,\,i^{-n}\(\lan E_\pm\cdot\frac\psi r, \,\, d_H(\eta\cdot E_\pm\cdot\frac{\ol\psi}r)\ran_s\)=&
\mp\frac12 i^{-2}\lan\frac\psi r, \,\,\,\frac{\ol\psi}r\ran_s
\end{align*}
We also have 
\bgn{align*}
\Re\,\,i^{-n}\(\lan E_\pm\cdot\frac\phi r, \,\, d_H(\zeta\cdot E_\pm\cdot\frac{\ol\phi}r)\ran_s\)=&
\mp\frac12 i^{-2}\lan\frac\phi r, \,\,\,\frac{\ol\phi}r\ran_s 
\end{align*}
Hence 
we have 
\bgn{align}\label{scalr curv of Hopf even}
\Re\,&\(
\,i^{-n}\lan E_\pm\cdot\frac{\til\psi} r, \,\, d_H(\eta\cdot E_\pm\cdot\frac{\ol{\til\psi}}r)\ran_s
+\,\,i^{-n}\lan E_\pm\cdot\frac{\til\phi} r, \,\, d_H(\zeta\cdot E_\pm\cdot\frac{\ol{\til\phi}}r)\ran_s\)\\
=&\vol_{E_\pm}\notag
\end{align}
\bgn{proposition}
The generalized scalar curvature of 
$(\J_{E_\pm\cdot\phi}, \,\, \J_{E_\pm\cdot\psi})$ is $1,$
where the volume form is given by $\vol_{E_\pm}.$
 \end{proposition}
 \bgn{proof}
 The result follows from (\ref{scalr curv of Hopf even}).
 \end{proof}

  \section{Review of results on Lie groups, due to Alekseev, Bursztyn and Meinrenken}\label{Review} 
 We give a review of results of the paper \cite{ABM} on pure spinors and Dirac structures of Lie groups, which we will need to obtain generalized K\"ahler structures on compact Lie groups. 
The results in \cite{ABM} are also effective to calculate the scalar curvature of 
the generalized K\"ahler structures on compact Lie groups. 
  \subsection{The isomorphism $\TT \cong G\times (\mathfrak g\oplus \ol{\mathfrak g})$}
Let $G$ be a Lie group, and let $\mathfrak g$ be its Lie algebra. We denote by $\xi^R, \xi^L\in {\mathfrak X}(G)$ the left-, right-invariant vector field on $G$ which are equal to $\xi\in \mathfrak g\cong T_eG$ at the group unit. 
Let $\t^L, \t^R\in\Omega(G)\otimes\mathfrak g$ be the left-, right-Maurer-Cartan forms, i.e., 
$\iota(\xi^L)\t^L=\iota(\xi^R)\t^R=\xi.$
They are related by $\t^R_g=\Ad_g(\t^L_g),$ for all $g\in G.$
the corresponding infinitesimal action is given by 
the vector fields 
$$
\A_{ad}(\xi)=\xi^L-\xi^R
$$
Suppose that the Lie algebra $\mathfrak g$ of $G$ carries an {\it invariant inner product}.
By this we mean an $\Ad$-invariant, non-degenerate symmetric bilinear form $B$. Equivalently, $B$ defines a bi-invariant pseudo Riemannian metric on $G.$
Given $B$, we can define the Ad-invariant $3$-form $H\in \Ome^3(G),$
$$
H:=\frac1{12}B(\t^L, \,\,[\t^L, \t^L])
$$
Since $H$ is bi-invariant, it is closed. and it defines an $H$-twisted Courant bracket 
$[\,,\,]_H$ on $G.$
The conjugation action $\A_{ad}$ extends to an action of "the double group" $D=G\times G$ on $G$, by 
$$
\A: D\to \Diff(G),\quad \A(a, a')=l_a\circ r_{(a')^{-1}},
$$
where $l_a(g)=ag$ and $r_a(g)=ga.$
The corresponding infinitesimal action 
$$
\A: \mathfrak d \to \mathfrak X(G), \qquad \A(\xi, \xi')=\xi^L-(\xi')^R
$$
lifts to a map 
$$
s:\mathfrak d\to \Gam(\TT), \qquad 
s(\xi, \xi')=s^L(\xi)-s^R(\xi'),
$$
where 
$$
s^L(\xi )=\xi^L\oplus B(\t^L, \xi), \quad s^R(\xi')=-(\xi')^R\oplus B(\t^R, \xi)
$$
Let us equip $\mathfrak d$ with the bilinear form $B_{\mathfrak d}$ given by 
$+B$ on the first $\mathfrak g$-summand and $-B$ on the second summand. 
Thus $\mathfrak d=\mathfrak g\oplus\ol{\mathfrak g}$ is an example of a Lie algebra with invariant split bilinear form. 
\bgn{proposition}\text{\rm\cite{ABM}}\label{ABM1}
The map $s: \mathfrak d\to \Gam(TG\oplus T^*G)$ is $D$-equivariant, and satisfies 
\bgn{align*}
&(1)\quad  \lan s(\zeta_1), \, s(\zeta_2)\ran_{\tt}=B_{\mathfrak d}(\zeta_1, \zeta_2), \\
&(2)\quad [s(\zeta_1), \, s(\zeta_2)]_{H}=s([\zeta_1, \,\, \zeta_2]_\mathfrak d)
\end{align*}
for all $\zeta_1, \zeta_2\in\mathfrak d$
\end{proposition}

\bgn{proof} 
For completeness of our paper, we will give a proof. \\
(1) It suffices to show that 
$$ \lan s(\zeta), \, s(\zeta)\ran_{\tt}=B_{\mathfrak d}(\zeta, \zeta)$$ for all $\zeta=(\xi, \xi').$
Since $s(\zeta)$ is given by 
$$
s(\zeta)=\xi^L-\xi^R+B(\t^L, \xi)+B(\t^R, \xi')
$$
one has 
\bgn{align*}
\lan s(\zeta), \,\,s(\zeta)\ran_{\tt}=&B(\t^L, \xi)(\xi^L)+B(\t^R,\xi')(\xi^L)\\
-&B(\t^L, \xi)({\xi'}^R)-B({\t'}^R, \xi)({\xi'}^R)
\end{align*}
From the definition of the Maurer-Cartan form, one has $B(\t^L, \xi)(\xi^L)=B(\xi, \xi)$ and 
$B(\t^R, \xi')({\xi'}^R)=B(\xi', \xi').$
Since $\xi^L=(\Ad_g\xi)^R$ and $\xi^R=(\Ad_{g^{-1}}\xi)^L$ at $g\in G,$
we  have 
\bgn{align*}B(\t^R, \xi)(\xi^L)=&B(\t^R, \xi')(\Ad_g\xi)^R=B(\Ad_g\xi, \, \xi') \\
B(\t^L, \xi)({\xi'}^R)=&B(\t^L, \xi)(\Ad_{g^{-1}}(\xi'))^L=B(\Ad_{g^{-1}}(\xi'), \xi)
\end{align*}
Since $B$ is Ad-invariant, it follows $B(\Ad_g\xi, \, \xi')=B(\Ad_{g^{-1}}(\xi'), \xi).$
Thus we have 
$$
\lan s(\zeta), \,\,s(\zeta)\ran_{\tt}=B(\xi, \xi)-B(\xi', \xi')
$$
(2) We shall show that $[s(\zeta_1), \, s(\zeta_2)]_{H}=s([\zeta_1, \,\, \zeta_2])$ for all 
$\zeta_1=(\xi_1,\xi'_1), \zeta_2=(\xi_2,\xi'_2).$
Let $s^L$ be the following 
\bgn{align*}
s^L(\xi_1)=&\xi_1^L+\t_1^L,\qquad
s^L(\xi_2)=\xi_2^L+\t_2^L,
\end{align*}
where $\t_1^L=B(\t^L, \xi_1)$ and $\t_2^L=B(\t^L, \xi_2).$
We also denote by $s^R$ the following 
\bgn{align*}
s^R({\xi'}_1)=&-{\xi'}_1^R+\t_1^R,\qquad
s^R({\xi'}_2)=-{\xi'}_2^R+\t^R_2,
\end{align*}
where $\t_1^R=B(\t^R, {\xi'}_1)$ and $\t_2^R=B(\t^R, {\xi'}_2)$.
Then one has
\bgn{align*}
s(\zeta_1)=&s^L(\xi_1)+s^R({\xi'}_1),\quad
s(\zeta_2)=s^L(\xi_2)+s^R({\xi'}_2)
\end{align*}
From the definition of $H$-twisted bracket, one has
\bgn{align*}
[s^L(\xi_1), \, s^L(\xi_2)]_H=&[\xi_1^L,\,\xi_2^L]+\frac12\(L_{\xi_1^L}\t_1^L-L_{\xi_2^L}\t_1^L\)\\
-&\frac12\(i_{\xi_2^L}d\t^L_1-i_{\xi_1^L}d\t^L_2\)-i_{\xi_2^L}i_{\xi^L_1}H
\end{align*}
Since $i_{\xi^L_1}\t_2^L=B(\xi_1, \xi_2)$ and $di_{\xi^L_1}\t_2^L=0,$
then we have 
\bgn{align*}
[s^L(\xi_1), \, s^L(\xi_2)]_H=&[\xi_1^L,\,\xi_2^L]-\(i_{\xi_2^L}d\t^L_1-i_{\xi_1^L}\t^L_2\)-i_{\xi_2^L}i_{\xi^L_1}H
\end{align*}
Since $i_{\xi_2}d\t_1^L=-B(\t^L \,,[\xi_1, \xi_2]_{\mathfrak g})$, we have 
\bgn{align*}
[s^L(\xi_1), \, s^L(\xi_2)]_H=&[\xi_1^L,\,\xi_2^L]+2B(\t^L, [\xi_1, \xi_2]_{\mathfrak g})-i_{\xi_2^L}i_{\xi^L_1}H
\end{align*}
Since $i_{\xi_3^L}i_{\xi_2^L}i_{\xi^L_1}H=- B([\xi_1, \xi_2], \,\,\xi_3),$
we have 
\bgn{align*}
[s^L(\xi_1), \, s^L(\xi_2)]_H=&[\xi_1^L,\,\xi_2^L]+B(\t^L, [\xi_1, \xi_2]_{\mathfrak g})
\end{align*}
$[s^R({\xi'}_1), \, s^R({\xi'}_2)]_H$ is calculated by using $[{\xi'}^R_1, {\xi'}^R_2]=-[{\xi'}_1, {\xi'}_2]^R$
and 
$i_{\xi_2^R}d\t_1^R=B(\t^R, [\xi_1, \xi_2]_{\mathfrak g}),$
\bgn{align*}
[s^R(\xi_1), \, s^R(\xi_2)]_H=&[{\xi'}_1^R,\,{\xi'}_2^R]+2B(\t^R, [{\xi'}_1, {\xi'}_2]_{\mathfrak g})+
i_{{\xi'}_2^R}i_{{\xi'}^R_1}H\\
=&-[{\xi'}_1, {\xi'_2}]_{\mathfrak g}^R+B(\t^R, [{\xi'}_1, {\xi'}_2]_{\mathfrak g})
\end{align*}
Then we have 
\bgn{align}
[s^L(\xi_1), \, s^L(\xi_2)]_H=&s^L([\xi_1, \xi_2]_{\mathfrak g}),\quad
[s^R({\xi'}_1), \, s^R({\xi'}_2)]_H=s^R([{\xi'}_1, {\xi'}_2]_{\mathfrak g})
\end{align}
We also have 
\bgn{align*}
[s^L(\xi_1), \, s^R({\xi'}_2)]_H=&0,\qquad
[s^R({\xi'}_1), \, s^R({\xi}_2)]_H=0
\end{align*}
Thus we obtain the result.
\end{proof}
\subsection{ The isomorphism $\w T^*G\cong G\times Cl(\mathfrak g)$}
Let us now assume that the adjoint action $\Ad: G\to O(\mathfrak g)$ lifts to a group homomorphism 
\bgn{equation}\label{75}
\tau: G\to \text{Pin}(\mathfrak g)\subset Cl(\mathfrak g)
\end{equation}
If $G$ is connected and $\pi_1(G)$ is torsion free, then this is automatic. 
Note that (\ref{75}) is consistent with our previous notation 
$\tau(\xi)=q(\lam(\xi)), $ since 
$$
\tau(\xi)=\frac{d}{dt}\Big|_{t=0}\tau(\exp(t\xi))
$$
We will write $N(g)=N(\tau(g))=\pm1$ for the image under the norm homomorphism, and 
$|g|=|\tau(g)|$ for the parity of $\tau (g).$
Since $\tau(g)$ lifts $\Ad_g$, one has $(-1)^{|g|}=\det\Ad_g.$
The definition of the Pin group implies that conjugation by $\tau(g)$ is the twisted adjoint action, 
$$
\tau(g)x\tau(g^{-1})=\wtil\Ad_g(x):=(-1)^{|g|\,|x|}\Ad_g(x)
$$
This twisted adjoint action extends to an action of the group $D$ on $Cl(\mathfrak g)$, 
$$
\A^{\cl}(a, a')(x)=\tau(a)x\tau((a')^{-1})
$$
Let us now fix a generator $\mu\in \det(\mathfrak g)$, and consider 
the corresponding star operator $\star: \w\mathfrak g^*\to \w\mathfrak g.$
The star operator satisfies 
$$
\Ad_g\circ\star =(-1)^{|g|}\star\circ\Ad_{g^{-1}}^*
$$

We use the left-invariant forms to identify $\w T^*G\cong G\times \w\mathfrak g^*$. 
Applying $\star$ point-wise, we obtain an isomorphism $q\circ\star:\w T^*_gG  \to Cl(\mathfrak g)$ for each $g\in G.$
Let us define the linear map 
\bgn{equation}\label{mathcal R: Cl frak g to Ome G }
\mathcal R: Cl(\mathfrak g)\to \Ome(G), \quad \mathcal R(x)\Big|_{g}=(q\circ\star)^{-1}(x\tau(g))
\end{equation}
We denote by $\mu^*$ the dual generator defined by $\iota((\mu^*)^T)\mu=1,$ and 
let $\mu_G$ be the left-invariant volume form on $G$ defined by $\mu^*.$
Let $q: \w\mathfrak g\to Cl(\mathfrak g)$ be the natural inclusion from the skewsymmetric algebra to the Clifford algebra. 
We define $\Xi\in \w^3\mathfrak g$ by 
$$
i_{\xi_3}i_{\xi_2}i_{\xi_1}\Xi =B(\xi_1, \,\,[\xi_2, \xi_3])
$$
We denote by $q(\Xi)\in Cl(\mathfrak g)$ the image of the $3$-form $\Xi$ by the map $q.$ 
We define the Clifford differential $d^{cl}: Cl(\mathfrak g)\to Cl(\mathfrak g)$ by the action of $\Xi$
$$
d^{cl}(x)=-[q(\Xi), \,\,]_{Cl}.
$$
\bgn{proposition}\text{\rm\cite{ABM}}
The map (\ref{mathcal R: Cl frak g to Ome G }) has the following properties: 
\bgn{enumerate}
\item[\rm (a)] $\mathcal R$ intertwines the Clifford actions, in the sense that 
$$
\mathcal R(\rho^{cl}(\zeta)x)=\rho(s(\zeta))\mathcal R(x), \qquad \forall x\in Cl(\mathfrak g), \, \zeta
\in \mathfrak d
$$
Up to a scalar function, $\mathcal R$ is uniquely characterized by this property.
\item[\rm (b)]
$\mathcal R$ intertwines differentials: 
\bgn{equation}\label{mathcal Rdcl(x)=d+H}
\mathcal R(d^{cl}(x))=(d+H)\mathcal R(x) 
\end{equation}
\item[\rm (c)]$\mathcal R$ satisfies the following $D$-equivariant  property: For any 
$h=(a, a')\in D$, and 
at any given point $g\in G$, and for all $x, x'\in Cl(\mathfrak g)$, 
$$
\A(h^{-1})^*\mathcal R(x)=(-1)^{|a|(|g|+|x|)}\mathcal R(\A^{cl}(h)x)
$$
\item[\rm (d)]
$\mathcal R$ relates the bilinear pairings on the Clifford modules $Cl(\mathfrak g)$ and $\Ome(G)$ as follows:
At any given point $g\in G,$ and for all $x, x'\in Cl(\mathfrak g)$, 
\bgn{equation}\label{R(x) R(x')Clfrak g}
\lan\mathcal R(x), \,\, \mathcal R(x')\ran_s=(-1)^{|g|(\dim G+1)}N(g)(x, x')_{Cl(\mathfrak g)}\mu_G
\end{equation}
Here the pairing $(.\,.)_{Cl(\mathfrak g)}$ is viewed as scalar-valued, using the trivialization of $\det(\mathfrak g)$ defined by $\mu_G$.
\end{enumerate}
Notice that the signs in part (c), (d) disappear if $G$ is connected.
\end{proposition}

 \section{Scalar curvature of generalized K\"ahler structures on compact Lie groups}\label{GK on compact Lie group}
Typical generalized K\"ahler structures on compact Lie groups are known \cite{Gua_2003}.
We shall give a different construction of generalized complex structures on 
compact Lie groups by using Proposition \ref{ABM1}.
 \subsection{Construction of generalized K\"ahler structures on compact Lie groups}
 Let $G$ be a connected compact Lie group of even dimension, and let $\mathfrak g$ be its Lie algebra.
 Suppose that $\mathfrak g$  carries a bi-invariant metric $B$ which is positive definite in this section.
 Let $\mathfrak g^\C$ be the complexification $\mathfrak g\otimes \C$ of $\mathfrak g.$
 The metric $B$ is linearly extended to the bilinear form on $\mathfrak g^\C$ over $\C.$
 We denote by $\mathfrak h$ a Cartan subalgebra of $\mathfrak g^\C$ with 
 $\dim_{\Bbb C}\mathfrak h=2l$. 
 Then we have the root decomposition: 
 $$
 \mathfrak g^\C =\mathfrak h\oplus \bigoplus_\a \mathfrak g_\a
 $$
 Choosing a system of positive roots, and for each positive root $\a>0$
 we have $\{e_\a\in \mathfrak g_\a\, |\,  \a>0\}$. We have the triangular decomposition
 $$
 \mathfrak g^\C =\mathfrak h \oplus \bigoplus_{\a>0}\mathfrak g_\a\oplus\bigoplus_{\a<0}\mathfrak g_\a
 $$ 
 Let $I_{\mathfrak h}$ be an almost  complex structure  on the real part $\mathfrak h^{\Bbb R}=\mathfrak g\cap \mathfrak h$ of $\mathfrak h$ which is compatible with the bi-invariant metric $B.$
Then $I_{\mathfrak h}$ gives the decomposition of $\mathfrak h$ into the following $(1,0)$ and $(0,1)$-components:
 $$ \mathfrak h =\frak h^{1,0}\oplus \frak h^{0,1}.$$
We define an almost complex structure on $I_{\mathfrak g}$ on $\mathfrak g$ by taking the following $(1,0)$ and 
  $(0,1)$-components
 $$
 \mathfrak g^\C =\mathfrak g^{1,0}\oplus\mathfrak g^{0,1},
 $$ 
 where $\mathfrak g^{1,0}=\mathfrak h^{1,0}\oplus \bigoplus_{\a>0}\mathfrak g_\a$ and $\mathfrak g^{0,1}=\mathfrak h^{0,1}\oplus\bigoplus_{\a<0}\mathfrak g_\a.$
 The complex conjugation is given by $\ol e_\a=e_{-\a}.$
 \bgn{proposition}
 $\mathfrak g^{1,0}$ is closed under the bracket of the Lie algebra $\mathfrak g.$
 The metric $B$ is a symmetric form of type $(1,1).$
 \end{proposition}
 \bgn{proof}  
Since $[e_\a, e_\b]\subset\frak g_{\a+\b},$ for $e_\a\in\frak g_\a,$ and $ e_\b\in\frak g_\b$,
the space of positive roots $\oplus_{\a>0}\frak g_\a$ is closed under the bracket of the Lie algebra $\frak g.$
Since $[h, e_\a]=\a(h)e_\a$ and $\frak h$ is commutative, 
$\frak g^{1,0}$ is closed under the bracket of the Lie algebra $\frak g.$
The metric $B$ is ad-invariant, it follows that 
$B([h, e_\a], e_\b)+B(e_\a, [h, e_\b])=0$ for $h\in\frak h$ and 
$e_\a\in\frak g_\a, e_\b\in\frak g_\b.$
Then one has $\a(h)B(e_\a, e_\b)+\b(h)B(e_\a, e_\b)=0.$
If $\a, \b$ are positive, $\a+\b\neq 0.$
Thus it follows $B(e_\a, e_\b)=0$ for all positive roots $\a, \b.$ 
Since $B$ is real, it follows $B(e_{-\a}, e_{-\b})=0.$
Thus $B$ is of type $(1,1).$
 \end{proof}
We choose a unitary basis $\{ u_i\}$ of $\frak g^{1,0}$ with respect to $B$ and denote by $\{\ol u_i\}$ the conjugate basis of $\frak g^{0,1}$.
We recall the Lie algebra $\frak d=\frak g\oplus \ol{\frak g}$ which is the direct sum of $\frak g$ and $\frak g.$ 
Let $\frak d$ equipped with the bilinear form
$B_{\frak d}$ given by $+B$ on the first summand and $-B$ on the second summand. 

\bgn{definition}
We define $\L_\phi$ by the direct sum $\frak g^{1,0}\oplus\frak g^{1,0}$ which is a subspace of $\frak d.$
We also define $\L_\psi$ by the direct sum $\frak g^{1,0}\oplus \frak g^{0,1}.$
\end{definition}
Then one has
\bgn{proposition}
The subspace $\L_\phi$ is maximal isotropic, which gives a decomposition $\mathfrak d\otimes\C=\L_\phi\oplus\ol\L_\phi.$
The subspace $\L_\psi$ is also maximal isotropic , which gives a decomposition 
$\mathfrak d\otimes\C=\L_\psi\oplus \ol\L_\psi.$
\end{proposition} 
\bgn{proof}
Since $B$ is of type $(1,1),$
we see that 
$$B_{\frak d}((u_i, u_j),\,\, (u_i, u_j))=B(u_i, u_i)-B(u_j, u_j)=0,$$
for all $(u_i, u_j)\in \L_\phi.$
Thus $\L_\phi$ is maximal, isotropic. 
One also has
$$B_{\frak d}((u_i, \ol u_j),\,\, (u_i, \ol u_j))=B(u_i, u_i)-B(\ol u_j, \ol u_j)=0,$$
Thus $\L_\psi$ is also maximal isotropic.
\end{proof}
\bgn{proposition}\label{Both Lphi and Lpsi are closed}
Both $\L_\phi$ and $\L_\psi$ are closed under the bracket of the Lie algebra $\frak d.$
\end{proposition}
\bgn{proof}
Since $\frak g^{1,0}$ is closed under the bracket, 
thus $\L_\phi=\frak g^{1,0}\oplus \frak g^{1,0}$ is closed.
Since $\frak g^{0,1}$ is also closed under the bracket,
thus $\L_\psi=\frak g^{1,0}\oplus\frak g^{0,1}$ is closed also.
\end{proof}
We define $\L_\phi^+$ and $\L_\phi^-$ by  
$$
\L_\phi^+ =\L_\phi\cap \L_\psi, \qquad \L_\phi^-=\L_\phi\cap \ol \L_\psi
$$
We also denote by $\ol\L_\phi^\pm$ the conjugate of $\L_\phi^\pm$, respectively.
Then we have that 
\bgn{align}
\L_\phi^+=&\frak g^{1,0}\oplus 0, \qquad \L_\phi^-=0\oplus \frak g^{1,0}\\
\ol\L_\phi^+=&\frak g^{0,1}\oplus 0, \qquad \ol\L_\phi^-=0\oplus \frak g^{0,1}
\end{align}
We define real subspaces $C_\pm $ by 
$$
C_\pm = (\L_\phi^\pm\oplus\ol\L_\phi^\pm)\cap \frak g 
$$
Then we see that 
\bgn{lemma}\label{C+=frak goplus 0}
$C_+=\frak g\oplus 0$ and $C_-=0\oplus\frak g$.
\end{lemma}
\bgn{proof}
The result follows from our definition of $C_\pm.$
\end{proof}
Thus we have 
\bgn{proposition}\label{the pair of isotropic subspaces}
The pair of isotropic subspaces $\L_\phi$ and $\L_\psi$ gives a generalized K\"ahler structure on the vector space
$\frak d=\frak g
\oplus \ol{\frak g}.$
\end{proposition}
\bgn{proof}
The result follows from Lemma \ref{C+=frak goplus 0}.
\end{proof}
 \bgn{proposition}\label{GK on G}
The pair $(\L_\phi, \L_\psi)$ gives rise to a generalized K\"ahler structure 
$(\J_\phi, \J_\psi)$ on 
a Lie group $G.$
\end{proposition} 
\bgn{proof}
Applying Proposition \ref{ABM1} and Proposition  \ref{Both Lphi and Lpsi are closed}, 
we have $\J_\phi$ and $\J_\psi$ give integrable generalized complex structures on $G$.
Then the  result follows from  Proposition \ref{the pair of isotropic subspaces}.
\end{proof}

\subsection{Calculation of the generalized scalar curvature of generalized K\"ahler structures on compact Lie groups}
Our observation of calculations of the generalized scalar curvature of generalized K\"ahler structures on the Hopf surfaces can be extended to the cases of compact Lie groups.
We shall show that the generalized scalar curvature of $(G,\J_\phi, \J_\psi)$ is a constant.
Let $(\J_\phi, \J_\psi)$ be a generalized K\"ahler structure on $G$ as in Proposition \ref{GK on G}.
Then we shall show the following:
\bgn{proposition}\label{S=constant}
The generalized scalar curvature $S(\J_\phi, \J_\psi)$ is a constant.
\end{proposition}
We need several lemmas to show Proposition \ref{S=constant}. 
The proof will be given at the end of this section.
In order to show that $S(\J_\phi, \J_\psi)$ is a constant,
we determine  nondegenerate, pure spinors corresponding to $\J_\phi, \J_\psi$, which
are given by nondegenerate, pure spinors of $Cl(\frak g)$.
As in (\ref{mathcal R: Cl frak g to Ome G }), the map
$$
\mathcal R :  Cl(\frak g)\to \Ome(G)
$$
gives a map from nondegenerate, pure spinors of $Cl(\frak g)$ to 
nondegenerate, pure spinors of $\Omega(G).$
{Note that  there is no loss of generality, we can assume nondegenerate, pure spinors $\mathcal R(\phi), \mathcal R(\psi)$ are globally defined on G.}
We shall give an explicit description of $(\phi, \psi)$ which gives $(\J_\phi, \J_\psi)$ on $G.$
Let $\{t_i\}$ be a unitary basis of $\frak h^{1,0}$ and $\{\ol t_i\}$ the conjugate basis of $\frak h^{0,1}.$
We denote by $\{\t_\a\}_{\a>0}$ a unitary basis of $\oplus_{\a>0}\frak g_\a$ and 
by $\{\ol\t_\a\}_{\a>0}$ the conjugate basis of $\oplus_{\a<0}\frak g_\a$, where note that
$\ol\t_\a=\t_{-\a}.$
Then we have
\bgn{lemma}\label{phiand psi are given by}
$\phi$ and $\psi$ are given by 
\bgn{align*}
\phi=&\(\prod_{i=1}^l t_i\)\cdot\(\prod_{\a>0}\t_\a\)\\
\psi=&\(\prod_{i=1}^lt_i\cdot \ol t_i\)\cdot \(\prod_{\a>0}\t_\a\cdot\ol\t_\a\)
\end{align*}
\end{lemma}
\bgn{proof}
Since $t_i\cdot t_j=0$ for all $i,j=1,\cdots, l$ and $\t_\a\cdot\t_\b=0$ for all positive roots $\a, \b$,
it follows $t_i\cdot\phi=\phi\cdot t_i=0$ and $\t_\a\cdot\phi=\phi\cdot \t_\a=0.$
Thus $\ker\phi =\frak g^{1,0}\oplus \frak g^{1,0}.$
Then it also follows $\ker\psi=\frak g^{1,0}\oplus\frak g^{0,1}.$
\end{proof}
We define $\vol_{\frak g}$ by 
$$
\vol_{\frak g}=i^{-n}\(\prod_{i=1}^lt_i\cdot \ol t_i\)\cdot \(\prod_{\a>0}\t_\a\cdot\ol\t_\a\).
$$
Then we have  
$$
\vol_{\frak g}=i^{-n}\lan \phi,\ol\phi\ran_s=i^{-n}\lan \psi, \ol\psi\ran_s.
$$
The volume form on $G$ is given by 
$$
\vol_G:=i^{-n}\lan \mathcal R(\phi),\,\,\mathcal R(\ol\phi)\ran_s=i^{-n}\lan \mathcal R(\psi),\,\, \mathcal R(\ol\psi)\ran_s
$$
As in (\ref{mathcal Rdcl(x)=d+H}), $\mathcal R$ intertwines differentials: 
$\mathcal R(d^{cl}(x))=(d+H)\mathcal R(x)$. 
By using this, we shall calculate $d^{cl}\phi$ and $d^{cl}\psi.$
The Clifford differential $d^{cl}$ is given by $d^{cl}:=-[q(\Xi), \,]_{Cl}$ as before,
where $\Xi$ is the Cartan $3$-form of $\frak g$ and 
the map $q: \w^\bullet \frak g^\C\to Cl(\frak g^\C)$ gives 
$q(\Xi)\in Cl(\frak g^\C).$
\bgn{lemma}\label{the 3-form Xi is given by}
The $3$-form $\Xi$ is given by 
\bgn{align*}
\Xi=&\sum_{ \a>0}  2P_\a\w \t_\a\w\ol\t_{\ol\a}
+\sum_{\a+\b=  \gam}C_{\a,\b, \ol\gam}\t_\a\w\t_\b\w \ol\t_{\ol\gam}
+\sum_{\a+\b=\gam}C_{\ol\a, \ol\b, \gam}\ol\t_{\ol\a}\w\ol\t_{\ol\b}\w\t_\gam, 
\end{align*}
where 
$ P_\a=C_{i,\a,\ol\a}t_i+C_{\ol i, \a, \ol\a}\ol t_i\in \frak h$ is pure imaginary and 
$C_{i, \a, \ol\a}$, $C_{\a, \b,\ol\gam}$ and $C_{\ol\a, \ol\b, \gam}$ are structure constants of $\frak g^\C$
and the summation runs over all triples of positive roots $(\a, \b,\gam)$ satisfying $\a+\b=\gam.$ 
\end{lemma}
\bgn{proof}Recall the real $3$-form $\Xi$ is given by 
$$
\Xi(x,y,z)=B([x,y], z), \quad x,y,z\in \frak g.
$$
Since $[e_\a, e_\b]\in \frak g_{\a+\b}$ and $B$ is of type $(1,1)$, 
it follows $\Xi(\t_\a, \t_\b, \t_{\gam})=0$ for all positive roots $\a, \b, \gam.$
Since $\Xi$ is real, one has $\Xi(\ol\t_\a, \ol\t_\b, \ol\t_{\gam})=0$ for all positive roots $\a, \b, \gam.$ 
If  $[e_\a, e_\b]\in \frak h$, then $\a+\b=0.$ 
 Thus $\Xi$ is of type $(1,2)$ and of type $(2,1).$
 Both  $\Xi(t_i, e_\a, e_\b)$ and $\Xi(\ol t_i, e_\a, e_\b)$ give the first term and 
 $\Xi(e_\a, e_\b, \ol e_\gam)$ and $\Xi(\ol e_\a, \ol e_\b, e_\gam)$ give 
 the second and third term. 
 Since $\frak h$ is commutative, $\Xi(t_i, t_j, e_\a)=0$. Then we have the result.
\end{proof}
Let $P_\a=C_{i,\a,\ol\a}t_i+C_{\ol i, \a, \ol\a}\ol t_i\in \frak h$ as before.
We denote by  $ P$ the sum $\sum_{a>0} P_\a$ and 
by $e(P)$ the diagonal element  $(P,  P)$ of $\frak g\oplus\ol{\frak g}$ which acts on Cl$(\frak g^\C)$ by the Clifford multiplication
 $\rho^{cl},$ that is, 
$\rho^{cl}(e(P))x=P\cdot x-(-1)^{|x|}x\cdot P$.
We also denote by $a(P)$ the anti-diagonal element $(P, -P)$ which acts on $x\in Cl(\frak g^\C)  $ by 
$\rho^{cl}(a(P))x=P\cdot x+(-1)^{|x|}x\cdot P$
\bgn{lemma}\label{dclphi=rhoe(P)phi}
Then the Clifford differential of $\phi$ and $\psi$ are given by 
$$
d^{\cl}\phi=-\rho^{cl}(a(P))\phi, \qquad 
d^{\cl}\psi =-\rho^{cl}(e(P))\psi,
$$
where note that $e(P), a(P)$ are  pure imaginary elements of $\frak g^\C\oplus\ol{\frak g^\C}.$
\end{lemma}
\bgn{proof}
First we consider the Clifford differential $d^{cl}\phi$. 
Since $\{\t_\a\}_{\a>0}$ is a unitary basis of the space of positive roots and $\ol e_\a=e_{-\a},$ 
one has 
$$
\t_\a\cdot\t_\b+\t_\b\cdot\t_\a =
\bgn{cases}
&1 \qquad \a=-\b, \\
&0 \qquad \a\neq -\b
\end{cases}
$$
In particular, $\t_\a\cdot\t_\a=0$ for every positive root $\a.$
Recall that $q(\t_1\w\t_2\w\t_3)$ is given by 
$$
\t_1\w\t_2 \w\t_{ 3}=\frac16\sum_{\sig\in \frak S_3}(-1)^{|\sig|}\(\t_{\sig(1)}\cdot\t_{\sig(2)}\cdot \t_{\sig(3)}
\),
$$
where $\sig$ runs over all permutation of $1,2,3$ and $|\sig|$ denotes the signature of $\sig.$
Thus both right and left multiplications of $q(\t_\a\w\t_\b\w \ol\t_{\ol\gam})$ on
$\(\prod_{\a>0}\t_\a\)$ vanish since $\a, \b\neq \gam.$
It also follows 
that both right and left multiplications of $q(\ol\t_{\ol\a}\w\ol\t_{\ol\b}\w\t_{\gam})$ on
$\(\prod_{\a>0}\t_\a\)$ vanish since $\a, \b\neq \gam.$
We also have $q(\t_\a\w\ol\t_\a)=\frac12(\t_\a\cdot\ol\t_\a-\ol\t_\a\cdot\t_\a).$
Then it follows that 
\bgn{align}
q(\t_\a\w\ol\t_\a)\cdot\t_\a=&\frac12(\t_\a\cdot\ol\t_\a-\ol\t_\a\cdot\t_\a)\cdot\t_\a=\frac12\t_\a\\
\t_\a\cdot q(\t_\a\w\ol\t_\a)=&\frac12\t_\a\cdot(\t_\a\cdot\ol\t_\a-\ol\t_\a\cdot\t_\a)=-\frac12\t_\a \label{sign -}
\end{align}
Then it follows that the right multiplications of $2q(\t_\a\w\ol\t_\a)$ on $\(\prod_{\a>0}\t_\a\)$
is trivial and the left one is $(-1)$-times.
Thus one has 
 $$[q(2P_\a\w\t_\a\w\ol\t_\a), \phi]_{Cl}=P_\a\cdot\phi+(-1)^{|\phi|}\phi\cdot P_\a=
 \rho^{cl}(a(P_\a))\phi
 $$
As in Lemma \ref{the 3-form Xi is given by},
$\Xi$ is given by 
\bgn{align*}
\Xi=&\sum_{ \a>0}  2P_\a\w \t_\a\w\ol\t_{\ol\a}
+\sum_{\a, \,\b \neq  \gam>0}C_{\a,\b, \ol\gam}\t_\a\w\t_\b\w \ol\t_{\ol\gam}
+\sum_{\a,\,\b \neq\gam>0}C_{\ol\a, \ol\b, \gam}\ol\t_{\ol\a}\w\ol\t_{\ol\b}\w\t_\gam, 
\end{align*}
Thus 
one has 
$-d^{cl}\phi =\rho^{cl}(e(\Xi))\phi=\sum_\a \rho^{cl}(a(P_\a))\phi=\rho^{cl}(a(P))\phi$.
We also have

\bgn{align}
q(\t_\a\w\ol\t_\a)\cdot(\t_\a\cdot\ol\t_\a)=&\frac12(\t_\a\cdot\ol\t_\a-\ol\t_\a\cdot\t_\a)\cdot (\t_\a\cdot\ol\t_\a)=\frac12 (\t_\a\cdot\ol\t_\a)\\
(\t_\a\cdot\ol\t_\a)\cdot q(\t_\a\w\ol\t_\a)=&\frac12(\t_\a\cdot\ol\t_\a)\cdot(\t_\a\cdot\ol\t_\a-\ol\t_\a\cdot\t_\a)=\frac12(\t_\a\cdot\ol\t_\a)\label{sign +}
\end{align}
(Note that the sign of (\ref{sign +}) is different from the case of (\ref{sign -}).)
Thus one has 
$-d^{cl}\psi =\rho^{cl}(e(\Xi))\psi=\sum_\a [P_\a, \psi]_{Cl}=\rho^{cl}(e(P))\psi$
\end{proof}

Next we shall calculate 
$d^{cl}\rho^{cl}(a(P))\psi$ and $d^{cl}\rho^{cl}(e(P))\phi.$
First we shall show the following,
\bgn{lemma}\label{rhocle(eXi)rhocla(P)}
\bgn{align}
&\rho^{cl}(e(\Xi))\rho^{cl}(a(P))+\rho^{cl}(a(P))\rho^{cl}(e(\Xi))=\rho^{cl}(a([\Xi, \, P]_{Cl})\\
&\rho^{cl}(e(\Xi))\rho^{cl}(e(P))+\rho^{cl}(e(P))\rho^{cl}(e(\Xi))=\rho^{cl}(e([\Xi, \, P]_{Cl})
\end{align}
\end{lemma}
\bgn{proof}
We shall show the Lemma in a general form.
Let $E=(u_1, u_2)$ be a real element of $\frak h\oplus \ol{\frak h}$ which satisfies 
$B(E, E)_{\frak d}=B(u_1, u_1)-B(u_2, u_2)\neq 0$. Then $E$ is an element of Pin group Pin$(\frak d)$. 
We denote by $e(\Xi)$ the diagonal $(\Xi, \Xi)\in Cl(\frak g\oplus\ol{\frak g})$. 
$Cl(\frak g)$ is a $Cl(\frak g\oplus\ol{\frak g})$-module and we denote by $\rho^{cl}$ the action of 
$Cl(\frak g\oplus\ol{\frak g}) $ on $Cl(\frak g)$.
Then for $x\in \frak g$ 
\bgn{align*}
\rho^{cl}(e(\Xi))\rho^{\cl}(E)x=&\rho^{cl}(e(\Xi))\(u_1\cdot x -(-1)^{|x|}x\cdot u_2\)\\
=&\Xi\cdot u_1\cdot x-(-1)^{|x|+1}u_1\cdot x\cdot \Xi\\
-&(-1)^{|x|}\Xi\cdot x\cdot u_2+(-1)^{2|x|+1}x\cdot u_2\cdot\Xi\\
\\
\rho^{cl}(E)\rho^{\cl}(e(\Xi))x=&\rho^{cl}(E)\(\Xi\cdot x -(-1)^{|x|}x\cdot\Xi\)\\
=&u_1\cdot \Xi\cdot x -(-1)^{|x|}u_1\cdot x\cdot\Xi\\
-&(-1)^{|x|+1}\Xi\cdot x\cdot u_2+(-1)^{2|x|+1}x\cdot\Xi\cdot u_2
\end{align*}
Since $[e(\Xi), E]_{Cl}=\([\Xi, u_1]_{Cl}, \,\, [\Xi, u_2]_{Cl}\)=
\(\Xi\cdot u_1+u_1\cdot\Xi, \,\,\Xi\cdot u_2+u_2\cdot\Xi\)$ is even,  one has
\bgn{align}
\rho^{cl}([e(\Xi), E]_{Cl})=\rho^{cl}(E)\rho^{cl}(e(\Xi))+\rho^{cl}(e(\Xi))\rho^{cl}(E)
\end{align}
Applying $E=e(P)$ or $a(P)$, we obtain the result.
\end{proof}
As in before, $d^{cl}$ is given by the action $-\rho^{cl}(e(\Xi))$.
From Lemma \ref{rhocle(eXi)rhocla(P)} and $-d^{cl}\psi =\rho^{cl}(e(\Xi))\psi=\rho^{cl}(e(P))\psi$, we have
\bgn{align}\label{dclrhoa(P)psirhocle(Xi)}
-d^{cl}\rho^{cl}(a(P))\psi =&\rho^{cl}(e(\Xi))\rho^{cl}(a(P))\psi\\
=&-\rho^{cl}(a(P))\rho^{cl}(e(\Xi))\psi+\rho^{cl}(a([\Xi, P]_{Cl}))\psi\\
=&-\rho^{cl}(a(P))\rho^{cl}(e(P))\psi+\rho^{cl}(a([\Xi, P]_{Cl}))\psi
\end{align}
We also have 
\bgn{align}\label{dclrhocle(P)phi=rhocle(Xi)}
-d^{cl}\rho^{cl}(e(P))\phi =&\rho^{cl}(e(\Xi))\rho^{cl}(e(P))\phi\\
=&-\rho^{cl}(e(P))\rho^{cl}(e(\Xi))\phi+\rho^{cl}(e([\Xi, P]_{Cl}))\phi\\
=&-\rho^{cl}(e(P))\rho^{cl}(a(P))\phi+\rho^{cl}(e([\Xi, P]_{Cl}))\phi
\end{align}
\bgn{lemma}\label{Rei-nlanpsirhocla(P)}
\bgn{align*}
\Re \,\,i^{-n}\lan \psi,\,\,\rho^{cl}(a(P))\rho^{cl}(e(P))\ol\psi\ran_s=&B(P,P)i^{-n}\lan \psi, \,\,\ol\psi\ran_s\\
\Re \,\,i^{-n}\lan\phi,\,\,\rho^{cl}(e(P))\rho^{cl}(a(P))\ol\phi\ran_s=&B(P,P)i^{-n}\lan \phi, \,\,\ol\phi\ran_s,
\end{align*}

\end{lemma}
\bgn{proof} Since $P$ is pure imaginary, the result follows.
\end{proof}

\bgn{lemma}\label{rhocleXi,Pclphi}
\bgn{align}
\rho^{cl}(e([\Xi, P]_{Cl})\phi=&2B(P,P)\phi\\
\rho^{cl}(a([\Xi, P]_{Cl})\psi=&2B(P,P)\psi
\end{align}
\end{lemma}
\bgn{proof}
We have 
\bgn{align*}
[\Xi, P]_{Cl}=&\sum_\a [P_\a\cdot(\t_\a\cdot\ol\t_\a-\ol\t_a\cdot\t_\a), \,\,P]_{Cl}\\
=&\sum_\a B(P_\a, P)\(\t_\a\cdot\ol\t_\a-\ol\t_a\cdot\t_\a\)
\end{align*}
Then we have
\bgn{align*}
\rho^{cl}(e([\Xi, P]_{Cl})\phi=&\sum_\a B(P_\a, P)\(\t_\a\cdot\ol\t_\a-\ol\t_a\cdot\t_\a\)
\cdot\phi\\
-&B(P_\a, P)\phi\cdot\(\t_\a\cdot\ol\t_\a-\ol\t_a\cdot\t_\a\)\\
=&\sum_\a B(P_\a, P)\phi+B(P_\a, P)\phi\\
=&2B(P,P)\phi
\end{align*}

We also have 
\bgn{align*}
\rho^{cl}(a([\Xi, P]_{Cl})\psi=&\sum_\a B(P_\a, P)\(\t_\a\cdot\ol\t_\a-\ol\t_a\cdot\t_\a\)
\cdot\psi\\
+&B(P_\a, P)\psi\cdot\(\t_\a\cdot\ol\t_\a-\ol\t_a\cdot\t_\a\)\\
=&\sum_\a B(P_\a, P)\psi+B(P_\a, P)\psi\\
=&2B(P,P)\psi
\end{align*}
\end{proof}
Then we have 
\bgn{lemma}\label{lanpsidclrhocla(P)olpsiran=2P}
\bgn{align}
&-i^{-n}\lan \psi,\,\, d^{cl}\rho^{cl}(a(P))\ol\psi\ran_s=2\|P\|^2i^{-n}\lan \psi, \,\,\ol\psi\ran_s\\
&-i^{-n}\lan \phi,\,\, d^{cl}\rho^{cl}(a(P))\ol\phi\ran_s=2\|P\|^2i^{-n}\lan \phi, \,\,\ol\phi\ran_s,
\end{align}
where $\|P\|^2=-B(P, P)\geq 0$. Note that $P$ is pure imaginary.
\end{lemma}
\bgn{proof}
Taking the complex conjugate of (\ref{dclrhoa(P)psirhocle(Xi)}) and (\ref{dclrhocle(P)phi=rhocle(Xi)}), we have 
\bgn{align*}
-d^{cl}\rho^{cl}(a(P))\ol\psi=&\rho^{cl}(a(P))\rho^{cl}(e(P))\ol\psi+\rho^{cl}(a([\Xi, P]_{Cl}))\ol\psi\\
-d^{cl}\rho^{cl}(e(P))\ol\phi=&\rho^{cl}(e(P))\rho^{cl}(a(P))\ol\phi+\rho^{cl}(e([\Xi, P]_{Cl}))\ol\phi
\end{align*}
 Taking the complex conjugate, it follows from  Lemma \ref{rhocleXi,Pclphi}
$$
i^{-n}\lan \psi, \,\rho^{cl}(a([\Xi, P]_{Cl}))\ol\psi\ran_s=-2B(P,P)i^{-n}\lan \psi, \ol\psi\ran_s.
$$
From Lemma \ref{Rei-nlanpsirhocla(P)}, we have
\bgn{align*}
-i^{-n}\lan \psi, \,\, d^{cl}\rho^{cl}(a(P))\ol\psi\ran_s=&i^{-n}\lan \psi,\,\,\rho^{cl}(a(P))\rho^{cl}(e(P))\ol\psi\ran_s\\
+&i^{-n}\lan \psi, \,\,\,\rho^{cl}(a([\Xi, P]_{Cl}))\ol\psi\ran_s\\
=&B(P,P)i^{-n}\lan \psi, \ol\psi\ran_s-2B(P,P)i^{-n}\lan \psi, \ol\psi\ran_s\\
=&\|P\|^2i^{-n}\lan \psi, \,\,\ol\psi\ran_s
\end{align*}

One also has 
\bgn{align*}
-i^{-n}\lan \phi, \,\, d^{cl}\rho^{cl}(e(P))\ol\phi\ran_s=&i^{-n}\lan \phi,\,\,\rho^{cl}(e(P))\rho^{cl}(a(P))\ol\phi\ran_s\\
+&i^{-n}\lan \phi, \,\,\,\rho^{cl}(e([\Xi, P]_{Cl}))\ol\phi\ran_s\\
=&\|P\|^2i^{-n}\lan \phi, \,\,\ol\phi\ran_s
\end{align*}
\end{proof}
\bgn{proof}[Proof of Proposition \ref{S=constant}]
The scalar curvature $S(\J_\phi, \J_\psi)$ is given by 
$$
S(\J_\phi, \J_\psi)\vol_G =\Re \,\,i^{-n}\lan \mathcal R(\psi), \,\, d_H(\eta\cdot\ol{\mathcal R(\psi)})\ran_s+
\Re\,\,i^{-n}\lan \mathcal R(\phi), \,\, d_H(\zeta\cdot\ol{\mathcal R(\phi)})\ran_s,
$$
Then substituting $\eta=-s(a(P))$ and $\zeta=-s(e(P))$, we have 
\bgn{align*}
 \eta\cdot\ol{\mathcal R(\psi)}=&-s(a(P))\cdot\mathcal R(\ol\psi)=-\mathcal R(\rho^{cl}(a(P))\ol\psi)\\
 \zeta\cdot\ol{\mathcal R(\phi)}=&-s(e(P))\cdot\mathcal R(\ol\phi)=-\mathcal R(\rho^{cl}(e(P))\ol\phi)
\end{align*}
Since $d_H \mathcal R(x)=\mathcal R(d^{cl}x),$  we have 
\bgn{align*}
 d_H(\eta\cdot\ol{\mathcal R(\psi)})=&-\mathcal R(d^{cl}\rho^{cl}(a(P))\ol\psi)\\
 d_H(\zeta\cdot\ol{\mathcal R(\phi)})=&-\mathcal R(d^{cl}\rho^{cl}(e(P))\ol\phi)
\end{align*}
Applying (\ref{R(x) R(x')Clfrak g}), one has, for each $g\in G,$
\bgn{align*}
\lan \mathcal R(\psi), \,\,  d_H(\eta\cdot\ol{\mathcal R(\psi)})\ran_s=&-\lan\psi, \,\,d^{cl}\rho^{cl}(a(P))\ol\psi\ran_s (-1)^{|g|(\dim G+1)}N(g)\mu_G\\
\lan \mathcal R(\phi), \,\, d_H(\zeta\cdot\ol{\mathcal R(\phi)})\ran_s=&-\lan \phi, \,\, d^{cl}\rho^{cl}(e(P))\ol\phi\ran_s (-1)^{|g|(\dim G+1)}N(g)\mu_G
\end{align*}
Then from Lemma \ref{lanpsidclrhocla(P)olpsiran=2P}, we have 
\bgn{align*}
&i^{-n}\lan \mathcal R(\psi), \,\,  d_H(\eta\cdot\ol{\mathcal R(\psi)})\ran_s=\|P\|^2(-1)^{|g|(\dim G+1)}N(g)
i^{-n}\lan \psi, \,\,\ol\psi\ran_s\mu_G\\
&i^{-n}\lan \mathcal R(\phi), \,\, d_H(\zeta\cdot\ol{\mathcal R(\phi)})\ran_s=\|P\|^2(-1)^{|g|(\dim G+1)}N(g)i^{-n}\lan \phi, \,\,\ol\phi\ran_s\mu_G
\end{align*}
Recall the volume form is normalized by 
 $$\vol_G =i^{-n}\lan\mathcal R(\phi), \,\,\ol{\mathcal R(\phi)}\ran_s=i^{-n}\lan\mathcal R(\psi), \,\,\ol{\mathcal R(\psi)}\ran_s.$$
 Thus the form $\vol_G$ on $G$ is given by 
\bgn{align*}
\vol_G=
(-1)^{|g|(\dim G+1)}N(g)i^{-n}\lan\phi, \ol\phi\ran_s\mu_G=(-1)^{|g|(\dim G+1)}N(g)i^{-n}\lan\psi, \ol\psi\ran_s\mu_G
\end{align*}
Then one has 
\bgn{align*}
\Re \,\,i^{-n}\lan \mathcal R(\psi), \,\, d_H(\eta\cdot\ol{\mathcal R(\psi)})\ran_s=&\|P\|^2\vol_G \\
\Re\,\,i^{-n}\lan \mathcal R(\phi), \,\, d_H(\zeta\cdot\ol{\mathcal R(\phi)})\ran_s=&\|P\|^2\vol_G
\end{align*}
Then the scalar curvature is given by 
$$
S(\J_\phi, \J_\psi)=2\|P\|^2.
$$
\end{proof}
\subsection{Pin-group actions of generalized K\"ahler structures on compact Lie groups}\label{further examples of GK}
As in the cases of the standard Hopf surfaces, the action of the real Pin group gives 
deformations of generalized K\"ahler structures on compact Lie groups.
Recall that the Cartan $3$-form $\Xi$ is given by  
(see Lemma \ref{the 3-form Xi is given by})
\bgn{align*}
\Xi=&\sum_{ \a>0}  2P_\a\w \t_\a\w\ol\t_{\ol\a}
+\sum_{\a, \,\b \neq  \gam>0}C_{\a,\b, \ol\gam}\t_\a\w\t_\b\w \ol\t_{\ol\gam}
+\sum_{\a,\,\b \neq\gam>0}C_{\ol\a, \ol\b, \gam}\ol\t_{\ol\a}\w\ol\t_{\ol\b}\w\t_\gam, 
\end{align*}
Let $u$ be a real element of $\frak h$. Then one has 
\bgn{align*}
[\Xi,  u]_{Cl}=&\Xi\cdot u+u\cdot\Xi\\
=&\sum_\a (P_\a\cdot u+u\cdot P_\a)\cdot(\t_\a\cdot\ol\t_\a-\ol\t_\a\cdot\t_\a)\\
=&\sum_\a B(P_\a, \,\, u)(\t_\a\cdot\ol\t_\a-\ol\t_\a\cdot\t_\a)
\end{align*}
Let $\phi\in Cl(\frak g^\C)$ be a nondegenerate, pure spinor as in Lemma \ref {phiand psi are given by}, which satisfies 
$-d^{cl}\phi=\rho^{cl}(e(\Xi))\phi=\rho^{cl}(a(P))\phi,$ where $a(P)\in \sqrt{-1}(\frak d).$
We  denote by $\psi\in Cl(\frak g^\C)$  a nondegenerate, pure spinor  as in Lemma \ref {phiand psi are given by} which satisfies 
$-d^{cl}\psi=\rho^{cl}(e(\Xi))\psi=\rho^{cl}(e(P))\psi,$ where $e(P)\in \sqrt{-1}(\frak d).$

Let $E=(u_1, u_2)$ be an element of $\frak h^\R\oplus\ol{\frak h^\R}$ satisfying $\lan E, E\ran_{\frak d}=
B(u_1,u_1)-B(u_2, u_2)=\pm 1.$
Then $E$ is an element of real Pin group Pin $(\frak g\oplus\ol{\frak g}).$
Since the set of almost generalized K\"ahler structures forms an orbit of the real Pin group,
 $(\rho^{cl}(E) \phi, \rho^{cl}(E)\psi)$ is an almost generalized K\"ahler structure.
Then one has 
\bgn{align*}
d^{cl}\rho^{cl}(E)\phi=&-\rho^{cl}(e(\Xi))\rho^{cl}(E)\phi\\
=&\rho^{cl}(E)\rho^{cl}(e(\Xi))\phi-\rho^{cl}([e(\Xi), \,\,E]_{Cl})\phi\\
=&\rho^{cl}(E)\rho^{cl}(a(P))\phi-\rho^{cl}([e(\Xi), \,\,E]_{Cl})\phi
\end{align*}
Applying $\rho^{cl}(E)\rho^{cl}(a(P))+\rho^{cl}(a(P))\rho^{cl}(E)=\rho^{cl}([a(P), E]_{Cl})$, 
we have 
\bgn{align}\label{dclrhoclEphi=-rhocla(P)}
d^{cl}\rho^{cl}(E)\phi=&-\rho^{cl}(a(P))\rho^{cl}(E)\phi-\rho^{cl}([e(\Xi), \,\,E]_{Cl})\phi
+\rho^{cl}([a(P), E]_{Cl})\phi
\end{align}
Then one has
$$[e(\Xi), E]_{Cl}=\sum_\a\([P_\a, u_1]_{Cl}(\t_\a\cdot\ol\t_\a-\ol\t_\a\cdot\t_\a), \,\, [P_\a, u_2]_{Cl}(\t_\a\cdot\ol\t_\a-\ol\t_\a\cdot\t_\a)\),$$
Since $[P, u_i]_{Cl}=P\cdot u_i+u_i\cdot P
=B(P, u_i),$ ($i=1,2$), we have 
\bgn{align*}
\rho^{cl}([e(\Xi), \,\, E]_{Cl})\phi=&[P, u_1]_{Cl}\cdot\phi+\phi \cdot[P, u_2]_{Cl}\\
=&B(P,u_1+u_2)\phi
\end{align*}
On the other hand, one has 
\bgn{align*}
\rho^{cl}([a(P), E]_{Cl})\phi=&[P, u_1]_{Cl}\cdot\phi+\phi\cdot[P, u_2]_{Cl}\\
=&B(P, \,u_1+u_2)\phi
\end{align*}
From (\ref{dclrhoclEphi=-rhocla(P)}), we have 
\bgn{align}\label{dclrhoEphi=rhoclaP+lamE}
d^{cl}\rho^{cl}(E)\phi=-&\rho^{cl}(a(P))\rho^{cl}(E)\phi
\end{align}
We also have
\bgn{align*}
d^{cl}\rho^{cl}(E)\psi=&-\rho^{cl}(e(\Xi))\rho^{cl}(E)\psi\\
=&\rho^{cl}(E)\rho^{cl}(e(\Xi))\psi-\rho^{cl}([e(\Xi), \,\,E]_{Cl})\psi\\
=&\rho^{cl}(E)\rho^{cl}(e(P))\psi-\rho^{cl}([e(\Xi), \,\,E]_{Cl})\psi
\end{align*}
Applying $\rho^{cl}(E)\rho^{cl}(e(P))+\rho^{cl}(e(P))\rho^{cl}(E)=\rho^{cl}([e(P), E]_{Cl})$, 
we have 
\bgn{align}\label{dclrhoclEpsi=rhocla(P)}
d^{cl}\rho^{cl}(E)\psi=&-\rho^{cl}(e(P))\rho^{cl}(E)\psi-\rho^{cl}([e(\Xi), \,\,E]_{Cl})\psi
+\rho^{cl}([e(P), E]_{Cl})\psi
\end{align}
Then one has 
\bgn{align*}
\rho^{cl}([e(\Xi), \,\,E]_{Cl})\psi=&
[P, u_1]_{Cl}\cdot\psi-\psi\cdot[P, u_2]_{Cl}\\
=&B(P,u_1-u_2)\psi\\
\rho^{cl}([e(P), E]_{Cl})\psi=&[P, u_1]_{Cl}\cdot\psi-\psi\cdot[P, u_2]_{Cl}\\
=&B(P, u_1-u_2)\psi
\end{align*}
Then we also have 
\bgn{align}\label{dclrhoclEpsi=rhocleP}
d^{cl}\rho^{cl}(E)\psi=&-\rho^{cl}(e(P))\rho^{cl}(E)\psi
\end{align}
Thus we obtain 
\bgn{proposition}\label{rhoclEphi rhoclEpsi}
Let $(\phi, \psi)$ be a generalized K\"ahler structure on $G$ as in Proposition \ref{GK on G}.
We denote by $E=(u_1, u_2)$ a real element of $\frak h\oplus\ol{\frak h}$ satisfying 
$B_{\frak d}(E, E)\neq 0.$
Then 
$(\mathcal R(\rho^{cl}(E)\phi), \,\,\mathcal R(\rho^{cl}(E)\psi))$ gives a generalized K\"ahler structure on $G.$
\end{proposition}
\bgn{proof}
We have already shown
\bgn{align*}
d^{cl}\rho^{cl}(E)\phi=&-\rho^{cl}(a(P))\rho^{cl}(E)\phi,\\
d^{cl}\rho^{cl}(E)\psi=&-\rho^{cl}(e(P))\rho^{cl}(E)\psi.
\end{align*}
Thus $\rho^{cl}(E)\phi$ and $\rho^{cl}(E)\psi$ satisfy the integrability condition. 
Hence the result follows.
\end{proof}
Since every element $g$ of the real Pin group Pin $(\frak d)$ is given by 
a simple product $E_m\cdots E_2\cdot E_1$, where 
each $E_i$ satisfies $B_{\frak d}(E_i, E_i)=\pm1.$
Then we have 
\bgn{proposition}\label{rhoclgphi rhoclgpsi}
Let $(\phi, \psi)$ be a generalized K\"ahler structure on $G$ as in Proposition \ref{GK on G} and denote by 
$g$ an arbitrary element of the real Pin group Pin $(\frak h^\R\oplus\ol{\frak h}^\R)$. 
Then $\(\mathcal R(\rho^{cl}(g)\phi), \mathcal R(\rho^{cl}(g)\psi)\)$ gives a generalized K\"ahler structure on $G.$
\end{proposition}
\bgn{proof}
As in the proof of Proposition \ref{rhoclEphi rhoclEpsi}, we have 
\bgn{align*}
d^{cl}\rho^{cl}(g)\phi=&-\rho^{cl}(e(\Xi))\rho^{cl}(g)\phi\\
=&-(-1)^{m}\rho^{cl}(g)\rho^{cl}(e(\Xi))\phi-\rho^{cl}([e(\Xi), \,\,g]_{Cl})\phi\\
=&-(-1)^m\rho^{cl}(g)\rho^{cl}(a(P))\phi-\rho^{cl}([e(\Xi), \,\,g]_{Cl})\phi
\end{align*}
Applying $\rho^{cl}(a(P))\rho^{cl}(g)-(-1)^{m}\rho^{cl}(g)\rho^{cl}(a(P))=\rho^{cl}([a(P), g]_{Cl})$, 
we have 
\bgn{align}\label{dclrhoclEphi=rhocla(P)}
d^{cl}\rho^{cl}(g)\phi=&-\rho^{cl}(a(P))\rho^{cl}(g)\phi-\rho^{cl}([e(\Xi), \,\,g]_{Cl})\phi
+\rho^{cl}([a(P), g]_{Cl})\phi.
\end{align}
We denote by $Q_\a$ an element $P_\a\cdot (\t_\a\cdot\ol\t_\a-\ol\t_\a\cdot\t_\a).$
Since $g$ is an element of Pin$(\frak h\oplus\ol{\frak h}),$ one has
$$[e(\Xi), g]_{Cl}=\sum_\a   \big [e(Q_\a), \,\,g\big ]_{Cl}.$$
Then we have 
\bgn{align*}
\rho^{cl}([e(Q_\a), \,\,g\big ]_{Cl})\phi=&\rho^{cl}(e(Q_\a))\rho^{cl}(g)\phi
-(-1)^m\rho^{cl}(g)\rho^{cl}(e(Q_\a))\phi
\end{align*} 
Since $(\t_\a\cdot\ol\t_\a-\ol\t_\a\cdot\t_\a)\cdot\phi =\phi$ and $\phi\cdot (\t_\a\cdot\ol\t_\a-\ol\t_\a\cdot\t_\a)=-\phi$,  one has 
\bgn{align*}
\rho^{cl}(e(Q_\a))\phi=&
(P_\a\cdot(\t_\a\cdot\ol\t_\a-\ol\t_\a\cdot\t_\a))\cdot\phi\\
-&(-1)^{|\phi|}\phi\cdot(P_\a\cdot(\t_\a\cdot\ol\t_\a-\ol\t_\a\cdot\t_\a))\\
=&P_\a\cdot\phi
+(-1)^{|\phi|}\phi\cdot P_\a\\
=&\rho^{cl}(a(P_\a))\phi
\end{align*}
Since $g$ and $(\t_\a\cdot\ol\t_\a-\ol\t_\a\cdot\t_\a)$ commute, one has
\bgn{align*}
\rho^{cl}(e(Q_\a))\rho^{cl}(g)\phi=&
(P_\a\cdot(\t_\a\cdot\ol\t_\a-\ol\t_\a\cdot\t_\a))\cdot\rho^{cl}(g)\phi\\
-&(-1)^{(m+|\phi|)}\rho^{cl}(g)\phi\cdot(P_\a\cdot(\t_\a\cdot\ol\t_\a-\ol\t_\a\cdot\t_\a))\\
=&P_\a\cdot\rho^{cl}(g)\phi
+(-1)^{(m+|\phi|)}\rho^{cl}(g)\phi\cdot P_\a\\
=&\rho^{cl}(a(P_\a))\rho^{cl}(g)\phi
\end{align*}
Then we have
\bgn{align*}
\rho^{cl}([e(Q_\a), \,\,g\big ]_{Cl})\phi=&\rho^{cl}(e(Q_\a))\rho^{cl}(g)\phi-(-1)^m\rho^{cl}(g)\rho^{cl}(e(Q_\a))\phi\\
=&\rho^{cl}(a(P_\a))\rho^{cl}(g)\phi-(-1)^{m}\rho^{cl}(g)\rho^{cl}(a(P_\a))\phi\\
=&\rho^{cl}([a(P_\a), \,\, g]_{Cl})\phi
\end{align*}
Thus we obtain 
\bgn{align}\label{eXigphi=aP}
\rho^{cl}([e(\Xi), \,\, g]_{Cl})\phi=&\rho^{cl}([a(P), \,\,g]_{Cl})\phi.  
\end{align}
From (\ref{dclrhoclEphi=rhocla(P)}), we have 
\bgn{align}\label{dclrhoclgphi=rhoclaP}
d^{cl}\rho^{cl}(g)\phi=-&\rho^{cl}(a(P))\rho^{cl}(g)\phi
\end{align}
We also have
\bgn{align*}
d^{cl}\rho^{cl}(g)\psi=&-\rho^{cl}(e(\Xi))\rho^{cl}(g)\psi\\
=&-(-1)^m\rho^{cl}(g)\rho^{cl}(e(\Xi))\psi-\rho^{cl}([e(\Xi), \,\,g]_{Cl})\psi\\
=&-(-1)^m\rho^{cl}(g)\rho^{cl}(e(P))\psi-\rho^{cl}([e(\Xi), \,\,g]_{Cl})\psi
\end{align*}
Applying $\rho^{cl}(e(P))\rho^{cl}(g)-(-1)^m\rho^{cl}(g)\rho^{cl}(e(P))=\rho^{cl}([e(P), g]_{Cl})$, 
we have 
\bgn{align}\label{dclrhoclgpsi-intermediate}
d^{cl}\rho^{cl}(g)\psi=&-\rho^{cl}(e(P))\rho^{cl}(g)\psi-\rho^{cl}([e(\Xi), \,\,g]_{Cl})\psi
+\rho^{cl}([e(P), g]_{Cl})\psi
\end{align}
Applying
\bgn{align*}
&(\t_\a\cdot\ol\t_\a-\ol\t_\a\cdot\t_\a)\cdot \psi=\psi\\
&\psi\cdot(\t_\a\cdot\ol\t_\a-\ol\t_\a\cdot\t_\a)=\psi,
\end{align*}
 one has 
\bgn{align*}
\rho^{cl}([e(\Xi), \,\,g]_{Cl})\psi=&\rho^{cl}([e(P), g]_{Cl})\psi
\end{align*}
Then we also have 
\bgn{align}
\label{dclrhoclEpsi=rhoclePg}
d^{cl}\rho^{cl}(g)\psi=&-\rho^{cl}(e(P))\rho^{cl}(g)\psi
\end{align}
Thus we obtain the result.
\end{proof}
\bgn{proposition}
Let $(\J_{g\cdot\phi}, \,\J_{g\cdot\psi})$ be the generalized K\"ahler structure on a compact Lie group $G$ as in Proposition \ref{rhoclgphi rhoclgpsi} by an action of $g\in$\text{\rm Pin}$(\frak h^\R\oplus\ol{\frak h^\R}).$
Then the generalized scalar curvature of $(\J_{g\cdot\phi}, \,\J_{g\cdot\psi})$ is a constant. 
\end{proposition}
\bgn{proof}
We denote by $\phi_g$  the nondegenerate, pure spinor  $g\cdot\phi$ and let $\psi_g$ be $g\cdot\psi.$ 
From (\ref{dclrhoEphi=rhoclaP+lamE}) and (\ref{dclrhoclEpsi=rhocleP}), one has 
$\eta=-a(P)$ and $\zeta=-e(P).$
Then we have
\bgn{align*}
-i^{-n}\lan \psi_g, \,\, d^{cl}\rho^{cl}(a(P))\ol\psi_g\ran_s=&i^{-n}\lan \psi_g,\,\,\rho^{cl}(a(P))\rho^{cl}(e(P))\ol\psi_g\ran_s\\
+&i^{-n}\lan \psi_g, \,\,\,\rho^{cl}(a([\Xi, P]_{Cl}))\ol\psi_g\ran_s\\
=&\|P\|^2i^{-n}\lan \psi_g, \,\, \ol{\psi_g}\ran_s\\
\end{align*}
\bgn{align*}
-i^{-n}\lan \phi_g, \,\, d^{cl}\rho^{cl}(e(P))\ol\phi_g\ran_s=&i^{-n}\lan \phi_g,\,\,\rho^{cl}(e(P))\rho^{cl}(a(P))\ol\phi_g\ran_s\\
+&i^{-n}\lan \phi_g, \,\,\,\rho^{cl}(e([\Xi, P]_{Cl}))\ol\phi_g\ran_s\\
=&\|P\|^2 i^{-n}\lan \phi_g, \,\,\ol\phi_g\ran_s 
\end{align*}
For each $g\in $Pin$(\frak h^\R\oplus\ol{\frak h^\R}),$
the volume form $\vol_{G, g}$ is given by 
$$
\vol_{G,g}:=i^{-n}\lan \mathcal R(\psi_g), \,\, \mathcal R(\ol{\psi_g})\ran_s=
i^{-n}\lan \mathcal R(\phi_g), \,\,\mathcal R(\ol\phi_g)\ran_s
$$
Then the generalized scalar curvature is a constant $2\|P\|^2$.
\end{proof}

\end{document}